\documentclass[]{informs3} 

\OneAndAHalfSpacedXI 

\usepackage{natbib}
 \bibpunct[, ]{(}{)}{,}{a}{}{,}%
 \def\bibfont{\small}%
\TheoremsNumberedThrough     

\EquationsNumberedThrough    

\usepackage{graphicx}
\usepackage[utf8]{inputenc}
\usepackage{enumitem}
\usepackage{multirow}
\usepackage{siunitx}
\usepackage{authblk}
\usepackage{url}
\usepackage{pgf}
\usepackage{float}
\usepackage{cite}
\usepackage{xspace}
\usepackage[font=scriptsize]{caption}
\usepackage{subcaption}
\usepackage{subfiles}
\usepackage{color}
\usepackage{booktabs}
\usepackage{algorithm}
\usepackage{algorithmic}
\usepackage{xr-hyper}
\usepackage[hidelinks]{hyperref}
\usepackage[nameinlink]{cleveref}
\usepackage{soul}
\usepackage{bbm}
\usepackage{soul}
\newtheorem{prop}{\textbf{Proposition}}
\newtheorem{properties}{\textbf{Property}}

\newtheorem{cor}{\textbf{Corollary}}

\newtheorem{rmk}{\textbf{Remark}}

\allowdisplaybreaks

\newcommand{\mynameref}[1]{\textcolor{blue}{\textit{\nameref{#1}}}}

\hypersetup{
     colorlinks = true,
     citecolor  = blue,
     linkcolor = blue,
     urlcolor = blue
}

\newcommand{\setword}[2]{%
  \phantomsection
  #1\def\@currentlabel{\unexpanded{#1}}\label{#2}%
}

\DeclareMathOperator*{\lexmin}{lex-min}

\begin{document}

\RUNAUTHOR{Guan, Basciftci, and Van Hentenryck}

\RUNTITLE{Bilevel Optimization and Heuristics for Transit Network Design with Latent Demand}

\TITLE{Bilevel Optimization and Heuristic Algorithms for Integrating Latent Demand \\ into the Design of Large-Scale Transit Systems}

\ARTICLEAUTHORS{%
\AUTHOR{Hongzhao Guan}
\AFF{H. Milton Stewart School of Industrial and Systems Engineering, Georgia Institute of Technology,\\ \EMAIL{hguan7@gatech.edu} \URL{}}
\AUTHOR{Beste Basciftci}
\AFF{Department of Business Analytics, Tippie College of Business, University of Iowa,\\ \EMAIL{beste-basciftci@uiowa.edu} \URL{}}
\AUTHOR{Pascal Van Hentenryck}
\AFF{H. Milton Stewart School of Industrial and Systems Engineering, Georgia Institute of Technology, \EMAIL{pascal.vanhentenryck@isye.gatech.edu} \URL{}}
} 

\ABSTRACT{
Capturing latent demand has a pivotal role in designing transit services as omitting these riders can lead to poor quality of service and/or additional costs. This paper explores this topic in the design of transit networks by considering the perspectives of both the transit agencies and riders. The paper presents a generic bilevel optimization model, namely the Transit Networks Design with Adoptions (TN-DA), that considers the network design decisions in the leader problem, and routing of the riders in the follower problem under the given network design, while allowing a black-box choice function for representing the adoption behavior of latent demand.
The paper then identifies structural properties of the optimal solution of the TN-DA problem, which are desirable for transit agencies for capturing adoption behavior of the riders. The paper further provides guideline metrics for the transit agencies based on these desired adoption properties. Due to the computational complexity of this bilevel problem, the paper proposes five efficient heuristic algorithms to solve large-scale instances, which leverage an iterative procedure by solving a simpler version of the TN-DA problem and integrating the evaluation of rider choices. These algorithms either satisfy the desired properties of the optimal solution or provide fast approximations. The paper presents extensive large-scale case studies on two different transit systems by utilizing real datasets: (i) On-demand Multimodal Transit Systems (ODMTS) and (ii) Scooters-Connected Transit Systems (SCTS). The results demonstrate that the heuristic algorithms can find high-quality solutions much faster than exact approaches over various instances, while satisfying key adoption properties of the optimal solutions.
}


\maketitle
%

\section{Introduction}

Transit agencies have a critical concern in designing their transportation networks: how to ensure that a transit network design will capture ridership accurately and, in particular, customer adoption. This interaction between transit agencies (or in general, network designers) and latent riders, who represent the potential riders that can adopt the transit system, can be considered as a game in the context of the transit network design problem. Essentially, this problems involves decision-making by both transit agencies and latent riders with potentially conflicting objectives. The objective of the transit agency revolves around boosting ridership by attracting latent riders while containing the rise in operational expenses. Simultaneously, latent riders expect to access high-quality services, which might incur substantial cost increases for the transit agency.

The decisions made by the transit agencies and riders are usually studied separately. More specifically, transit agencies focus on designing networks with estimated travel demand, a subject extensively explored in literature (e.g., \citet{borndorfer2007column, stiglic2018enhancing, maheo2019benders, cats2025}). However, omitting potential riders that can adopt the transit system during the network design can lead to suboptimal designs with lower rider adoption and higher costs. As an example, this phenomena is observed in designing On-Demand Multimodal Transit Systems, which is an emerging form of public transportation system that tightly integrates on-demand shuttles with fixed transit services such as high-frequency buses and urban rail systems. Specifically, the transit agency tends to have higher costs and can provide poor quality of service for certain riders if the behavior of latent riders are omitted during the network design resulting in non-effective utilization of on-demand shuttles \citep{Basciftci2020, basciftci2023capturing}. When it comes to modeling the choices of the riders, research leans towards the utilization of logit and machine learning models for demand analysis \citep{xie2003work, hagenauer2017comparative, zhao2020prediction}. However, it remains critical to integrate rider choices to the transit network design problems while taking into account the decision making processes of both the transit agency and riders. Consequently, this paper explores the possibilities of integrating the aforementioned two pieces into one single framework through a game theoretic approach. Specifically, the paper proposes a bilevel optimization framework which is called Transit Network Design with Adoptions ({\sc TN-DA}), with a rider choice modeling integrated as a black-box-like component. The optimal solution of the {\sc TN-DA} problem can be viewed as the equilibrium point between the transit agency and the latent riders. Furthermore, the {\sc TN-DA} framework is highly generalized and can be tailored to various network design problems.

Due to the bilevel nature of the problem with the combinatorial choice function of the latent riders, solving {\sc TN-DA} can be practically challenging, especially for large-scale instances. This paper thus explores whether heuristics can provide fast approximations to the optimal solutions of large-scale instances.

To this end, this paper first identifies the structural properties of the optimal solution of the TN-DA problem, in terms of the adoption behavior. These properties are desirable by the transit agencies in establishing the equilibrium between the transit agency and riders. Then, the paper presents an iterative solution framework to significantly reduce the problem complexity by splitting the original bilevel framework into two components: (i) a regular transit network design problem without rider adoption and (ii) an evaluation component that uses a choice model to capture the adoption behavior and decides the set of riders and the network to be used in the next iteration. The heuristic algorithms combine these two components by increasing incrementally the sets of riders to be considered in the design and/or by progressively adding more fixed routes in the design, while guaranteeing to satisfy desired adoption properties. Therefore, this paper then proposes five heuristic algorithms to efficiently approximate the optimal solution of the {\sc TN-DA} problem. The differences among these algorithms mostly occur at the transitions between these two primary steps, e.g, how they reshape the demand data and how they fix certain arcs in the network. 

The effectiveness of these five heuristic algorithms is evaluated across two distinct transit network design problems. In addition to the objective value that can be used for evaluating the computational performance of the algorithms, this paper further proposes new metrics that assess the quality of the network designs by measuring the trade-off between the transit agency and the potential new users based on the properties of optimal solutions.

The primary contributions of this study can be summarized as follows:
\begin{enumerate}
\item This paper introduces the {\sc TN-DA} framework, a generic bilevel optimization model that captures latent demand while designing transit network. The paper identifies adoption properties of optimal solutions, which inform heuristic design and new metrics: \textit{False Rejection Rate} and \textit{False Adoption Rate}, offering insights for evaluating and selecting transit network designs.
\item This paper proposes five heuristic algorithms to efficiently approximate optimal solutions of the {\sc TN-DA} problem that can satisfy desired adoption properties or provide rapid approximations to the optimal solution.
\item This paper demonstrates the benefits of the heuristic algorithms on the {\sc TN-DA} problem through large-scale case studies using real data on two transit systems: (i) On-demand Multimodal Transit Systems (ODMTS) and (ii) Scooter-Connected Transit Systems (SCTS). The results over various instances show that the heuristics are much more efficient than exact approaches by obtaining high-quality solutions in significantly shorter amounts of time while maintaining key adoption properties.
\end{enumerate}

The remainder of this paper is organized as follows: Section~\ref{sect:literature} first presents the related work in transit network design and its connection with travel mode choices, followed by a discussion on the application of heuristic algorithms to various mobility problems. Section~\ref{sect:original_problem} introduces the generalized TN-DA problem. The five iterative algorithms are introduced in Section~\ref{sect:algorithms}. Two detailed {\sc TN-DA} examples, namely {\sc ODMTS-DA} and {\sc SCTS-DA}, are fully demonstrated with realistic case studies in Section~\ref{sect:odmts_da}~and~\ref{sect:scts_da}. The conclusion and future research directions are discussed in Section~\ref{sect:conclusion}.
\section{Related Literature}
\label{sect:literature}

Many existing methods for addressing  transit network design problems are based on the assumption that the passenger demand is known and fixed \citep{schobel2012line}, whereas fewer number of studies have considered design/planning problems together with latent demand---potential customers who might be attracted by the proposed solutions. \citet{klier2008line} presented an optimization model without a fixed demand by employing a linear demand function of the expected travel time to model its impact on transit demand. In another study, while maximizing the number of expected passengers, the same authors employed a logit model to estimate the mode choice between the best transit route and an alternative travel mode \citep{klier2015urban}. Additionally, a series of studies presented by Canca et al. comprehensively addressed this problem. Their first study introduced a model that simultaneously determines the network design, line planning, capacity, fleet investment, and passengers mode and route choices \citep{canca2016general}. Later on, two follow-up studies proposed adaptive large neighborhood search metaheuristics to solve variations of this model over large-scale data, respectively \citep{canca2017adaptive, canca2019integrated}. Similarly, \citet{Laporte2005_Logit} aims to maximize trip coverage while considering rider behavior during the planning of rapid transit networks and \citet{Laporte2011SocioEconomic} provide a review on this area. Although these studies highlight the importance of integrating rider behavior during the network design, it remains critical to consider the decision making of both the transit agency and riders within a multi-level formulation.

A bilevel formulation is a mathematical program where certain variables are constrained by the solution of another optimization problem, as described by \citet{kleinert2021survey}. This approach is valuable for modeling hierarchical decision processes and has extensive applications in transportation. Among the studies that consider latent ridership, \citet{Basciftci2020, basciftci2023capturing} introduced rider adoption awareness to the design of On-Demand Multimodal Transit System (ODMTS). These studies indicate that a bilevel optimization framework needs to be utilized to model the problem to ensure a fair design for the transit system that provides accessible transit services to the areas with existing demand and potential demand. The optimization model includes a personalized choice model which associates adoption choices with the time and cost of trips in the ODMTS. These two studies and a follow-up study \citep{guan2024path} provide significant insights into the design of methods that solve the problem under different choice model assumptions. However, these methods still encounter computational challenges in efficiently solving large-scale problems. Building on these studies, this paper presents a generalized bilevel formulation applicable to a broad range of transit network design problems with latent demand by allowing more comprehensive choice and price models, while proposing heuristic algorithms providing desirable properties of the optimal designs with improved computational efficiency.

In general, iterative approaches are widely applied to the field of transportation science, especially for complex large-scale problems. For example, the following studies investigated iterative approaches on problems related to large-scale network design such as designing distribution networks and planning public transit for large urban areas \citep{gattermann2016iterative, schobel2017eigenmodel, cipriani2012transit, liu2021iterative}. Studies presented by \citet{siebert2013experimental} and \citet{polinder2022iterative} demonstrate a typical method for employing iterative algorithms in network design, which involves designing a network, establishing timetables, and then iteratively refining the design. When expanding the network design problem from city-wise to country-wise, another example of using iterative approaches to solve large-scale public transit problems is the study presented by \citet{liu2019iterative}. Similar to ODMTS that will be investigated in Section~\ref{sect:odmts_da}, the goal of the problem is to optimize both the operational costs and passenger travel convenience. However, instead of using time to directly measure travel convenience, a conversion from travel time to \textit{operational income} is utilized. Two case studies using the Chinese high-speed rail data presented by the authors demonstrate that iterative algorithms can efficiently solve large-scale public transit planning problems while providing profitable solutions. Furthermore, iterative approaches are utilized to solve more general network design problems such as the \textit{incomplete hub location} problems---selecting hub locations without the restriction of fully-connected networks \citep{dai2019hubbi} and the \textit{fleet assignment problem}---assigning a fleet while maximizing expected revenue \citep{dumas2009improving}. The latter study inspired \citet{kroon2015rescheduling}, and they proposed a three-step iterative approach to optimize the rolling stock to the disrupted passenger flows. Although these studies do not explicitly address latent demand in network design and omit investigating properties of optimal solutions, they highlight the efficiency of using iterative algorithms in designing large-scale logistic and urban transit networks. 

Besides designing public transit networks, iterative approaches are also applied to other topics in the field of transportation science when the original problem has a bilevel nature. \citet{yu2010genetic} first proposed a bilevel optimization model to determine the optimal bus frequencies to minimize the total passengers' travel time, and an iterative algorithm that consists of a genetic algorithm and a label-marking method is then proposed to solve this model. \citet{tawfik2021iterative} provided another example of applying iterative algorithms on bilevel problems by continuing on a previous study where a bilevel model formulates service design and pricing for freight networks \citep{tawfik2019bilevel}. \citet{Lodi2016} studied game theoretic frameworks between the transit agency and bus operators for optimizing certain contractual relationships in public transit systems, and employed iterative schemes as part of their solution procedure.

So far, not much attention has been paid to the role of iterative algorithms when integrating public transit planning and passenger mode choices. One of the earliest approaches is presented by \citet{lee2005transit}. They consider a transit network design problem with variable transit demand under a given fixed total demand that can be split into transit mode or personal vehicle mode, and this setting is also adopted by  \citet{fan2006optimal} within a genetic algorithm. The approaches used by these studies share two similarities with this paper---(i) the transit demand is variable such that riders can choose between transit services and personally-owned vehicles and (ii) at each iteration, the algorithm reshapes the transit demand then redesigns the transit network. On the other hand, this paper proposes a systematic approach for integrating rider preferences into the transit network design problem by introducing five iterative algorithms and leveraging the bilevel nature of the original problem. Three of these algorithms are driven by modifying trip sets considered in each iteration, whereas the other two are mainly driven by arc selection in the network. Furthermore, beside efficiently solving this large-scale problem, the iterative algorithms in this study also bring novelties by providing network designs with guaranteed adoption and rejection properties with insightful metrics that are practically useful to transit agencies.

\section{Transit Network Design with Adoptions Problem ({\sc TN-DA})}
\label{sect:original_problem}

This section introduces the general formulation of the Transit Network Design with Adoptions ({\sc TN-DA}) problem. The motivation for the {\sc TN-DA} problem stems from the desire to establish a transit network design $\mathbf{z}$ that concurrently serves the core riders (existing riders of the transit system) and captures latent demand with their adoption choice characteristics. Table~\ref{table:nomenclature_general} provides a summary of the nomenclatures that will be utilized in the remainder of this paper. Sections~\ref{subsect:tn_da_pre}~and~\ref{subsect:bi_level} describes the problem and introduces a generalized bilevel optimization framework for the {\sc TN-DA} problem, respectively.

\begin{table}[!t]%
\SingleSpacedXI
    \centering%
    \resizebox{\textwidth}{!}{
    \begin{tabular}{l l }%
    \toprule
    \textbf{Symbol} & \textbf{Definition} \\
    \toprule
    \textbf{General Sets}: &  \\
    $N$ & The set of all locations considered in the problem \\
    $H$ & The set of connecting points for the network, $H \subseteq N$ \\
    $T$ & The set of all trips \\
    $T_{core}$ & The set of core trips, i.e., existing and fixed transit demand. Assumed to use transit under any conditions. \\
    $T_{latent}$ & The set of latent trips, i.e., latent demand. They are the potential customers of the system. \\
    $\hat{T}$ & The set of input trips for the {\sc TN-DFD} problem. It is the union of $T_{core}$ and a subset of $T_{latent}$. \\
    $M$ & The set of possible modes, e.g., shuttles, buses, scooters, and walk. \\
    $Z$ & The set of design arcs that are considered in the design problem\\
    $\mathcal{Z}$ & The set of all feasible solutions\\
    \midrule
    \textbf{General Parameters}: & \\
    $d_{ij}^{m}$ & Travel distance between stops $i, j \in N$ under mode $m$ \\
    $t_{ij}^{m}$ & Travel time between stops $i, j \in N$ under mode $m$ \\
    $p^r$ & The number of riders in trip $r$ \\ 
    $g_\pi^r$ & Travel cost for a rider in trip $r$ under path $\pi$ \\
    $\varphi_\pi^r$ & Transit agency's revenue from a rider in trip $r$ under path $\pi$ \\
    $\zeta_\pi^r$ & Transit agency's net-cost to serve a rider in trip $r$ under path $\pi$ \\
    $t_\pi^r$ & Travel time for a rider in trip $r$ under path $\pi$ \\
    $t_{cur}^r$ & Travel time of a latent trip $r \in T_{latent}$ with its current travel mode \\
    $\alpha^r$ & Adoption factor of a latent trip $r \in T_{latent}$, used in Equation~\eqref{eq:choice_model_time_based} \\
    \midrule
    \textbf{Decision Variables}: & \\
    $z_{hl}$   & A binary variable indicating if a network arc between $h, l \in H$ is open \\
    $\delta^r$ & A binary variable indicating adoption decision of trip $r$, dynamically depending on $\mathbf{z}$ \\
    \midrule
    \textbf{Sets for Algorithms}: & \\
    $\overline{T}^k$ & The set of trips considered by the algorithms at iteration $k$ \\
    $T_{adp}^k$ & The set of latent trips that adopt network design $\mathbf{z}^k$ at iteration $k$ \\
    $C$ & The set of all adopted trips in previous iterations (Algorithms~\nameref{alg:rho_GRAD}~and~\nameref{alg:rho_GAGR}) \\
    $T_{\rho}^k$ & The set of latent trips that transferred to $C$ at iteration $k$ (Algorithms~\nameref{alg:rho_GRAD}~and~\nameref{alg:rho_GAGR})\\
    $R$ & The set of all rejected trips in previous iterations (Algorithms~\nameref{alg:eta_GRRE}~and~\nameref{alg:rho_GAGR})  \\
    $T_{m}^k$ & The set of latent trips that transferred to $R$ at iteration $k$ (Algorithms~\nameref{alg:eta_GRRE}~and~\nameref{alg:rho_GAGR})\\
    $D^k$ & The set of candidate sub-networks at iteration $k$ (Algorithms~\nameref{alg:arc_s1}~and~\nameref{alg:arc_s2}) \\
    $T_{next}^k$ & The set of latent trips for next iteration's expansion (Algorithms~\nameref{alg:arc_s1}~and~\nameref{alg:arc_s2})
    \\
    \midrule
    \textbf{Parameters for Algorithms}: & \\
    $\rho$ & The amount of trips transferred to $C$ at each iteration (Algorithms~\nameref{alg:rho_GRAD}~and~\nameref{alg:rho_GAGR})\\
    $\eta$ & The step-size applied to $m$ (Algorithms~\nameref{alg:eta_GRRE}~and~\nameref{alg:rho_GAGR})\\
    \midrule
    \textbf{Variables for Algorithms}: & \\
    $k$ & An integer that records algorithms' iterations \\
    $m$ & The amount of trips transferred to $R$ at each iteration (Algorithms~\nameref{alg:eta_GRRE}~and~\nameref{alg:rho_GAGR}) \\
    $B^k$ & The upper bound of design objective at iteration $k$ (Algorithms~\nameref{alg:arc_s1}~and~\nameref{alg:arc_s2}) \\
    $p$ & Stage indicator (Algorithms~\nameref{alg:arc_s1}~and~\nameref{alg:arc_s2}) \\
    $v^r$ & A customizable value that is used to rank trips (Algorithms~\nameref{alg:rho_GRAD},\nameref{alg:eta_GRRE}, and~\nameref{alg:rho_GAGR}) \\
    \bottomrule
    \end{tabular}
    }
    \caption{{\sc TN-DA} Notations. (For simplicity and consistency with others, notations such as $t_{\pi^r}$ are written as $t_\pi^r$.)
    }
    \label{table:nomenclature_general}
\end{table}

\subsection{Problem Description}
\label{subsect:tn_da_pre}
The {\sc TN-DA} problem revolves around a set of nodes, denoted by $N$, which represents all the locations involved in the problem. Among these nodes, a subset $H \subseteq N$ represents potential connecting points for network arcs, e.g., high-frequency buses in a transit system. A network is represented by $\mathbf{z} \in \{ 0, 1 \}^{|H| \times |H|}$, where the binary variable $z_{hl}$ represents whether an arc connecting nodes $h$ and $l$ in the set $H$ is open or not. For any nodes $i, j \in N$, $t_{ij}^{m}$ and $d_{ij}^{m}$ represent the time and distance under travel mode $m$, respectively, between these nodes, where $m \in M$. The set of modes $M$ can include various transportation modes considered in the problem, such as rail, bus, and walk.

To incorporate preferences of the riders into the planning process, the problem distinguishes between two sets of demand. The first set, denoted as $T_{core}$, represents riders who rely on transit systems and are guaranteed to use the newly designed system once it is deployed. A key assumption  is that $T_{core}$ remains fixed, meaning all individuals within this set are always transit users. The second set, denoted as $T_{latent}$, consists of riders who currently use other types of transportation and have the option to switch to the new transit system. One can consider an example where a new transit system, denoted as $\mathbf{z}$, is proposed to replace the existing system. Core riders are individuals who currently utilize the existing transit system, and they are expected to continue using public transportation when the new system is implemented. On the other hand, latent riders do not currently rely on the existing transit systems and instead choose to drive their own cars for their travel needs. These individuals have not yet adopted public transportation as their preferred mode of travel. However, the introduction of a newly-designed transit network has the potential to attract and persuade latent riders to adopt public transit. The complete set of trips, denoted by $T$, is formed by combining $T_{core}$ and $T_{latent}$. Each trip, $r \in T$, consists of an origin location, $or^r \in N$, a destination location, $de^r \in N$, and a group of $p^r$ riders who travel together, where $p^r \in \mathbb{Z}{+}$. Within the context of this paper, the term ``trip $r$" denotes a set of $p^r$ riders who undertake a journey collectively from their origin to destination. This implies that they share the same routes and possess common attributes and characteristics.

A path $\pi^r$ in the transit network for a given trip $r$ can be dynamically constructed based on the state of the network. The path should start at the trip origin node ${or}^r$, end at the trip destination node ${de}^r$, and may include zero or multiple intermediate nodes in between. The trip duration, cost, and revenue for a single rider to complete path $\pi^r$ can be calculated by considering the legs of the path. The duration, cost, and revenue values are represented as $t_\pi^r$, $g_\pi^r$, $\varphi_\pi^r$, respectively, and the net-cost $\zeta_\pi^r$ incurred by the transit agency for serving the path $\pi^r$ can be subsequently computed as $\zeta_\pi^r = g_\pi^r - \varphi_\pi^r$. It is important to note that in the given context, the terms \textit{cost}, \textit{revenue}, and \textit{net-cost} do not necessarily refer to the monetary value associated with serving trip $r$. Instead, they represent customizable measures that can be defined based on specific criteria or considerations, although the most common approach is to use monetary measures. 
\subsection{Bilevel Optimization Framework}
\label{subsect:bi_level}

Figure~\ref{fig:tn_da_formulation} presents a generalized bilevel optimization framework for capturing the travel mode adoption during the transit network design by considering the core riders ($T_{core}$) and potential riders ($T_{latent}$) that can adopt the new system. The leader problem~\eqref{eq:tn_da_problem_leader} corresponds to the network operator, which aims at designing a transit network, i.e., determining which arcs are connected to each other under adoption and problem specific considerations. The objective of the leader problem \eqref{eq:tn_da_leader_obj} minimizes the sum of (i) the investment cost of opening network arcs, which is denoted by a function $\text{\sc Invest}(\mathbf{z})$, (ii) the net-cost incurred by core riders, and (iii) the net-cost of latent riders that are adopting the new transit system. Constraint~\eqref{eq:tn_da_cst_arc} corresponds to the customized restrictions in constructing the transit network design.  Additional constraints on the arc variables $\mathbf{z}$ can be introduced based on specific problems; in this section, only the generic form is presented. Constraint~\eqref{eq:tn_da_cst_adoption} represents the mode choice function of a latent trip depending on path $\pi^r$.

\begin{figure}[ht]
\begin{subequations} 
\label{eq:tn_da_problem_leader}
\begin{alignat}{1}
\min_{\mathbf{z},  \delta^r} \quad & \text{\sc Invest}(\mathbf{z})+ \sum_{r \in T_{core}} p^r (g^r - \varphi^r) + \sum_{r \in T_{latent}} p^r \delta^r (g^r - \varphi^r) \label{eq:tn_da_leader_obj} \\
\text{s.t.} \quad & 
\mathbf{z} \in \mathcal{Z} \subseteq \{ 0, 1 \}^{|H| \times |H|}
\label{eq:tn_da_cst_arc} \\
& \delta^r = \mathcal{C}^r (\pi^r)  \quad \forall r \in T_{latent} \label{eq:tn_da_cst_adoption} \\
& \delta^r \in \{0,1\} \quad \forall r \in T_{latent} 
\end{alignat}
\end{subequations}
where $(\pi^r, g^r, \varphi^r)$ are a solution to the optimization problem $\text{\sc OPT-PATH}^r(\mathbf{z})$ as follows:
\begin{subequations} 
\label{eq:tn_da_problem_follower}
\begin{alignat}{1}
\pi^r  \in & \argmin \text{\sc OPT-PATH}^r(\mathbf{z}) \quad \forall r \in T
\label{eq:tn_dn_opt_path_cst}  \\
g^r = & \text{\sc PATH-COST}(\pi^r) \quad \forall r \in T \label{eq:tn_dn_opt_path_cost} \\
\varphi^r = & \text{\sc PATH-REVENUE}(\pi^r) \quad \forall r \in T \label{eq:tn_dn_opt_path_revenue} 
\end{alignat}
\end{subequations}
\caption{The Generalized Optimization Model for the Transit Network Design with Adoption ({\sc TN-DA}) Problem.}
\label{fig:tn_da_formulation}
\end{figure}


\paragraph{Adoption Choice Model}
In the {\sc TN-DA} problem, a choice model is employed to capture the adoption decision of latent riders, taking a path $\pi$ as its input. The binary output of $\mathcal{C}^r$ (i.e., $\delta^r$ in \eqref{eq:tn_da_leader_obj}~and~\eqref{eq:tn_da_cst_adoption}) indicates whether all riders in latent trip $r$ will adopt or reject the transit network design $\mathbf{z}$. The choice model $\mathcal{C}^r$ can be arbitrarily complex, which can range from a simple decision model to a sophisticated machine learning-based model. If available, information such as geographical and demographic data of riders can be incorporated into $\mathcal{C}^r$. The crucial aspect is that $\mathcal{C}^r$ is treated as a black-box, which implies the generalizability of the proposed framework.

A simple example of such a model is presented in Equation~\eqref{eq:choice_model_time_based}, namely {\sc Time-Based}. This particular model compares the trip duration $t_\pi^r$ with the travel duration $t_{cur}^r$ under the current travel mode of the riders, multiplied by a parameter $\alpha^r \geq 1$. Thus, if the travel duration under $\pi^r$ is less than the desired threshold of the riders, then this choice function returns 1, representing adoption of the designed system.
\begin{equation}
\label{eq:choice_model_time_based}
\text{{\sc Time-Based}}: {\cal C}^r(\pi) \equiv \mathbbm{1}(t_\pi^r \leq \alpha^r \ t^r_{cur})
\end{equation}

\paragraph{Subproblem}
Each trip $r$ in the {\sc TN-DA} problem is associated with an abstractly defined subproblem (or follower problem) called ${\text{\sc OPT-PATH}}^r$. This subproblem takes the solution of the leader problem $\mathbf{z}$ as its input and produces an optimal path $\pi^r$ specific to trip $r$. The subproblem ${\text{\sc OPT-PATH}}^r$ can be customized according to specific requirements. A common example of ${\text{\sc OPT-PATH}}^r$ is to find the shortest path from the origin ${or}^r$ to the destination ${de}^r$ under the network design $\mathbf{z}$, considering the duration as the primary criterion.  It is important to note that the optimal solution of the subproblem~\eqref{eq:tn_da_problem_follower} may not be unique. Therefore, it is necessary to apply constraints or rules to the subproblem~\eqref{eq:tn_da_problem_follower} to ensure that the optimal paths do not result in different outcomes from $\cal C^r$, such as through a lexicographic optimization approach that can account for the choice function within the objective of the subproblem. This could prevent scenarios where two paths have the same objective function value, but one is categorized as a rejecting path while the other is categorized as an adopting path by $\cal C^r$. Furthermore, the optimal path's cost $g^r$ and revenue $\varphi^r$, are computed with customized functions {\sc PATH-COST} and {\sc PATH-REVENUE}, respectively, and fed back as inputs to the leader problem, along with~$\pi^r$. To ensure the existence of $\pi^r$, a direct path between the origin and destination using an always available travel mode  can be included in the problem.  Additionally, subproblem components can be further customized with additional considerations, for instance, to ensure that travel time does not exceed a specified limit.

To this end, the intuition behind the bilevel model can be summarized as follows: The network designer first obtains the design by considering the net-cost of the trips considered. Secondly, the follower problems determine the paths for the riders during operations. The choice function in the leader problem determines, for each trip with a mode choice, whether they would adopt the proposed path or elect to use their previous travel mode. If such riders decide to use the new system, the operator receives a revenue for that trip.

An instance of the {\sc TN-DA} problem is the On-demand Multimodal Transit System Design with Adoptions (ODMTS-DA) problem, which was examined in previous studies by \citet{Basciftci2020, basciftci2023capturing, guan2024path}. After introducing the heuristic algorithms in the subsequent Section~\ref{sect:algorithms}, Section~\ref{sect:odmts_da} will revisit the ODMTS-DA problem and explore its connection to the algorithms in details. Another example of the {\sc TN-DA} problem, namely the Scooters-Connected Transit Systems Design with Adoptions ({\sc SCTS-DA}) problem, will be showcased in Section~\ref{sect:scts_da}. Both case studies concentrate on emerging technologies in the evolving landscape of public transportation.

\section{Heuristic Algorithms and Design Properties}
\label{sect:algorithms}

Section~\ref{sect:original_problem} presents the generalized {\sc TN-DA} problem within a bilevel  optimization framework. However, the presence of hierarchical decision processes in bilevel optimization problems, combined with the combinatorial nature of the {\sc TN-DA} problem with rider choices, introduce additional complexities, making it inherently difficult to solve in general. Therefore, this section presents heuristic algorithms to find high-quality solutions of the {\sc TN-DA} problem in reasonable time, especially for large-scale instances. The heuristic
algorithms use, as a key component, the Transit Network Design with Fixed
Demand ({\sc TN-DFD}) problem which is presented in Figure~\ref{fig:tn_dfd_formulation}.
Compared to the {\sc TN-DA}, the {\sc TN-DFD} removes the choice function
constraint \eqref{eq:tn_da_cst_adoption} and the decision variables of
the form $\delta^r$ that correspond to the user choices. The heuristic algorithms will use {\sc TN-DFD} with a set of trips
$\hat{T} \subseteq T$ that is always constructed from all core trips
and a subset of the latent trips.

\begin{figure}[ht]
\begin{subequations} 
\label{eq:tn_da_leader}
\begin{alignat}{1}
\min_{\mathbf{z}} \quad & \text{\sc INVEST}(\mathbf{z}) + \sum_{r \in \hat{T}} p^r (g^r - \varphi^r)  \label{eq:tn_dfd_leader_obj} \\
\text{s.t.} \quad & 
\eqref{eq:tn_da_cst_arc}\notag
\end{alignat}
\end{subequations}
where $(\pi^r, g^r, \varphi^r)$ are a solution to the optimization problem $\text{\sc OPT-PATH}^r(\mathbf{z})$ as follows:
\begin{subequations} 
\label{eq:tn_dfd_lower}
\begin{alignat}{1}
\eqref{eq:tn_dn_opt_path_cst}-\eqref{eq:tn_dn_opt_path_revenue} \quad \forall r \in \hat{T}
\end{alignat}
\end{subequations}

\caption{The Generalized Formulation for the Transit Network Design with Fixed Demand $\hat{T}$: {\sc TN-DFD}($\hat{T}$). }
\label{fig:tn_dfd_formulation}
\end{figure}

By leveraging the {\sc TN-DFD} problem, the overall idea of the heuristics can be summarized in Algorithm~\ref{alg:heuristic_general}. 
All algorithms start their first
iteration by solving the {\sc TN-DFD} problem with the core trip set $T_{core}$ to obtain an initial transit network design,
unless specified otherwise. 
Then, the set $\overline{T}^k$ denotes the set of trips considered by the algorithm at iteration $k$, i.e., $\hat{T} =
\overline{T}^k$ at iteration $k$, where the latent trips are gradually
added to $\overline{T}^k$ over the iterations depending on their adoption behavior under the resulting designs, while the trips in the core trip
set remain in $\overline{T}^k$ as $k$ increases. In certain algorithms, network design can be constructed incrementally by identifying a specific set of arcs to be fixed at each iteration. 

\begin{algorithm}[!t]
\SingleSpacedXI
\caption{TN-DA Heuristic (Generic Form)}
\label{alg:heuristic_general}
\begin{algorithmic}[1]
\STATE Set $k = 0$, $\overline{T}^0 = T_{core}$, no arcs are fixed, unless other input values are mentioned
\WHILE{Stop Criteria Not Met}
    \STATE Solve $\mathbf{z}^k = \text{\sc TN-DFD}(\overline{T}^k)$.
    \STATE Evaluate the design $\mathbf{z}^k$ and each latent trip $r \in T_{latent}$'s adoption decision $\delta^r$ using $\mathcal{C}^r$.
    \STATE Construct $\overline{T}^{k+1}$ using $T_{latent}$ and/or fix a specific set of arcs by applying certain rules while ensuring that $T_{core} \subseteq \overline{T}^{k+1}$.
    \STATE Update $k$.
\ENDWHILE
\RETURN $\mathbf{z}^{k}$ unless another design is more desirable under certain rules.
\end{algorithmic}
\end{algorithm}

Five heuristic algorithms are proposed later in this section which
can be divided into two categories: (1) trip-based iterative
algorithms that present rules for constructing the trip set $\hat{T}$
at each iteration and (2) arc-based iterative algorithms that focus on
obtaining network designs incrementally by fixing a set of arcs at
each iteration. Since the iterative algorithms only provide an
approximation to the optimal solution of the {\sc TN-DA} problem, especially on
large-scale instances, two key properties are introduced in this
section in order to assess a network design.
Ideally, one would like to obtain
the design that satisfies the following two properties. Here, the notation $\mathbf{z} = \mbox{{\sc TN-DFD}}(\hat{T})$ implies that the optimal solution of the $\mbox{{\sc TN-DFD}}(\hat{T})$ problem is the network design $\mathbf{z}$.
The set $T_{latent} \setminus \hat{T}$ represents the set of latent trips that are not considered during the design step, and the set $T_{latent} \cap \hat{T}$ represents the latent trips considered. 
\begin{properties}[Correct Rejection]
\label{property:reject}
Let $\mathbf{z} = \mbox{{\sc TN-DFD}}(\hat{T})$. Every user $r \in T_{latent} \setminus \hat{T}$ rejects the travel option in $\mathbf{z}$ obtained by $\text{\sc OPT-PATH}^r(\mathbf{z})$.
\end{properties}
\begin{properties}[Correct Adoption]
Let $\mathbf{z} = \mbox{{\sc TN-DFD}}(\hat{T})$. Every user $r \in  T_{latent} \cap \hat{T}$ adopts the travel option in $\mathbf{z}$ obtained by $\text{\sc OPT-PATH}^r(\mathbf{z})$.
\label{property:adopt}
\end{properties}

\noindent
When a design satisfies these two properties, it represents an
equilibrium between the transit agency and the latent riders. The following
proposition links the optimal solution $\mathbf{z}^*$ of {\sc TN-DA} and
Properties~\ref{property:reject}~and~\ref{property:adopt}.

\begin{prop}
\label{prop:exact_alg_prop}
Assume $\mathbf{z}^*$ is the optimal solution of the {\sc TN-DA} problem
(Figure~\ref{fig:tn_da_formulation}). There exists an input trip set $\hat{T}^*
\subseteq T$ such that $\mathbf{z}^*=\mbox{{\sc TN-DFD}}(\hat{T}^*)$ and
$\mathbf{z}^*$ satisfies both {\sc Correct Rejection} and {\sc Correct Adoption}
properties (Properties~\ref{property:reject}~and~\ref{property:adopt}).
\end{prop}
\proof{Proof:}
Given $\mathbf{z}^*$ and $\boldsymbol{\delta}^*$ as the optimum solution of the {\sc TN-DA} problem, $\hat{T}^*$ can be defined as $T_{core} \cup \{r \in T_{latent}: (\delta^{r})^* = 1\}$.
\Halmos\endproof

\noindent
The goal of the heuristic algorithms is to find the set $\hat{T}^*$,
which is computationally challenging to obtain. Thus, the proposed algorithms aim at achieving one of two different goals:
(1) efficiently finding an approximation of $\hat{T}^*$ that satisfies
the {\sc Correct Rejection} property or {\sc Correct Adoption} property, or (2)
rapidly finding an approximation to $\mathbf{z}^*$ regardless of these
properties. The rest of this section is organized as follows. Section~\ref{subsect:common_concepts} introduces key concepts for evaluating the designs found by the heuristic algorithms. Sections
\ref{subsect:trip_based_algs}~and~\ref{subsect:arc_based_algs}
describe the trip-based iterative algorithms and arc-based iterative
algorithms, respectively. The discussion of the two properties is also
presented throughout Sections~\ref{subsect:trip_based_algs}~and~\ref{subsect:arc_based_algs} as the iterative algorithms are introduced.

\subsection{Key Concepts for Design Evaluation}
\label{subsect:common_concepts}
This section summarizes multiple common features shared
by all proposed iterative algorithms.  Since the {\sc TN-DFD} is employed to approximate the optimal solution of the {\sc TN-DA}, the objective function~\eqref{eq:tn_da_leader_obj} of the {\sc TN-DA} is
utilized to evaluate any network design $\mathbf{z}$ produced by the
{\sc TN-DFD}, using Equation~\eqref{eq:design_obj}. Regardless of the set $\hat{T}$ used to obtain the design $\mathbf{z}$, the input trip set to evaluate $\mathbf{z}$ is always the full trip
set $T$. Under the design $\mathbf{z}$, this evaluation is performed using the choice decisions $\delta^r$ from $\mathcal{C}^r$ and the optimal path $\pi^r$ associated with each trip $r$.
\begin{equation}
  \text{\sc Eval}(\mathbf{z}) =  \text{\sc Invest}(\mathbf{z}) + \sum_{r \in T_{core}} p^r (g_\pi^r - \varphi_\pi^r) + \sum_{r \in T_{latent}} p^r \delta^r (g_\pi^r - \varphi_\pi^r) \label{eq:design_obj}
\end{equation}

To evaluate the quality of the network designs found by the iterative
algorithms, two additional metrics, i.e., {\sc False Rejection Rate} ($\%R_{false}$) and {\sc False Adoption Rate}
($\%A_{false}$), are introduced to measure the violation of the
{\sc Correct Rejection} and {\sc Correct Adoption} properties (Properties~\ref{property:reject}~and~\ref{property:adopt}). These two metrics are
defined as follows:
\begin{equation*}
    \%R_{false} = \cfrac{| \{r \in T_{latent} \setminus \hat{T}: r \text{ adopts design } \mathbf{z}\}|}{|T_{latent}| } \times \%100
\end{equation*}

\begin{equation*}
\%A_{false} = \cfrac{| \{r \in  T_{latent} \cap \hat{T}: r \text{ rejects design } \mathbf{z}\}|}{|T_{latent}|} \times \%100
\end{equation*}

\noindent
Here, $\%R_{false}$ indicates the percentage of latent trips that are excluded during the network design although they adopt the resulting design, implying false rejection, whereas $\%A_{false}$ corresponds to the percentage of latent trips that are included during the network design although they reject the resulting design, implying false adoption. Core trips are excluded from the calculation of these metrics, because they are always included in $\hat{T}$ and not associated with a choice function. 
The subsequent two Remarks highlight that obtaining zero {\sc False Rejection Rate} or {\sc False Adoption Rate} is neither a necessary nor a sufficient condition for finding an optimal solution of {\sc TN-DA} as these conditions depend on the trip set considered in the {\sc TN-DFD} problem. Nevertheless, they establish desired adoption properties by identifying an equilibrium between the transit agency and latent riders based on the latent riders considered during the network design.

\begin{rmk}
\label{rmk:0p_but_not_opt}
It is possible that the solution of the {\sc TN-DFD} with a trip set
$\hat{T} \neq \hat{T}^*$ is a design $\mathbf{z} \neq \mathbf{z^*}$
that satisfies both {\sc Correct Rejection} property and {\sc Correct Adoption}
property. Thus, satisfying both properties under an input trip set
$\hat{T}$ does not necessarily imply that $\mathbf{z}$ is the optimal
solution of {\sc TN-DA}.
\end{rmk}

\begin{rmk}
\label{rmk:opt_but_not_0p}
It is possible that the optimal solution of the {\sc TN-DFD} with a trip
set $\hat{T} \neq \hat{T}^*$ is the optimal design $\mathbf{z}^*$ of
the {\sc TN-DA}.  In other words, $\mathbf{z}^*$ may be discovered by a
trip set $\hat{T} \neq \hat{T}^*$ with $\%R_{false} \neq 0$ or
$\%A_{false} \neq 0$ under $\hat{T}$.
\end{rmk}

\subsection{Trip-based Iterative Algorithms}
\label{subsect:trip_based_algs}

This section introduces the first category of iterative
algorithms---the \textit{trip-based algorithms.} In general, each
iteration of this type of algorithm can be summarized by three
fundamental steps. The first step solves the
{\sc TN-DFD} problem with $\hat{T} = \overline{T}^k$. The second step
evaluates the design with $\text{\sc Eval}(\mathbf{z})$ (Equation \eqref{eq:design_obj}).
The third step computes the new $\overline{T}^{k+1}$ for the next iteration.
All trip-based iterative algorithms have identical first and
second steps, with the only difference occurring in the third step. At the third step, the trip-based algorithms need to compute a customizable value for each trip $r$, denoted as $v^r$ (such as the net-cost $\zeta_\pi^r = g_\pi^r - \varphi_\pi^r$), based on the optimal path of that trip under the given network design. This value is then utilized in constructing the trip set to be considered in the subsequent iterations.

Three trip-based iterative algorithms are proposed. The first two are
presented in Sections~\ref{subsubsect:rho_GRAD} and
\ref{subsubsect:m_GRRE}. Another one combining the previous two
algorithms is explained in Section~\ref{subsubsect:rho_GAGR}, which
enables to explore more intermediate network designs, despite of its
much longer expected running time. Section
\ref{subsubsect:TripBasedAlgorithms_Properties} presents some
properties of the proposed trip-based algorithms.

\subsubsection{Greedy Adoption Algorithm (\mynameref{alg:rho_GRAD})}
\label{subsubsect:rho_GRAD}
The first trip-based iterative algorithm starts with the core trips
and greedily adds adopting trips at each iteration. Algorithm
\mynameref{alg:rho_GRAD} maintains a set $C$ of trips and transfers
latent trips from $T_{latent}$ to $C$ at each iteration. Specifically, at
a given iteration $k$, the solution $\mathbf{z}^k$ identifies a set of
adopting trips $T_{adp}^k \in T_{latent}$. It transfers $\rho$ trips from
$T_{adp}^k$ to $C$, choosing those with the least $v^r$ values. This set of
$\rho$ trips is denoted by $T_{\rho}^k$. The algorithm
terminates when there are no latent trips to be absorbed by set $C$
at the end of an iteration, i.e., $|T_{\rho}^k| = 0$. Algorithm
\mynameref{alg:rho_GRAD} satisfies the following property.

\begin{prop}
\label{prop:alg_1_prop}
The last design $\mathbf{z}^k$ found by Algorithm \mynameref{alg:rho_GRAD} with trip set $\hat{T} = \overline{T}^k$ satisfies the {\sc Correct Rejection} property (Property~\ref{property:reject}), i.e., $\% R_{false} = 0$.
\end{prop}

\proof{Proof:}
The last design $\mathbf{z}^k$ is found by solving {\sc TN-DFD} with $\hat{T} = T_{core} \cup C$. On completion, $|T_{\rho}^k| = 0$ and no latent trip in $T_{latent} \setminus C$  adopts design $\mathbf{z}^k$. The result follows.
\Halmos\endproof

\begin{algorithm}[!t]
\SingleSpacedXI
\caption{$\rho$-GRAD}
\label{alg:rho_GRAD}
\begin{algorithmic}[1]
\STATE Set $k = 0$, $\overline{T}^0 = T_{core}$, $C = \emptyset$, unless other input values are specified.
\STATE \textbf{Require}: step-size parameter $\rho$
\WHILE{True}
\STATE Set $\mathbf{z}^k = \mbox{{\sc TN-DFD}}(\overline{T}^k)$.
\STATE Set $T_{check} = T_{latent} \setminus C$.
\FORALL{$r \in T_{check}$}
\STATE Set $\pi^r$ = $\argmin \text{\sc OPT-PATH}^r(\mathbf{z}^k)$.
\STATE Set $\delta^r = {\cal C}^r(\pi^r)$.
\IF{$\delta^r = 1$}
\STATE Compute $v^r$ for trip $r$.
\ENDIF
\ENDFOR
\STATE Set $T_{adp}^k = \{ r \in T_{check} \mid \delta^r = 1 \}$.
\STATE Set $T_{\rho}^k$ as the set of $\rho$ trips from $T_{adp}^k$ with the lowest $v^r$.
\IF{$T_{\rho}^k = \emptyset$}
\STATE break.
\ELSE
\STATE Set $C = C \cup T_{\rho}^k$.
\STATE Set $\overline{T}^{k+1} = T_{core} \cup C$.
\STATE Set $k = k + 1$.
\ENDIF
\ENDWHILE
\RETURN $\mathbf{z}^{k}$
\end{algorithmic}
\end{algorithm}


\subsubsection{Greedy Rejection Algorithm (\mynameref{alg:eta_GRRE})}
\label{subsubsect:m_GRRE}

Algorithm~\mynameref{alg:eta_GRRE} greedily excludes trips from $T_{latent}$ and
introduces a growing set $R$ of rejecting trips. In addition, a
variable $m$ is used to indicate how many adopting trips from $T_{latent}
\setminus R$ are included in the set $\overline{T}^k$ and $m$ is
increased by $\eta$ at each iteration. Using $T_{latent} \setminus R$ directly
would lead to many rejections, so the algorithm proceeds carefully
when choosing $\overline{T}^k$ early on. Algorithm
\mynameref{alg:eta_GRRE} terminates when it meets the following two
conditions at an iteration $k \geq 2$:
\begin{enumerate}
    \item the design is stable, i.e., $\mathbf{z}^{k} = \mathbf{z}^{k-1}$;
    \item  $m - \eta \geq |T_{adp}^k|$.
\end{enumerate}
Note that condition 1 alone is not adequate. Indeed, in early stages,
the algorithm may tend to find the same design because $m$ is still
relatively small. As a result, condition 2 is firstly employed to guarantee that the algorithms delve into a greater number of iterations. Nevertheless, it is crucial to acknowledge that beyond the point where $m$ exceeds $|T_{adp}^k| + \eta$, the algorithm will commence producing designs similar to those from earlier iterations. This is likely due to the decrease in $\hat{T}$ over iterations when many trips are included in $R$, which might also lead to oscillating results between iterations. Thus, condition 2 serves the purpose of striking a balance between these two effects. The designs found by Algorithm
\mynameref{alg:eta_GRRE} may not satisfy either the {\sc Correct Rejection}
property or the {\sc Correct Adoption} property. Instead, this algorithm is
particularly designed to provide a rapid approximation.

\begin{rmk}
Because of possible oscillations in network designs between iterations, the \mynameref{alg:eta_GRRE} algorithm reports the design with the minimal objective instead of the last design. These oscillations occur as the set $\overline{T}^{k}$ does  not monotonically increase. Note that if the last design $\mathbf{z}^{k}$ is reported, then the algorithm satisfies the {\sc Correct Adoption} property (Property~\ref{property:adopt}) given that $\mathbf{z}^{k} = \mathbf{z}^{k-1}$.
\end{rmk}
\begin{algorithm}[!t]
\SingleSpacedXI
\caption{$\eta$-GRRE}
\label{alg:eta_GRRE}
\begin{algorithmic}[1]
\STATE Set $k = 0$, $m = 0$, $R = \emptyset$, $\overline{T}^0 = T_{core}$, and $o^* = \infty$, unless other input values are specified.
\STATE \textbf{Require}: step-size parameter $\eta$
\WHILE{True}
\STATE $\mathbf{z}^k = \mbox{{\sc TN-DFD}}(\overline{T}^k)$.
\IF{$\text{\sc Eval}(\mathbf{z}^k) < o^*$}
\STATE Set $o^* = \text{\sc Eval}(\mathbf{z}^k)$.
\STATE Set $\mathbf{z}_{min} = \mathbf{z}^k$.
\ENDIF
\STATE Set $T_{check} = T_{latent} \setminus R$.
\FORALL{$r \in T_{check}$}
\STATE Set $\pi^r$ = $\argmin \text{\sc OPT-PATH}^r(\mathbf{z}^k)$ 
\STATE Set $\delta^r = {\cal C}^r(\pi^r)$.
\IF{$\delta^r = 1$}
\STATE Compute $v^r$ for trip $r$.
\ELSE
\STATE Set $R = R \cup \{r\}$.
\ENDIF
\ENDFOR
\STATE Set $T_{adp}^k = \{ r \in T_{check} \mid \delta^r = 1 \}$.
\STATE Set $m = m + \eta$.
\STATE Set $T_m^k = $ the set of $m$ trips from $T_{adp}^k$with the lowest $v^r$.
\IF{the termination criteria is met}
\STATE break.
\ELSE
\STATE Set $\overline{T}^{k+1} = T_{core} \cup T_m^k$;
\STATE Set $k = k + 1$.
\ENDIF
\ENDWHILE
\RETURN $\mathbf{z}_{min}$
\end{algorithmic}
\end{algorithm}

\subsubsection{Greedy Adoption Algorithm with Greedy Rejection Subproblem (\mynameref{alg:rho_GAGR})}
\label{subsubsect:rho_GAGR}

Algorithm~\mynameref{alg:rho_GAGR} is the same as Algorithm~\mynameref{alg:rho_GRAD}, except that it uses Algorithm~\mynameref{alg:eta_GRRE} as its subroutine instead of solving the {\sc TN-DFD} problem to obtain a
network design. This design decision makes it possible for Algorithm~\mynameref{alg:rho_GAGR} to explore many more designs than Algorithm~\mynameref{alg:rho_GRAD}. Algorithm~\mynameref{alg:rho_GAGR} shares the same
terminating criteria with Algorithm~\mynameref{alg:rho_GRAD}. The full algorithm is presented in Appendix~\ref{sect:heuristic_algorithm_appendix}.

\subsubsection{Additional Properties of Trip-based Algorithms}
\label{subsubsect:TripBasedAlgorithms_Properties}

This section analyzes the change in the net-cost of the trips, as different sets of trips are considered during the network design. To this end, recall that $\zeta^r(\mathbf{z}) = g^r(\mathbf{z}) - \varphi^r(\mathbf{z})$ is defined as the optimum net-cost of the trip $r$ under the network design~$\mathbf{z}$. 

\begin{prop}
\label{prop:additional_prop_trip}
Let $\mathbf{z}^1 = \mbox{{\sc TN-DFD}}(T^1)$ and $\mathbf{z}^2=\mbox{{\sc TN-DFD}}(T^2)$ with $T^1 \subseteq T^2$. Then $\sum\limits_{r \in T^2 \setminus T^1} p^r (\zeta^r(\mathbf{z}^2) - \zeta^r(\mathbf{z}^1)) \leq 0$.
\end{prop}
\proof{Proof:}
Note that the optimum objective function value of the network design $\mathbf{z}^2$ under the set of trips $T^2$ can be written as $\sum\limits_{h,l \in H} \beta_{hl} \mathbf{z}^2_{hl} + \sum\limits_{r \in T^2} p^r \zeta^r(\mathbf{z}^2)$. Since $\mathbf{z}^1$ is a feasible network design for this  problem,  the following relationship holds:
\begin{align*}
\sum_{h,l \in H} \beta_{hl} \mathbf{z}^2_{hl} + \sum_{r \in T^2} p^r \zeta^r(\mathbf{z}^2) \leq \sum_{h,l \in H} \beta_{hl} \mathbf{z}^1_{hl} + & \sum_{r \in T^2} p^r \zeta^r(\mathbf{z}^1) \\
\sum_{h,l \in H} \beta_{hl} \mathbf{z}^2_{hl} + \sum_{r \in T^1} p^r \zeta^r(\mathbf{z}^2) + \sum_{r \in T^2 \setminus T^1} p^r \zeta^r(\mathbf{z}^2) \leq \sum_{h,l \in H} \beta_{hl} \mathbf{z}^1_{hl} + & \sum_{r \in T^1} p^r \zeta^r(\mathbf{z}^1) + \sum_{r \in T^2 \setminus T^1} p^r \zeta^r(\mathbf{z}^1) \\
\sum_{h,l \in H} \beta_{hl} (\mathbf{z}^2_{hl} - \mathbf{z}^1_{hl}) + \sum_{r \in T^1} p^r (\zeta^r(\mathbf{z}^2) - \zeta^r(\mathbf{z}^1)) + & \sum_{r \in T^2 \setminus T^1} p^r (\zeta^r(\mathbf{z}^2) - \zeta^r(\mathbf{z}^1))  \leq 0
\end{align*}
\noindent
Following an analogous relationship by comparing the network designs under the set of trips $T^1$ results in
\begin{align*}
\sum_{h,l \in H} \beta_{hl} \mathbf{z}^2_{hl} + \sum_{r \in T^1} p^r \zeta^r(\mathbf{z}^2) & \geq \sum_{h,l \in H} \beta_{hl} \mathbf{z}^1_{hl} + \sum_{r \in T^1} p^r \zeta^r(\mathbf{z}^1) \\
\sum_{h,l \in H} \beta_{hl} (\mathbf{z}^2_{hl} - \mathbf{z}^1_{hl}) + \sum_{r \in T^1} p^r (\zeta^r(\mathbf{z}^2) - \zeta^r(\mathbf{z}^1)) & \geq 0
\end{align*}
\noindent
Combining the above derivations, it follows that $\sum\limits_{r \in T^2 \setminus T^1} p^r (\zeta^r(\mathbf{z}^2) - \zeta^r(\mathbf{z}^1)) \leq 0$.
\Halmos\endproof

\begin{cor}
If the set $T^2 \setminus T^1$ is a singleton $\{r\}$, then $\zeta^r(\mathbf{z}^2) \leq \zeta^r(\mathbf{z}^1)$.
\end{cor}

\noindent
These results have interesting consequences for Algorithms
\mynameref{alg:rho_GRAD} and \mynameref{alg:rho_GAGR}. With a step size of 1, the added trip will not have a net-cost in the last design worse than those in the prior iteration. With larger step size, at least one of the added trips will satisfy this property. However, this
property is not guaranteed for all the added trips in that
iteration. Note that, as opposed to Algorithms
\mynameref{alg:rho_GRAD} and \mynameref{alg:rho_GAGR}, this property does not necessarily hold for Algorithm \mynameref{alg:eta_GRRE} as $\overline{T}^{k} \subseteq \overline{T}^{k+1}$ is not ensured in any iteration $k$.

\subsection{Arc-based Iterative Algorithms}
\label{subsect:arc_based_algs}

This section introduces the second type of iterative algorithms---the \textit{arc-based iterative algorithms}. 
These algorithms begin with a pre-defined transit network design, denoted as $\mathbf{z}_{fixed}$, which can be either an empty network or an existing network such as a well-developed rail system. 
At the beginning of each iteration, these algorithms solve the {\sc TN-DFD} problem and obtain a network design building on the open arcs in $\mathbf{z}_{fixed}$.   
Then, a sub-network from the resulting design is incrementally added in each iteration to enhance the overall network. 
The input trip set considered for the
{\sc TN-DFD} problem is then expanded, considering the adoption decisions of the riders. These steps are iterated until no improving sub-network can be added to the design. 
Since these algorithms select the best sub-network with respect to the improvement in the overall objective, arc-based algorithms produce a sequence
of network designs whose evaluated objectives are monotonically decreasing.  As
a result, they find primal solutions quickly, which is not guaranteed
in trip-based algorithms. In summary, arc-based algorithms satisfy the
following properties:
\begin{enumerate}
\item The design objective, computed by function $\text{\sc Eval}(\mathbf{z})$ (Equation~\eqref{eq:design_obj}),  monotonically decreases.
\item The transit network design $\mathbf{z}_{fixed}$ monotonically increases.
\end{enumerate}

This paper proposes two arc-based algorithms. Section
\ref{subsubsect:arc_base_greedy_alg} presents the Arc-based Greedy
Algorithm \mynameref{alg:arc_s1}, and Section~\ref{subsubsect:arc_base_extended_greedy_alg}
presents the algorithm \mynameref{alg:arc_s2} that extends
\mynameref{alg:arc_s1} to overcome the potential challenges faced by
\mynameref{alg:arc_s1} in terms of run time and solution quality. Both algorithms are parametrized by one or
more functions that expand the set of users at each iteration.
Section~\ref{subsubsect:arc_base_trip_rule} introduces
expansion rules to construct these {\sc Expand} functions used by the algorithms.

\subsubsection{Arc-based Greedy Algorithm (\mynameref{alg:arc_s1})}
\label{subsubsect:arc_base_greedy_alg}

\begin{algorithm}[!t]
\SingleSpacedXI
\caption{arc-S1}
\label{alg:arc_s1}
\begin{algorithmic}[1]
\STATE Set $k = 0$, $\mathbf{z}_{fixed} = \Vec{0}$, $\overline{T}^0 = T_{core}$, $B^0 = \infty$.
\STATE \textbf{Require}: {\sc Expand} function for increasing the set of users considered.
\WHILE{True}
\STATE Set $\mathbf{z}_{temp}^k$ = {\sc TN-DFD}($\overline{T}^k$,$\{ (h,l) \mid {{z}_{fixed}}_{hl} = 1 \ \> \  h, l \in H\}$).
\STATE Set $\mathbf{z}^k_{unfixed} = \mathbf{z}_{temp}^k - \mathbf{z}_{fixed}$.
\STATE Set $D^k$ to the set of all candidate sub-networks in $\mathbf{z}^k_{unfixed}$.
\IF{$D^k = \emptyset$}
\STATE break.
\ENDIF
\FORALL{$\mathbf{d} \in D^k$}
\IF{$(\mathbf{z}_{fixed} + \mathbf{d}) \in \mathcal{Z}$ (i.e. if the design is feasible)}
\STATE Set $obj_{\mathbf{d}} = \text{\sc Eval}(\mathbf{z}_{fixed} + \mathbf{d})$.
\ENDIF
\ENDFOR
\STATE Set ${obj}^k = \min\limits_{\mathbf{d} \in D^k} \{ obj_{\mathbf{d}} \}$ and $\mathbf{d}^k = \argmin\limits_{\mathbf{d} \in D^k} \{ obj_{\mathbf{d}} \}$.
\IF {${obj}^k < B^{k}$}
\STATE Set $B^{k+1} = {obj}^k$.
\ELSE
\STATE break.
\ENDIF
\STATE Set $\mathbf{z}_{fixed}$ = $\mathbf{z}_{fixed} + \mathbf{d}^k$.
\STATE Set $\mathbf{z}^k$ = $\mathbf{z}_{fixed}$
\STATE Set $T_{next}^k$ = {\sc Expand}($T_{latent},\mathbf{z}^{k}$).
\STATE Set $\overline{T}^{k+1} = \overline{T}^{k} \cup T_{next}^k$.
\STATE Set $k = k + 1$.
\ENDWHILE
\RETURN $\mathbf{z}^{k}$
\end{algorithmic}
\end{algorithm}

Algorithm~\mynameref{alg:arc_s1} describes the arc-based greedy
algorithm. It maintains an upper bound $B^k$ that represents the best design objective found up to iteration $k$. Algorithm~\mynameref{alg:arc_s1} also uses a
generalization of the {\sc TN-DFD}, where the second parameter represents
the set of arcs that must be fixed to 1 in the {\sc TN-DFD} solution
(line 4).

During each iteration, after solving such a generalized {\sc TN-DFD}, the
resulting design $\mathbf{z}_{temp}^k$ is decomposed into two parts:
$\mathbf{z}_{fixed}$ and $\mathbf{z}_{unfixed}^k$ (line 5). The
partial design $\mathbf{z}_{fixed}$ denotes all the arcs which were
constrained to be opened in the {\sc TN-DFD}, while
$\mathbf{z}^k_{unfixed}$ is the remaining assignment in
$\mathbf{z}_{temp}^k$. 
Then, the sub-network candidate pool $D^k$ needs to be carefully determined while ensuring that the addition of sub-networks does not lead to infeasibility in future iterations by considering
the specific characteristics and requirements of the problem at hand. In particular, each iteration attempts to add a new sub-network
$\mathbf{d}^k \in D^k$ to the fixed design $\mathbf{z}_{fixed}$. Algorithm \mynameref{alg:arc_s1} does this greedily, adding the sub-network
$\mathbf{d}^k$ which provides the greatest decrease in $B^k$ (lines
10--18).  Therefore, the series $B^0, \ldots, B^k$ decreases
monotonically over the iterations, while the design $\mathbf{z}_{fixed}$ expands (line 21), where the computational time of each
iteration tends to decrease as the number of variables decreases. Algorithm
\mynameref{alg:arc_s1} also expands the trips to consider in the
{\sc TN-DFD} (lines 23--24).
The algorithm terminates when the design objective bound is not
decreasing, i.e., $B^k = B^{k+1}$, which can happen under the following two cases:
\begin{enumerate}
    \item $D^k$ is empty, i.e., it does not contain any new candidate sub-networks.
    \item No sub-networks in $D^k$ provides a better objective value.
\end{enumerate}

\subsubsection{Arc-based Extended Greedy Algorithm (\mynameref{alg:arc_s2})}
\label{subsubsect:arc_base_extended_greedy_alg}

The {\sc Expand} function used in Algorithm~\ref{alg:arc_s1} can have a significant impact on the performance of the algorithm, as it provides the balance in avoiding the set of trips in $\overline{T}^k$ from growing too quickly or too slowly over the iterations. A rapid expansion of the set $\overline{T}^k$ can introduce computational challenges when solving the {\sc TN-DFD} problem.  Conversely, a slow expansion of $\overline{T}^k$ may lead to premature termination of the algorithm, resulting in poor solutions. To overcome this difficulty, \mynameref{alg:arc_s2} implements a two-phase extension of \mynameref{alg:arc_s1} by using two
{\sc Expand} functions. Initially, \mynameref{alg:arc_s2}
expands the set of trips slowly until convergence using
function {\sc Expand$_1$}. Then \mynameref{alg:arc_s2} enters
its second phase where the trip set increases faster. Both phases
implement \mynameref{alg:arc_s1}, except that the second phase starts
with the fixed design and the set of trips produced by phase 1. The full algorithm is presented in Appendix~\ref{sect:heuristic_algorithm_appendix}.

\subsubsection{Trip Expansion Rules}
\label{subsubsect:arc_base_trip_rule}

The arc-based algorithms have been examined with various expansion rules to build the {\sc Expand} functions over different transit systems. The initial rule, referred to as the ``Full Adoption Expansion Rule", considers \textit{all latent trips that adopt $\mathbf{z}^k$} to be included in the subsequent iteration. This rule is intuitive and can be generally applied to various types of transit systems. In the following Sections~\ref{sect:odmts_da}~and~\ref{sect:scts_da}, the paper will further demonstrate how expansion rules can be carefully designed within a specific context by leveraging the underlying problem. This will showcase the flexibility of the approach in adapting to different transit system scenarios while maintaining the underlying principles discussed earlier.  When using this expansion rule, Algorithm \mynameref{alg:arc_s1} satisfies the {\sc Correct Rejection} property. This result extends to Algorithm \mynameref{alg:arc_s2}.
\begin{prop}
\label{prop:arc_s1_proFR}
The transit network design returned by Algorithm \mynameref{alg:arc_s1} with Full Adoption Expansion Rule satisfies the {\sc Correct
Rejection} property (Property~\ref{property:reject}), i.e., $\% R_{false} = 0$.
\end{prop}

\proof{Proof:}
When Algorithm \mynameref{alg:arc_s1} terminates at iteration $k$,
the design $\mathbf{z}_{fixed}$ does not change. That means that the
{\sc TN-DFD} was executed with $\hat{T}^k = \overline{T}^{k-1} \cup T_{next}^{k-1}$.
Moreover, $T_{next}^{k-1}$ is the set of all latent trips that adopt
$\mathbf{z}_{fixed}$. Hence, all latent trips that are not
in $\hat{T}^k$ reject design $\mathbf{z}^{k}$ and Property~\ref{property:reject} holds.
\Halmos\endproof

\begin{cor}
\label{cor:alg_arc_s2_second_rule_a}
The transit network design returned by Algorithm \mynameref{alg:arc_s2} with
the Full Adoption Expansion Rule in the second phase satisfies the {\sc Correct
Rejection} property (Property~\ref{property:reject}).
\end{cor}
\section{TN-DA Example 1: ODMTS Design with Adoptions Problem}
\label{sect:odmts_da}

This section presents a detailed example of the {\sc TN-DA} problem with three realistic case studies, focusing on a specific type of transit network known as On-demand Multimodal Transit Systems (ODMTS). Therefore, within this particular setting, the {\sc TN-DA} problem is referred to as {\sc {\sc ODMTS-DA}} (ODMTS Design with Adoptions) Problem.

\begin{figure}[!ht]
    \centering
    \includegraphics[width=0.4\textwidth]{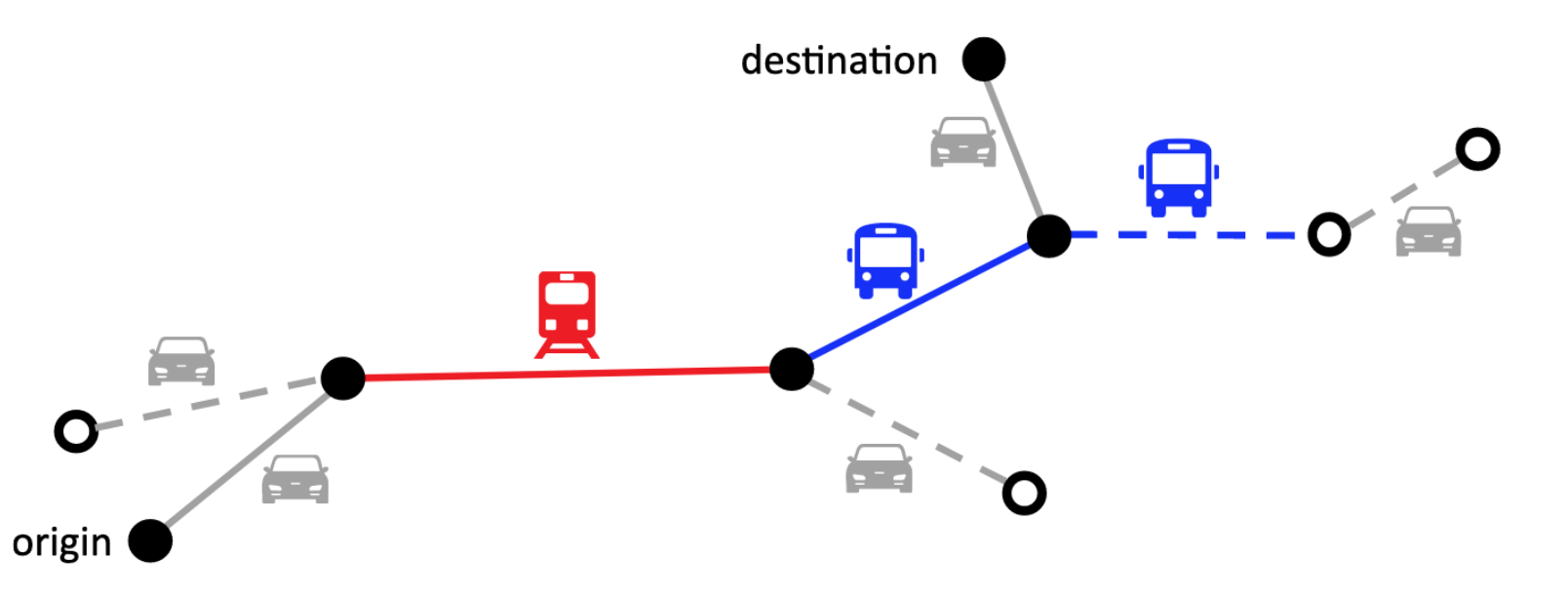}
    \caption{Picture that illustrates the On-demand Multimodal Transit Systems (ODMTS). The solid circles represent transit hubs. The solid line represents a typical ODMTS path offered by the transit agency.}
    \label{fig:odmts_general}
\end{figure}

As depicted in Figure~\ref{fig:odmts_general}, ODMTS are operated around a number of hubs that serve as stops for buses and rail. The on-demand shuttles are primarily employed as the feeders to the hubs
and represent an effective solution to the first-and-last-mile
problem faced by most of the transit agencies. Several simulation
studies have demonstrated that ODMTS can provide considerable cost and
convenience benefits in different cities under diverse circumstances
\citep{maheo2019benders, dalmeijer2020transfer,
auad2021resiliency}. In practice, ODMTS were also successfully deployed in the cities of Atlanta \citep{van2023marta} and Savannah \citep{cat2024} in 2022 and 2024---2025, respectively. Using the {\sc ODMTS-DA} as an example of the {\sc TN-DA}, the goal of this section is to demonstrate the practical advantages and outcomes of the heuristic algorithms in a real-world transit network design problem.
  
In the context of {\sc ODMTS-DA}, two
previous studies addressed this problem by modeling it with
bilevel optimization frameworks
\citep{Basciftci2020,basciftci2023capturing}. In particular, the agency proposes an ODMTS service and optimal paths for each user under the current design by considering both their duration and cost; the potential riders then evaluate the convenience of these paths, and decide whether to adopt the ODMTS service. It is important to note that this is a fundamental distinction between ODMTS and conventional transit systems, as the transit agencies provide the designated paths for riders by using on-demand shuttles and transit network, and potential riders make their adoption decisions accordingly.

The remainder of this section is structured as follows: Sections~\ref{subsect:odmts_da_bilevel}---\ref{subsect:odmts_da_heuristics} study {\sc ODMTS-DA} problem using the bilevel optimization framework and explain the employment of heuristic algorithms to solve {\sc ODMTS-DA}. Section \ref{sec:CaseStudyOverview} provides an overview of the three real-world case studies, which are categorized by size as medium, large, and extra-large. Finally, Sections \ref{sec:MediumInstancepsilanti}---\ref{subsect:atlanta_extra_large_instance} present the results over the corresponding three cases.

\subsection{ODMTS-DA Bilevel Formulation}
\label{subsect:odmts_da_bilevel}

The analogy between {\sc ODMTS-DA} and the generalized framework {\sc TN-DA} presented in Section~\ref{sect:original_problem} can be summarized as follows:
\begin{itemize}
    \item[--] Investment Function $\textsc{Invest}$: Sums of cost of opening each design arc, where the cost of opening $z_{hl}$ is denoted as $\beta_{hl}$. The unit of $\beta_{hl}$ should be consistent with $g^r$ and $\varphi^r$.
    \item[--] Subproblem $\text{\sc OPT-PATH}^r$: Finds the shortest path $\pi^r$ for trip $r$ from its origin to destination in terms of the  cost $g^r$. It is important to emphasize that this subproblem aims to emulate the functionality of a rider mobile app (a realistic rider app example can be found in the Appendix of \citet{van2023marta}), suggesting an ODMTS path by considering a convex combination of the agency’s monetary cost and the rider’s travel time. As a result, the solution to $\text{\sc OPT-PATH}^r$ is not necessarily the shortest path solely in terms of travel duration for trip $r$. Furthermore, the path can be multimodal, such as shuttle $\rightarrow$ bus $\rightarrow$ shuttle, or simply a single shuttle leg, namely \textit{direct shuttle}. However, ODMTS does not allow a path with only two shuttle legs involved, i.e., shuttle $\rightarrow$ shuttle. Note that the existence of direct shuttle trips will guarantee the subproblem is always feasible.
    \item[--] Function {\sc PATH-COST}: Computes the weighted cost $g^r$ of path $\pi^r$, based on the travel time and on-demand shuttle operating cost.
    \item[--] Function {\sc PATH-REVENUE}: Returns a constant revenue value $\varphi^r$, because the ticket price is the same for all paths in ODMTS in practice.
\end{itemize}
The detailed bilevel formulation of {\sc ODMTS-DA} can be found in Appendix~\ref{sect:odmts_da_appendix}.

\subsection{Existing Methods for Solving ODMTS-DA}
To solve the {\sc ODMTS-DA} problem,
\citet{basciftci2023capturing} proposed an exact decomposition algorithm whose resulting
design can be considered as the equilibrium point of the game between
the transit agency and the potential transit users. The effectiveness
of this exact algorithm was demonstrated in a case study conducted in
Ann Arbor and Ypsilanti area in Michigan, USA.  Nevertheless, this
exact algorithm becomes computationally challenging on large-scale
instances. Furthermore, due to the
complexity of the exact algorithm, it is even practically challenging
to obtain high-quality sub-optimal solutions on large-scale problems, demonstrating the need to develop effective heuristics for providing fast approximations to the optimal solutions of large-scale instances with desired adoption properties of the resulting network designs. 

As an alternative approach, another method, known as {\sc P-Path}, reformulates the corresponding bilevel problem into a single-level Mixed-Integer Programming (MIP) problem~\citep{guan2024path}. {\sc P-Path} is capable of finding optimal solutions efficiently for medium-sized {\sc ODMTS-DA} instances and within a reasonable time frame for large-scale instances. However, its computational resource requirements (e.g., memory) can be substantial because of the size of the resulting MIPs. Consequently, heuristic algorithms are beneficial, as they not only provide computational efficiency but might also offer a high-quality initial solution for {\sc P-Path}.

\subsection{Heuristic Algorithms for ODMTS-DA}
\label{subsect:odmts_da_heuristics}
This section illustrates how heuristic algorithms can be utilized to solve the {\sc ODMTS-DA} problem and presents additional technical details to be utilized as part of the arc-based algorithms tailored to this problem setting. Building upon the general concepts discussed in Section~\ref{sect:algorithms}, the heuristic algorithms make use of the ODMTS Design with Fixed Demand ({\sc ODMTS-DFD}) problem and the choice model separately. The {\sc ODMTS-DFD} problem is a specific form of the {\sc TN-DFD} problem that is tailored to the ODMTS settings. The full formulation for {\sc ODMTS-DFD} can be found in Appendix~\ref{sect:odmts_da_appendix}.
To evaluate a design found by {\sc ODMTS-DFD}, the evaluation function $\text{\sc Eval}(\mathbf{z})$ is adopted from the objective function of {\sc ODMTS-DA} where $\pi^r$ and $\delta^r$ can all be deduced from $\mathbf{z}$, and the revenue term $\varphi^r$ for core trips are dropped because it remains as a constant value:
\begin{equation}
  \text{\sc Eval-ODMTS}(\mathbf{z}) =  \sum_{h,l \in H} \beta_{hl} z_{hl} + \sum_{r \in T_{core}} p^r g_\pi^r + \sum_{r \in T_{latent}} p^r \delta^r (g_\pi^r - \varphi^r) \label{eq:design_obj_odmts}
\end{equation}

To apply trip-based algorithms to {\sc ODMTS}---Algorithms~\mynameref{alg:rho_GRAD},~\mynameref{alg:eta_GRRE},~and~\mynameref{alg:rho_GAGR}, the customizable value $v^r$ should be clearly defined in order to rank trips. The net profit of serving $p^r$ riders of each trip $r$
under a given network design $\mathbf{z}$ constitutes a critical factor in determining which latent trips should be considered. To this end, ${fee}
(\pi^r)$ corresponds to the monetary shuttle operating cost of the trip $r$
under the path $\pi^r$, and its monetary net-cost value can be defined as $\textit{net-spending}(\pi^r) = \textit{fee}(\pi^r) - \textit{income}(\pi^r)$ considering the income of utilizing
the ODMTS under path $\pi^r$ (which is a constant in this particular setting). It is important to reemphasize that $\textit{net-spending}$, $\textit{fee}$, and $\textit{income}$ represent monetary values, while cost $g_\pi^r$ and revenue $\varphi^r$ in the optimization frameworks are weighted values between monetary cost and convenience.

For arc-based algorithms---Algorithms~\mynameref{alg:arc_s1}~and~\mynameref{alg:arc_s2}, the set of all candidate sub-networks $D^k$ at each iteration $k$ is restricted to include only directed bus cycles. Given that $\textbf{z}_{temp}^k$,
$\textbf{z}_{unfixed}^k$, and $\textbf{z}_{fixed}$ all satisfy
constraint~\eqref{eq:odmts_da_upperLevelConstr1}, the algorithm of
\citet{johnson1975finding} can be used to find the set $D^k$ of
directed cycles (lines 6 and 7 in Algorithms~\mynameref{alg:arc_s1}~and~\mynameref{alg:arc_s2}, respectively). The arc-based algorithms have been tested with a number of expansion rules to construct the {\sc Expand} based on the specific nature of ODMTS. The expansion rules used at iteration $k$ includes:
\begin{enumerate}[label=(\textbf{\alph*})]
    \item All latent trips that adopt $\textbf{z}^k$. This is the Full Adoption Expansion Rule. \label{rule:all_adopt_odmts}
    \item All latent trips that adopt $\textbf{z}^k$ and are not served by direct shuttle, i.e., they use fixed-routes. \label{rule:use_fixed_route_odmts}
    \item All latent trips that adopt $\textbf{z}^k$ and satisfies $\text{UB}^r \leq \alpha^r \cdot t_{cur}^r$, where the term $\text{UB}^r$ is an upper bound on ODMTS travel time. The derivation of this bound can be found in Appendix~\ref{sect:odmts_da_appendix}. \label{rule:upper_bound_odmts}
\end{enumerate}
The definition of rule~\ref{rule:upper_bound_odmts} may seem overly complex; however, it leads to a useful proposition presented below, with its proof also provided in Appendix~\ref{sect:odmts_da_appendix}. 
\begin{prop}
\label{prop:arc_s1_proFA}
When using {\sc Time-Based} as Choice Model $\mathcal{C}^r$ (Equation~\eqref{eq:choice_model_time_based}) in {\sc ODMTS-DA}, the ODMTS design found by Algorithm~\mynameref{alg:arc_s1} with expansion rule~\ref{rule:upper_bound_odmts} satisfies the Correct Adoption property (Property~\ref{property:adopt}), i.e., $\%A_{false} = 0$.
\end{prop}

Additionally, when applying Algorithm \mynameref{alg:arc_s2} in practice, it is worth noticing that the first rule should not be rule~\ref{rule:all_adopt_odmts} as this rule includes all adopting latent trips, covering the trips considered in other expansion rules.

\subsection{Overview on ODMTS-DA Case Studies}
\label{sec:CaseStudyOverview}

This section summarizes the three {\sc ODMTS-DA} case studies. Table~\ref{table:odmts_da_3_case_studies} provides a comparison of the settings across these cases by providing instance and solution approach information. 
\begin{table}[!ht]
\centering
\resizebox{0.9\textwidth}{!}{ 
\begin{tabular}{l l r r r r r l l l}
\toprule
Size & Location & $|N|$ & $|H|$ & $|Z|$ & $|T_{core}|$ & $|T_{latent}|$ & Congestion & Choice Model $\mathcal{C}^r$& Methods \\ \midrule
Medium & Ypsilanti & 1267 & 10 & 90 & 937 & 566 & No & {\sc Time-Based} & Exact Alg., {\sc P-Path}, Heuristics \\ \midrule
Large & Atlanta & 2426 & 58 & 828 & 15,478 & 36,283 & No & {\sc Time-Transfer-Based}& {\sc P-Path}, Heuristics \\ \midrule
Extra Large & Atlanta &  2426 & 66 & 2034 & 15,478 & 36,283 & Yes & {\sc Time-Based} & Exact Alg., Heuristics \\ \bottomrule
\end{tabular}
}
\caption{A summary on the three ODMTS-DA case studies}
\label{table:odmts_da_3_case_studies}
\end{table}

The first case study utilizes a real dataset from
AAATA (\url{https://www.theride.org} Last Visited Date: September 30,
2025), the transit agency serving the greater Ypsilanti and Ann
Arbor region of Michigan, USA. This case study assumes that the agency will replace the existing transit system entirely by implementing an ODMTS. Since the exact algorithm and {\sc P-Path} can both successfully find the optimal design in a reasonable amount of time on this regular-sized dataset, this case
study aims at demonstrating the workability of the heuristics proposed in Section~\ref{sect:algorithms} by comparing
computational results obtained from different methods in terms of solution quality, adoption properties, and run time. Additional information on the dataset and experimental settings can be found in Appendix~\ref{sect:michigan_instance_appendix}.

The second case study conducted with real travel data represents a regular workday in Atlanta, Georgia, to evaluate the performance of the heuristic algorithms on large-scale instances. This case study assumes that the agency will retain the rail system, remove the existing bus services, and integrate the rail system with newly designed ODMTS. As a result, all existing four rail lines are preserved, and all 38 rail stations are designated as ODMTS hubs (set $H$). Visualizations and data processing details are available in Appendix~\ref{sect:atlanta_data_appendix}. As shown by \citet{guan2024path}, {\sc P-Path} can efficiently find the optimal solution for large-scale {\sc ODMTS-DA} instances within a reasonable time frame and clearly outperforms the exact algorithm in computational efficiency. Therefore, in this case study, only {\sc P-Path} and the heuristic algorithms are compared. Furthermore, to showcase that {\sc ODMTS-DA} can take the choice model $\mathcal{C}^r$ as a black-box, another choice model, namely {\sc Time-Transfer-Based}, is applied. In this choice model, a path $\pi$ is evaluated based on both travel duration and the number of required transfers, which is denoted by $l_\pi^r$. For instance, a ``shuttle $\rightarrow$ rail $\rightarrow$ bus'' path involves 2 transfers. Riders also have a transfer tolerance ($l_{ub}^r$), which is typically a low number and is set to 3 for all riders in this case study. The $\alpha$ value for all latent riders is set to 1.5.
\begin{equation}
\text{\sc Time-Transfer-Based}: \mathcal{C}^r(\pi) \equiv \mathbbm{1}(t_\pi^r \leq \alpha^r \ t^r_{cur}) \land  \mathbbm{1}(l_\pi^r \leq l_{ub}^r)
\end{equation}

The third case study uses the same data as the second case study, however it incorporates a larger set for $Z$ by considering more arcs between hubs (see details in Appendix~\ref{sect:atlanta_data_appendix}). Preliminary investigation indicates that {\sc P-Path} cannot handle this setting due to its massive memory requirement, exceeding 500 GB with commercial solver Gurobi 9.5. However, the exact algorithm, which is based on a decomposition approach, can partially solve this instance. Therefore, in this case, the exact algorithm is chosen as the baseline for evaluating the heuristic algorithms. The choice model is switched back to {\sc Time-Based} to showcase the connection with Proposition~\ref{prop:arc_s1_proFA}. Additionally, this case study incorporates congestion modeling to enhance its realism.

The design objectives reported in the following subsections are based on {\sc Eval-ODMTS}. Note that, although
only a subset of trips is used to obtain the design in all algorithms,
the design evaluation is computed with all core trips and latent
trips to account for all riders considered. 

\subsection{Results: Medium Instance in Ypsilanti}
\label{sec:MediumInstancepsilanti}

Table~\ref{tab:yp_exp_4_alg_perform} compares the computational performances of all solution methods along with the quality of their resulting network designs. The column `\# Total Itr.' reports the number of iterations, corresponding to the number of ODMTS-DFD problems solved during the runs. The column `\# Outer Itr.' indicates the number of iterations for~\mynameref{alg:rho_GAGR}, where its corresponding total iterations represent the number of iterations of~\mynameref{alg:eta_GRRE} (its subproblem). First, it can be seen that Algorithm~\mynameref{alg:rho_GAGR} requires significantly more computational time, whereas Algorithm~\mynameref{alg:eta_GRRE} takes about the same amount of computational time as Algorithm~\mynameref{alg:rho_GRAD}. Note that a practical approach to reduce the running time is to increase the step sizes of the trip-based algorithms. On the other hand, arc-based algorithms take less time compared to the trip-based algorithms. Secondly, four out of seven heuristic runs find the optimal solution, whereas Algorithm~\mynameref{alg:arc_s1} (with both rule~\ref{rule:all_adopt_odmts} and rule~\ref{rule:upper_bound_odmts}) and Algorithm~\mynameref{alg:arc_s2} (with rules~\ref{rule:upper_bound_odmts} and~\ref{rule:all_adopt_odmts}) report solutions with small optimality gaps. Despite having the same gap, the two arc-based algorithms report different results, due to their different expansion rules. Moreover, Propositions~\ref{prop:alg_1_prop}, \ref{prop:arc_s1_proFR}, and~\ref{prop:arc_s1_proFA} can be related to the results on {\sc False Rejection Rate} and {\sc False Adoption Rate}. In general, for transit agencies, small $\%R_{false}$ or $\% A_{false}$ values are indications of advantageous ODMTS designs when the optimal solution is unavailable. Three additional experiments with different parameter settings are provided in Appendix~\ref{sect:michigan_instance_appendix}, demonstrating the performances of the proposed algorithms while providing network designs with desired adoption properties.

\label{subsect:ypsi_case_study}
\begin{table}[!ht]
\centering%
\resizebox{\textwidth}{!}{
\begin{tabular}{l r r r r r r r r r r}%
    \toprule
    Method&Step Size&Trip Rule&\# Outer Itr.&\# Total Itr.&Run Time (min)&Design Obj & \% Opt Gap. &\% Gap &$\%R_{false}$&$\%A_{false}$\\%
    \midrule%
    Exact Alg.& - & - & - & - & 360.29 & 29962.79 & 1.85 & - & 0.00 & 0.00\\%
    \midrule
    {\sc P-Path} & - & - & - & - & 3.34 & 29962.79 & - & 0.00 & 0.00 & 0.00 \\%
    \midrule
    Alg.~\mynameref{alg:rho_GRAD}&10&{-}&1&26&2.48&29962.79&-&0.00&0.00&0.53\\%
    \midrule%
    Alg.~\mynameref{alg:eta_GRRE}&10&{-}&1&26&2.49&29962.79&-&0.00&11.13&0.53\\%
    \midrule%
    Alg.~\mynameref{alg:rho_GAGR}&10&{-}&26&357&47.14&29962.79&-&0.00&0.00&1.41\\%
    \midrule%
    Alg.~\mynameref{alg:arc_s1}&{-}&\ref{rule:all_adopt_odmts}&1&8&1.05&30019.14&-&0.19&0.00&34.28\\%
    \midrule%
    Alg.~\mynameref{alg:arc_s1}&{-}&\ref{rule:upper_bound_odmts}&1&7&1.13&30020.26&-&0.19&13.78&0.00\\%
    \midrule%
    Alg.~\mynameref{alg:arc_s2}&{-}&\ref{rule:use_fixed_route_odmts},~\ref{rule:all_adopt_odmts}&1&8&1.06&29962.79&-&0.00&0.00&2.65\\%
    \midrule%
    Alg.~\mynameref{alg:arc_s2}&{-}&\ref{rule:upper_bound_odmts},~\ref{rule:all_adopt_odmts}&1&8&0.95&30020.26&-&0.19&0.00&0.00\\%
    \bottomrule%
\end{tabular}%
}
\caption{Performance comparison on the medium instance. The exact algorithm reports optimality gap after 6 hours. The heuristics report the gap between themselves and the optimal solution ({\sc P-Path}).}
\label{tab:yp_exp_4_alg_perform}
\end{table}

\subsection{Results: Large Instance in Atlanta}

Table~\ref{table:atl_alg_perform_large} presents the results from heuristic algorithms on a large-scale instance, using the optimal solution from the {\sc P-Path} model as a benchmark. For trip-based iterative algorithms, three step sizes, i.e., 1000,
2000, and 3000, were tested, and the step size associated with the
best design objective is reported. The results show that while the heuristics do not reach the optimal solution, they still produce high-quality designs in a very short amount of time, whereas the optimal solution requires nearly five days of computation. Note that due to the structure of the {\sc ODMTS-DFD} problem, the heuristics can utilize multi-threading techniques; here, 12 threads are employed. Additionally, using the solution from Algorithm~\mynameref{alg:arc_s2} as a warm-start for {\sc P-Path} (because it provides the best solution in a shorter amount of time) can help reduce overall computation time. However, as pointed out by \citet{guan2024path}, starting {\sc P-Path} with a heuristic solution does not necessarily guarantee a faster path to the optimal solution.

One may wonder why heuristics are necessary if an optimal solution is available for instances up to a certain size, even though it is time-consuming to obtain. One of the reasons is that a transit agency may need to quickly develop a transit plan for special events, and not all locations have access to powerful computing resources capable of solving a large-scale {\sc P-Path} model. In these circumstances, the heuristics can be useful. Furthermore, as noted by \citet{guan2024path}, {\sc P-Path} may spend a substantial amount of computing time on closing the final 1\% of the optimality gap during the branch-and-bound processes. The solution reported here required 73\% of its total computation time (around 85 out of the 116.48 hours) to address this last 1\%, while the heuristics can achieve a better solution within just a few minutes.

\begin{table}[!ht]%
\centering%
\resizebox{\textwidth}{!}{
\begin{tabular}{l r r r r r r r r r}%
\toprule
Method&Best Step Size&Trip Rule&\# Outer Itr.&\# Total Itr.&Run. Time (hour)&Design Obj&\%Gap&$\%R_{false}$&$\%A_{false}$\\%
\midrule%
{\sc P-Path} & - & - & - & -& 116.48 & 214463.71 & 0.00 & 0.00 & 0.00 \\ \midrule
{\sc P-Path} Warm & - & - & - & -& 88.19 & 214463.71 & 0.00 & 0.00 & 0.00 \\ \midrule
Alg.~\mynameref{alg:rho_GRAD}&3000&{-}&1&8&0.07&214906.96&0.21&0.00&7.90\\%
\midrule%
Alg.~\mynameref{alg:eta_GRRE}&1000&{-}&1&17&0.15&221265.96&3.17&4.72&0.84\\%
\midrule%
Alg.~\mynameref{alg:rho_GAGR}&2000&{-}& 10 &69&0.66&214906.96&0.21&0.40&6.70\\%
\midrule%
Alg.~\mynameref{alg:arc_s1}&{-}&\ref{rule:all_adopt_odmts}&1&20&0.72&214846.58&0.18&0.00&14.37\\%
\midrule%
Alg.~\mynameref{alg:arc_s1}&{-}&\ref{rule:upper_bound_odmts}&1&13&0.31&242229.26&12.95&31.43&0.00\\%
\midrule%
Alg.~\mynameref{alg:arc_s2}&{-}&\ref{rule:use_fixed_route_odmts},\ref{rule:all_adopt_odmts}&1&22&0.60&214846.58&0.18&0.00&7.29\\%
\midrule%
Alg.~\mynameref{alg:arc_s2}&{-}&\ref{rule:upper_bound_odmts},\ref{rule:all_adopt_odmts}&1&22&0.50&214906.96&0.21&0.00&11.60\\%
\bottomrule
\end{tabular}%
}
\caption{A comparison between {\sc P-Path} and the heuristics on the large instance. {\sc P-Path} is solved to optimal. Warm stands for warm-start, the reported time includes both the heuristic and the {\sc P-Path} computation time.}
\label{table:atl_alg_perform_large}
\end{table}

\subsection{Results: Extra-Large Instance in Atlanta}
\label{subsect:atlanta_extra_large_instance}
This section evaluates the computational results of the extra-large case study from various viewpoints along with the analyses on the resulting network designs. Additional discussions on results such as adoptions, operating cost, travel time, and travel distance can be found in Appendix~\ref{sect:atlanta_real_study_appendix}.

\paragraph{Algorithmic Performances}
Table~\ref{table:atl_alg_perform_extra_large} compares the performances of all
algorithms. The most obvious finding to emerge from the analysis is
that the heuristic algorithms can outperform the exact algorithm in a
reasonable amount of time. Moreover, after 48 hours, the exact
algorithm finds a solution with objective 213828.87, which is only
1.61\% better than the 24-hour running time result. Thus, the proposed
iterative algorithms can efficiently provide high-quality solutions to
the {\sc ODMTS-DA} problem. The only exception is Algorithm
\mynameref{alg:arc_s1} with rule~\ref{rule:upper_bound_odmts}, whose greater objective
value can be attributed to the strict nature of rule~\ref{rule:upper_bound_odmts}.

\begin{table}[ht!]%
\centering%
\resizebox{\textwidth}{!}{
\begin{tabular}{l r r r r r r r r r }%
\toprule
Method&Best Step Size&Trip Rule&\# Outer Itr.&\# Total Itr.&Run. Time (hour)&Design Obj.&\% Gap&$\%R_{false}$&$\%A_{false}$\\
\midrule%
Exact Alg.  &-  & - &-&-&26.20&220557.59&-&-&-\\%
\midrule%
Alg.~\mynameref{alg:rho_GRAD}&3000&{-}&1&9&1.38&209548.18&{-}4.99&0.00&7.18\\%
\midrule%
Alg.~\mynameref{alg:eta_GRRE}&3000&{-}&1&11&1.66&213409.00&{-}3.24&5.46&0.49\\%
\midrule%
Alg.~\mynameref{alg:rho_GAGR}&3000&{-}&8&53&9.72&209412.04&{-}5.05&1.41&4.07\\%
\midrule%
Alg.~\mynameref{alg:arc_s1}&{-}&\ref{rule:all_adopt_odmts}&1&20&5.58&208664.07&{-}5.39&0.00&20.16\\%
\midrule%
Alg.~\mynameref{alg:arc_s1}&{-}&\ref{rule:upper_bound_odmts}&1&11&0.74&246561.82&11.79&63.78&0.00\\%
\midrule%
Alg.~\mynameref{alg:arc_s2}&{-}&\ref{rule:use_fixed_route_odmts}, \ref{rule:all_adopt_odmts}&1&25&3.69&211253.00&{-}4.22&0.00&11.48\\%
\midrule%
Alg.~\mynameref{alg:arc_s2}&{-}&\ref{rule:upper_bound_odmts}, \ref{rule:all_adopt_odmts}&1&27&4.75&208294.23&{-}5.56&0.00&18.48\\%
\bottomrule%
\end{tabular}%
}
\caption{Performance comparison between all algorithms on the extra-large instance. The exact algorithm reports the best feasible solution within a 24 hours running-time limit. '\% Gap' shows the gap between the objective value reported by the heuristics and the objective reported by the exact algorithm.}
\label{table:atl_alg_perform_extra_large}
\end{table}

The table also shows that Algorithms~\mynameref{alg:rho_GRAD}~and~\mynameref{alg:eta_GRRE} achieve the second and third fastest termination among all iterative algorithms, although they fail to discover a better
solution compared to the other approaches. These performances are likely to be related to the adoption set $C$ of Algorithm~\mynameref{alg:rho_GRAD} and the rejecting set $R$ of Algorithm~\mynameref{alg:eta_GRRE}. These sets are designed to speed up the design searching processes, yet it leads the
algorithm to overlook multiple intermediate ODMTS
designs. Note that multi-threading is disabled for this instance in order to produce a fair comparison on all algorithms. However, when the multi-threading techniques for subproblems in {\sc ODMTS-DFD} are allowed, the heuristics can be massively sped up. For instance, with 18 threads and a step size of 3000, Algorithm~\mynameref{alg:eta_GRRE} can terminate in 0.24 hours. On the other hand, the combined trip-based iterative algorithm can
indeed explore more ODMTS designs and provide a lower objective:
algorithm~\mynameref{alg:rho_GAGR} reports a better objective
than Algorithm~\mynameref{alg:eta_GRRE}.  However, it does not
outperform Algorithm~\mynameref{alg:rho_GRAD}. A possible
explanation for this result can be the choices of step sizes.

In this particular case study, Algorithm~\mynameref{alg:rho_GRAD} produces
a better balance between running time and design objective among all
the trip-based iterative algorithms. The reported $\%R_{false}$ value
is consistent with Proposition~\ref{prop:alg_1_prop}. In other words,
with ODMTS design with $\%R_{false} = 0$, all riders excluded during the
designing procedure reject the system. However, there will be
$\%A_{false}$ riders who reject the system even though they were
included in the design procedure, e.g., 7.18\% latent trips
reject the design when using Algorithm~\mynameref{alg:rho_GRAD} although
they are in $\hat{T}$. Thus, this type of ODMTS designs can be treated
as an upper limit when agencies plan to expand their transit network
to serve more passengers.

For the arc-based iterative algorithms, Algorithm
\mynameref{alg:arc_s1} has consistent $\%R_{false}$ and $\%A_{false}$
values with Propositions~\ref{prop:arc_s1_proFR} and
\ref{prop:arc_s1_proFA}. Algorithm~\mynameref{alg:arc_s1} with rule
\ref{rule:upper_bound_odmts} is the only case that guarantees $\%A_{false} = 0$ in
this study. That is, all riders used to design the ODMTS service will
adopt it. Although this is a relatively conservative approach, the
transit agency can start to expand the services based on this ODMTS
design. Secondly, the extended version of Algorithm
\mynameref{alg:arc_s1}, i.e., Algorithm~\mynameref{alg:arc_s2}, also
found high-quality solutions with a considerable drop in running time
and more iterations compared to Algorithm
\mynameref{alg:arc_s1} with trip rule~\ref{rule:all_adopt_odmts}. This result can be explained by the fact that the
expansion rule of Algorithm~\mynameref{alg:arc_s1} is significantly
looser than the first rule employed by Algorithm
\mynameref{alg:arc_s2}. Thus, the size of the $\hat{T}$ in Algorithm
\mynameref{alg:arc_s2} for each iteration is substantially smaller,
resulting in shorter running time per iteration. Closer inspection of
the table shows that the solution found by Algorithm
\mynameref{alg:arc_s2} (with rules~\ref{rule:upper_bound_odmts} and~\ref{rule:all_adopt_odmts})
outperforms all other algorithms and successfully establishes the smallest
upper bound on the design objective in a relatively short
duration. Furthermore, the results found by Algorithm
\mynameref{alg:arc_s2} are also consistent with Corollary
\ref{cor:alg_arc_s2_second_rule_a}.

In summary, the results compare the algorithmic performances of all
heuristic algorithms. A few practical recommendations for choosing
the algorithms might then be summarized: (i) use Algorithms
\mynameref{alg:rho_GRAD} or \mynameref{alg:eta_GRRE} with relatively large step size to conduct a
fast exploration process, (ii) employ Algorithm~\mynameref{alg:arc_s1}
or Algorithm~\mynameref{alg:arc_s2} with certain trip expansion rules
in order to build the upper bound of the problem, and (iii) apply
either Algorithm~\mynameref{alg:rho_GRAD} or Algorithm
\mynameref{alg:arc_s2} when the balance between solution quality and
running time is demanded.

\paragraph{ODMTS Designs}

\begin{figure}[!ht]
    \centering
    \begin{subfigure}[b]{0.25\textwidth}
        \includegraphics[width=\textwidth]{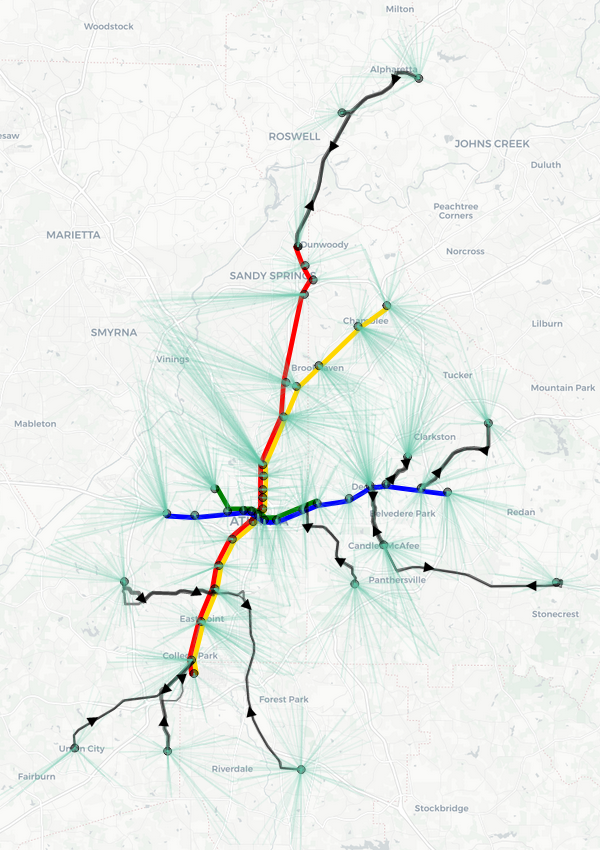}
        \tiny\caption{Only with Core Trips}
        \label{subfig:design_core}
    \end{subfigure}
    \begin{subfigure}[b]{0.25\textwidth}
        \includegraphics[width=\textwidth]{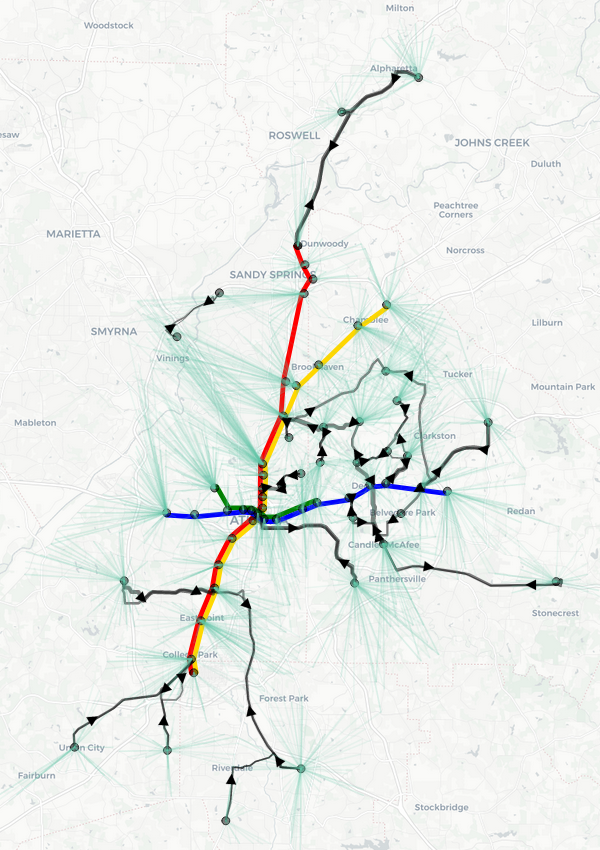}
        \caption{Exact Alg. in 24 hours}
        \label{subfig:design_exact_alg}
    \end{subfigure}
    \begin{subfigure}[b]{0.25\textwidth}
        \includegraphics[width=\textwidth]{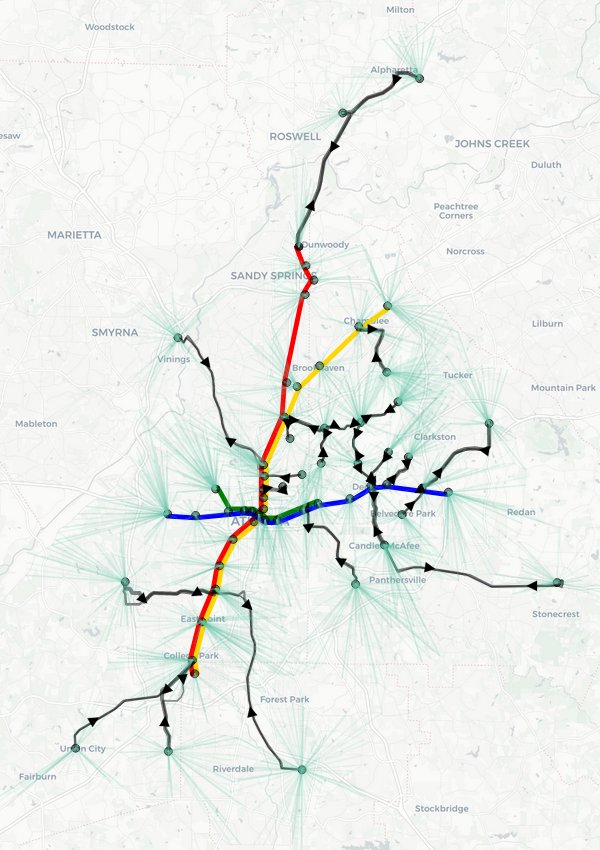}
        \caption{Alg.~\mynameref{alg:rho_GRAD}}
        \label{subfig:design_alg_rho_GRAD}
    \end{subfigure}
        \begin{subfigure}[b]{0.25\textwidth}
        \includegraphics[width=\textwidth]{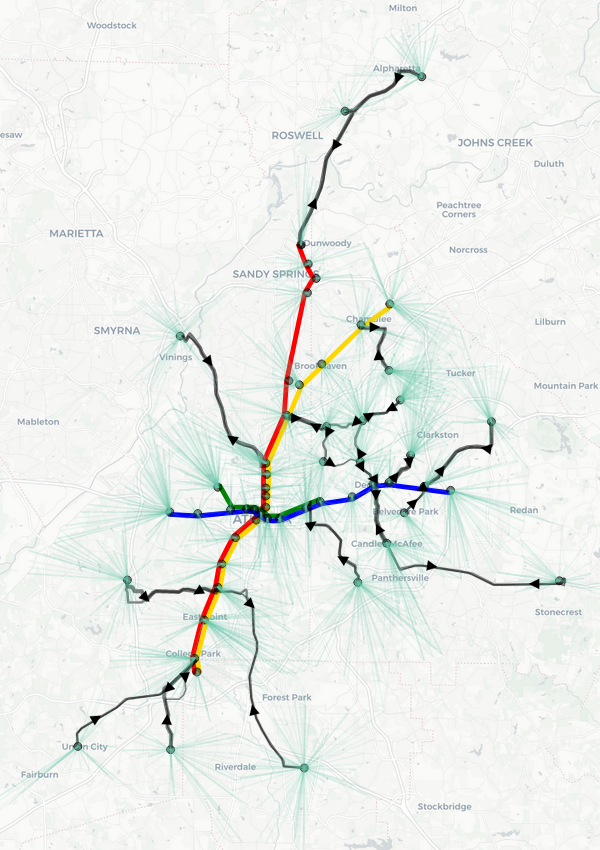}
        \caption{Alg.~\mynameref{alg:eta_GRRE}}
        \label{subfig:design_alg_MIAR}
    \end{subfigure}
    \begin{subfigure}[b]{0.25\textwidth}
        \includegraphics[width=\textwidth]{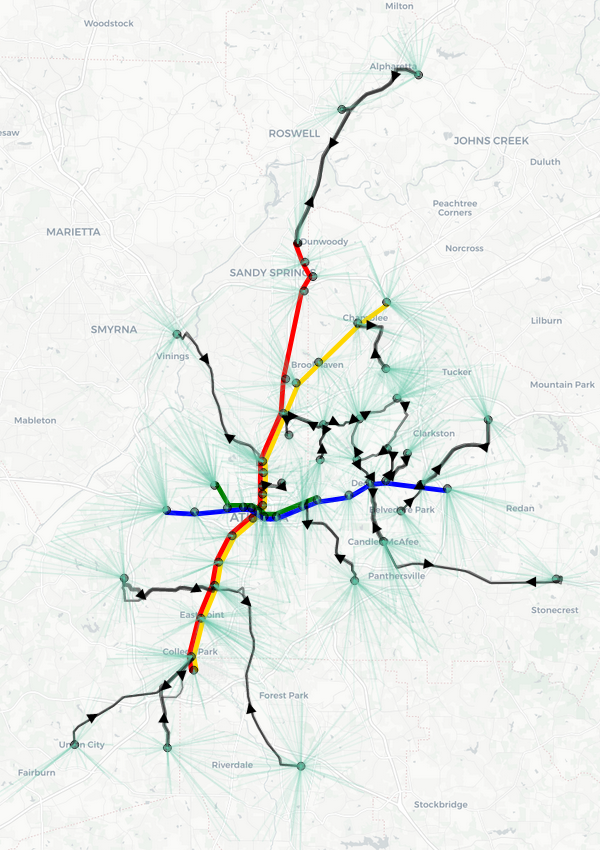}
        \caption{Alg.~\mynameref{alg:rho_GAGR}}
        \label{subfig:design_alg_GAGR}
    \end{subfigure}
    \begin{subfigure}[b]{0.25\textwidth}
        \includegraphics[width=\textwidth]{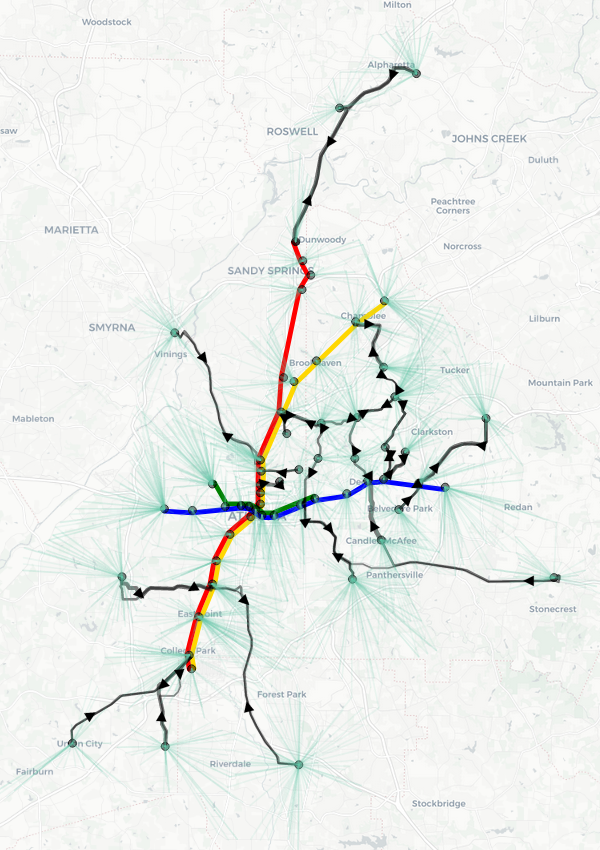}
        \caption{Alg.~\mynameref{alg:arc_s1} rule \ref{rule:all_adopt_odmts}}
        \label{subfig:design_alg_arc_s1_all}
    \end{subfigure}
        \begin{subfigure}[b]{0.25\textwidth}
        \includegraphics[width=\textwidth]{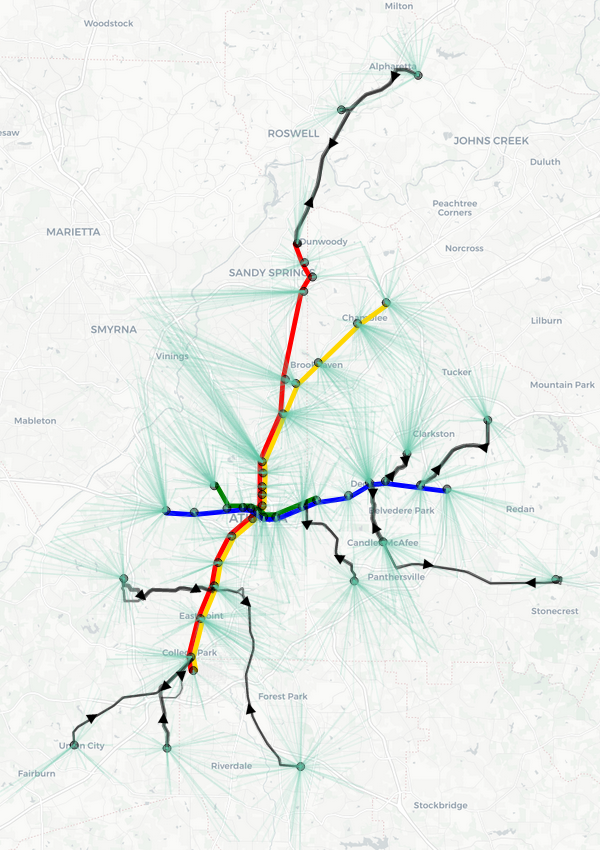}
        \caption{Alg.~\mynameref{alg:arc_s1} rule \ref{rule:upper_bound_odmts}}
        \label{subfig:design_alg_arc_s1_ub}
    \end{subfigure}
    \begin{subfigure}[b]{0.25\textwidth}
        \includegraphics[width=\textwidth]{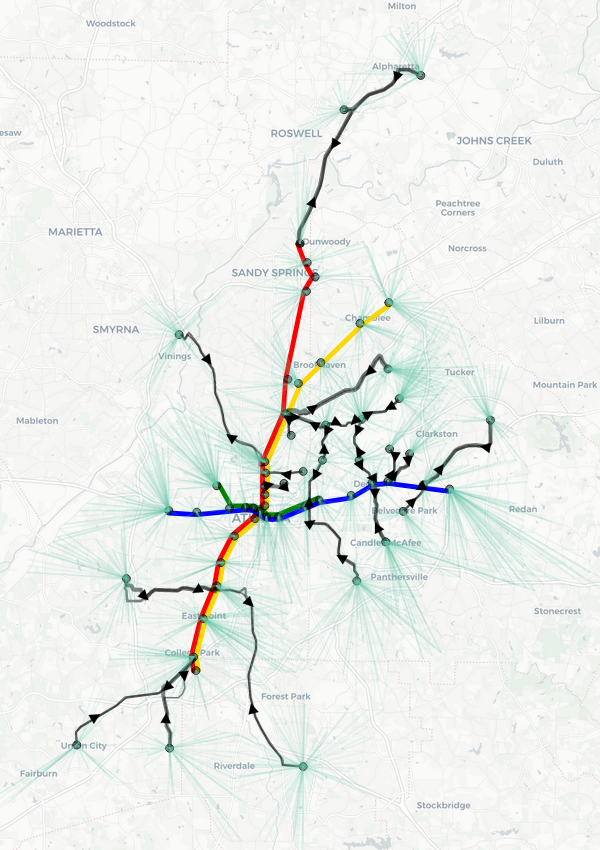}
        \caption{Alg.~\mynameref{alg:arc_s2} rules \ref{rule:use_fixed_route_odmts} \ref{rule:all_adopt_odmts}}
        \label{subfig:design_alg_arc_s2_mm_all}
    \end{subfigure}
    \begin{subfigure}[b]{0.25\textwidth}
        \includegraphics[width=\textwidth]{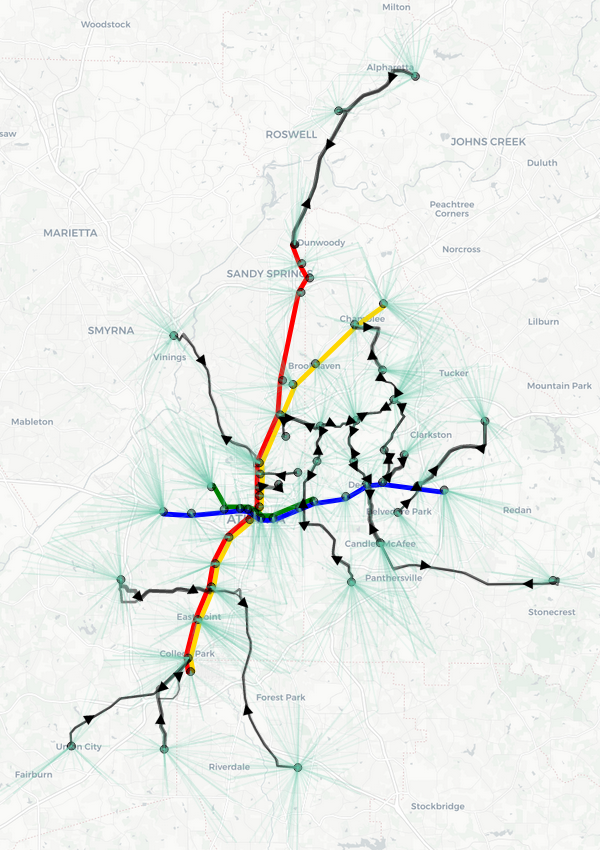}
        \caption{Alg.~\mynameref{alg:arc_s2} rules \ref{rule:upper_bound_odmts} \ref{rule:all_adopt_odmts}}
        \label{subfig:design_alg_arc_s2_ub_all}
    \end{subfigure}
\caption{The ODMTS designs found by the Algorithms for the extra-large case study. Black polylines and green straight lines are utilized to represent opened bus arcs and opened non-direct shuttle legs, respectively. Thickness of the black polylines indicates the number of users using
the arcs. The rail systems are in red, yellow, dark green, and blue.}
\label{fig:atl_all_designs}
\end{figure}

The ODMTS designs reported by the exact algorithm and all iterative algorithms are illustrated in Figure~\ref{fig:atl_all_designs}. In particular, Figure~\ref{subfig:design_core} presents the design when solely using the core trip set $T_{core}$, i.e., $\hat{T} = T_{core}$. It is thus the ODMTS design discovered by all trip-based heuristic algorithms at the first iteration. The other figures show the ODMTS designs obtained from the algorithms. The non-direct shuttle legs are on-demand shuttle services that connect a hub and a rider's origin stop or destination stop. For clearer interpretation of the maps, the main concentrations of latent trips in this case study are in the East Atlanta region (east of the red and yellow rail lines and north of the blue rail line), which is also where the primary differences in the designs appear. Additional details can be found in the latent demand heatmaps provided in Appendix~\ref{subsect:atlanta_latent_demand}. Compared to the design in Figure~\ref{subfig:design_core}, it can be seen that additional trips in $\hat{T}$ generally lead to a more complex network in the East Atlanta region. Closer inspection of all ODMTS designs show that the majority of bus arcs are used to deliver riders to the rail system. This might be explained by the fact that the rail system operates at a high frequency and connects the major districts such as Downtown (where all 4 rail lines intersect), Midtown (just north of Downtown, where the red and yellow lines overlap), and the city of Sandy Springs (the northernmost section of the red line). The simplest design
among Figure \ref{fig:atl_all_designs} is found by Algorithm
\mynameref{alg:arc_s1} with rule \ref{rule:upper_bound_odmts}, illustrated in Figure
\ref{subfig:design_alg_arc_s1_ub}, which is identical to the design with only core trips as input. The designs found by arc-based
algorithms with looser expansion rule can be seen in Figures~\ref{subfig:design_alg_arc_s1_all},~\ref{subfig:design_alg_arc_s2_mm_all},~and~\ref{subfig:design_alg_arc_s2_ub_all}, and these designs are
relatively more convoluted as expected. Additionally, Figure~\ref{subfig:design_exact_alg} presents the design reported by the
exact algorithm after 24 hours running time. Comparing to the other
designs, this design shows a few differences in the areas with latent riders, e.g., in the South and Northwest regions. This can
be explained by the 1\% or 8 hours terminating criteria for the master
problem in the ODMTS-DA problem (Figure \ref{fig:odmts_da_formulation}).
\section{TN-DA Example 2: SCTS Design with Adoptions Problem}
\label{sect:scts_da}

In order to show that {\sc TN-DA} can be tailored for different transit network design problems, this section introduces another example, namely Scooters-Connected Transit Systems Design with Adoptions ({\sc SCTS-DA}). As pointed out by \citet{mcqueen2022assessing} and \citet{yan2023evaluating}, for connecting to transit lines, e-scooters started to emerge as a solution to the first and last mile challenges. Different to the ODMTS example, the {\sc SCTS-DA} problem assumes there are two independent service operators involved: (1) a transit agency responsible for operating high-frequency buses and rails, and (2) a Transportation Network Company (TNC) that offers e-scooter services. This problem assumes a relatively simplified scenario where there is a sufficient number of e-scooters available throughout the city, and the TNC is responsible for managing and relocating the e-scooters as needed. In this context, the primary goal of the problem is to design an efficient high-frequency bus network that can be seamlessly connected with the e-scooter services while taking into account rider adoption considerations. 

\subsection{SCTS-DA Bilevel Formulation}

The analogy between {\sc SCTS-DA} and the generalized framework {\sc TN-DA}, can be summarized as follows:
\begin{itemize}
    \item[--] Investment Function $\textsc{Invest}$: Sums of cost of opening each design arc, where the cost of opening $z_{hl}$ is denoted as $\beta_{hl}$. The unit of $\beta_{hl}$ should be consistent with $g^r$ and $\varphi^r$.
    \item[--] Subproblem $\text{\sc OPT-PATH}^r$: Finds the shortest path $\pi^r$ for trip $r$ in terms of the cost $g^r$. Assume that the riders' primary goal is minimizing travel time; this is equivalent to the conventional shortest path based on travel time. A direct scooter trip between the origin and destination ensures that the subproblem is feasible.
    \item[--] Function {\sc PATH-COST}: Computes the weighted cost $g^r$ of path $\pi^r$. The cost $g^r$ is weighted based on $\pi^r$'s travel time of path $t_\pi^r$.
    \item[--] Function {\sc PATH-REVENUE}: Returns constant revenue value $\varphi^r$ because transit ticket price is flat for all users. Scooter fees are not considered as part of the revenue as the scooter services are not offered by the transit agency.
\end{itemize}
The detailed bilevel formulation of {\sc SCTS-DA} can be found in Appendix~\ref{sect:scts_da_appendix}. 

The decision model employed in this case study is referred to as {\sc Time-Based-Scooter} (Equation~\eqref{eq:choice_model_time_based_with_scooter}). It takes into account three factors: firstly, the travel duration within the transit system is required to be smaller than or equal to an $\alpha^r$ value multiplied by the travel duration of the current mode of transportation. Secondly, riders will only adopt the transit system if they actually utilize buses or trains, which will incur income for the transit agency. Lastly, riders only use the transit system if they travel on scooter less than a threshold value $d_{ub}^r$ given that scooters are not an efficient mode for longer distances. In the following equation, $x_{hl}$ and $y_{ij}$ denote if trip $r$ use design arc $(h,l)$ and scooter arc $(i,j)$. The travel distance between locations $i$ and $j$ is~$d_{ij}$.
\begin{equation}
\label{eq:choice_model_time_based_with_scooter}
\text{{\sc Time-Based-Scooter}}: {\cal C}^r(\pi) \equiv \mathbbm{1}(t_\pi^r \leq \alpha^r \ t^r_{cur}) \land \mathbbm{1}\bigg(\sum\limits_{h,l \in H} x_{hl}^r > 0 \bigg) \land \mathbbm{1}\bigg(\sum\limits_{i,j \in N} y_{ij}^r \cdot d_{ij} \leq d_{ub}^r\bigg)
\end{equation}

\subsection{Heuristic Algorithms for SCTS-DA}
This section presents the heuristic settings for {\sc SCTS-DA}; the formulations and the evaluation function are documented in Appendix~\ref{sect:scts_da_appendix}.

In the case of the trip-based iterative algorithms, the value used to rank trips, denoted as $v^r$, is determined based on the time traveled by e-scooters since e-scooter services typically charge riders based on the amount of time they utilize the scooters. 
For arc-based algorithms, the expansion rules used at iteration $k$ include:
\begin{enumerate}[label=(\roman*)]
    \item All latent trips that adopt $\textbf{z}^k$. This is the Full Adoption Expansion Rule. \label{rule:all_adopt_scts}
    \item All latent trips that adopt $\textbf{z}^k$ and use the backbone transit lines. Backbone refers to arcs that are preserved in the system without requiring design through solving {\sc TN-DA}, e.g., rail system in Atlanta. \label{rule:backbone_scts}
\end{enumerate}

The first rule considers an expansive case by including all latent riders adopting, whereas second rule considers a subset of these trips by reflecting a common scenario where scooters are used for first-and last-mile connections to primary transit corridors, typically rail lines. For instance, these transit stations are usually equipped with convenient parking areas for scooters in Atlanta.

\subsection{Results: SCTS-DA in Atlanta}
Similar to the large-scale {\sc ODMTS-DA} case studies explained in the previous section for Atlanta, this case study employs the same dataset and parameters for buses and rails (see Appendix~\ref{subsect:scts-da-settings}). For scooters, the average speed is assumed to be \SI{16.1}{km/h} (\SI{10}{mph}) when traveling in Atlanta, and the scooter distance bound $d_{ub}^r$ is set at \SI{24}{km} (\SI{15}{miles}) for all trips. Moreover, this case study also employs the congested matrices for cars and buses. Due to the similar structure of this problem to {\sc ODMTS-DA}, a modified version {\sc P-Path} is derived to seek the optimal solution of {\sc SCTS-DA}, and treated as a benchmark for this setting.

\begin{figure}[!t]
    \centering
    \begin{subfigure}[b]{0.28\textwidth}
        \includegraphics[width=\textwidth]{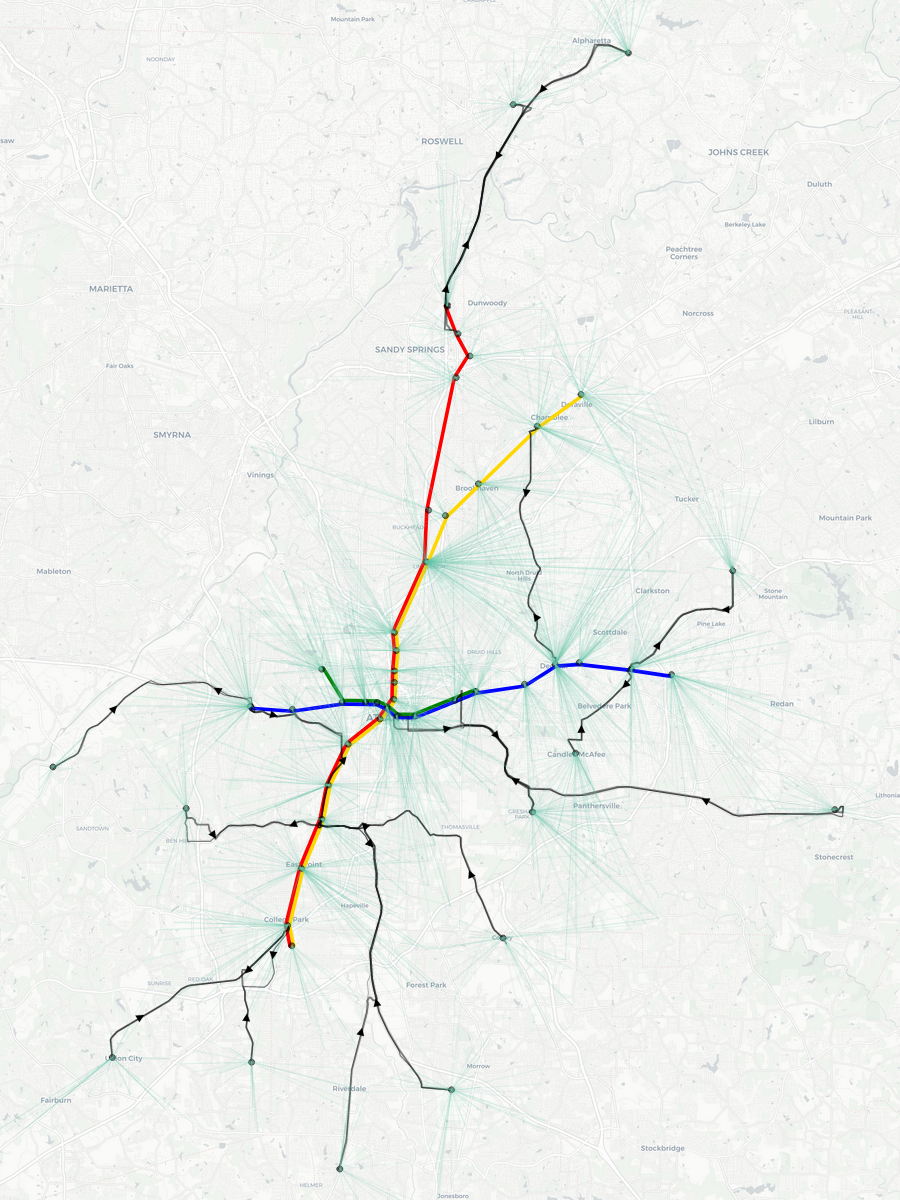}
        \caption{
            Trip-Based Algorithms
            \newline
        }
        \label{subfig:scooter_design_alg_eta_GRRE}
    \end{subfigure}
    \begin{subfigure}[b]{0.28\textwidth}
        \includegraphics[width=\textwidth]{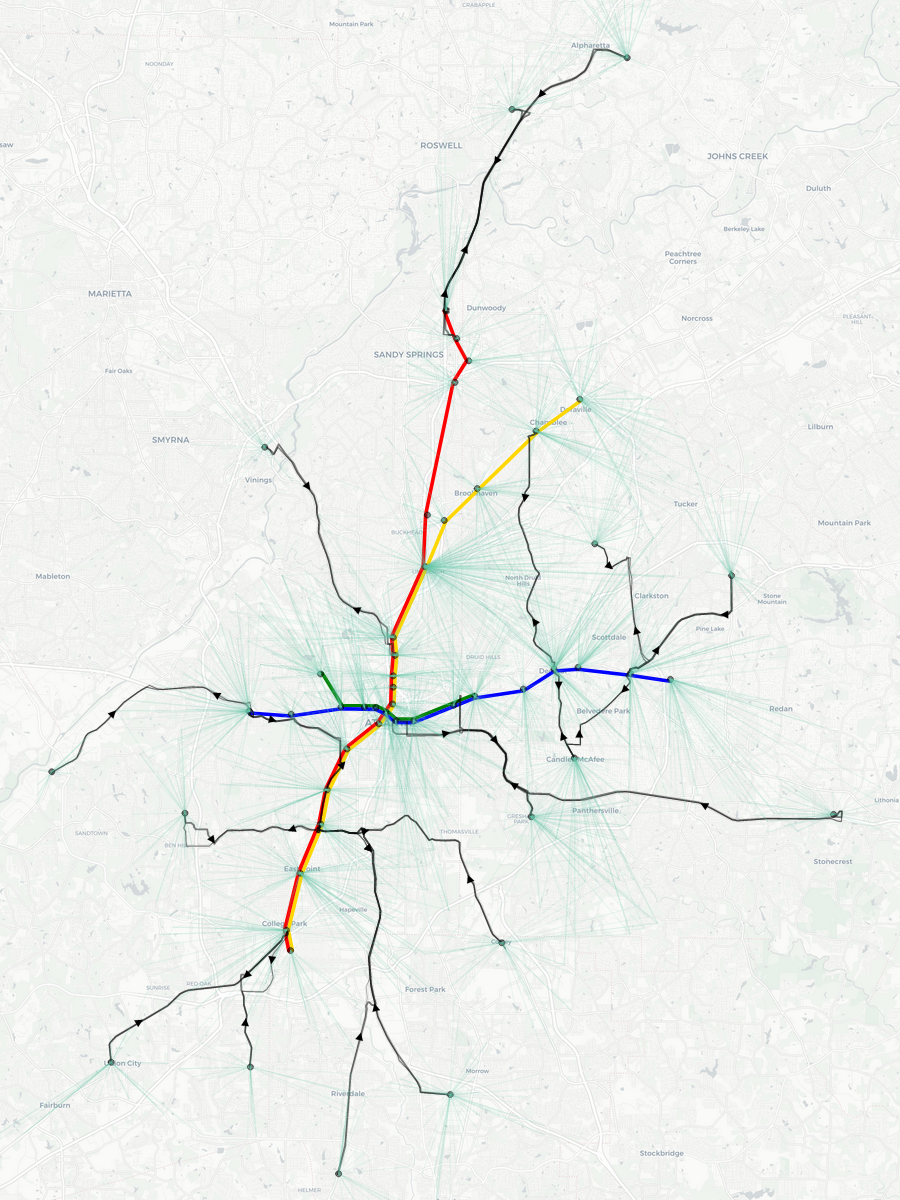}
        \caption{
            Arc-Based Algorithms and Modified {\sc P-Path}
        }
        \label{subfig:scooter_design_alg_arc_s2_backbone_alladopt}
    \end{subfigure}
\caption{The Designs Found by the Heuristic Algorithms and modified {\sc P-Path} for the SCTS-DA problem in Atlanta. Light green lines represent scooter connections.}
\label{fig:atl_scooter_all_designs}
\end{figure}

The computational results are reported in Table~\ref{table:atl_scooter_alg_perform}. Here, the two arc-based algorithms can find a solution that have zero {\sc False Rejection Rate} and {\sc False Adoption Rate}, which are the desired properties of an optimal solution. Although  Algorithms~\mynameref{alg:rho_GRAD}~and~\mynameref{alg:eta_GRRE} find the solution with higher objective, they terminate in a short duration. On the other hand, Algorithm~\mynameref{alg:rho_GAGR} is ineffective in producing advantageous results for this case study. In Figure~\ref{fig:atl_scooter_all_designs}, it can be observed that the designs obtained by the three trip-based heuristic algorithms (Figure~\ref{subfig:scooter_design_alg_eta_GRRE}) are identical, as do the designs discovered by the two arc-based algorithms (Figure~\ref{subfig:scooter_design_alg_arc_s2_backbone_alladopt}). Furthermore, it is worth mentioning that the designs resulting from the trip-based algorithms are quite similar to the design derived exclusively from core trips.

\begin{table}[t!]%
\centering%
\resizebox{\textwidth}{!}{
\begin{tabular}{l r r r r r r r r r }%
\toprule
Method&step size&trip rule&\# outer itr.&\# total itr.& time (hours)&design obj& \% Gap&$\%R_{false}$&$\%A_{false}$\\%
\midrule%
Modified {\sc P-Path} &{-}&{-}&{-}&{-}&4.58&1051031.55&0.00&0.00&0.00\\
\midrule
Alg.~\mynameref{alg:rho_GRAD}&500&{-}&1&12&0.11&1055265.30&0.40&11.61&0.00\\%
\midrule%
Alg.~\mynameref{alg:eta_GRRE}&500&{-}&1&12&0.10&1055265.30&0.40&11.61&0.00\\%
\midrule%
Alg.~\mynameref{alg:rho_GAGR}&500&{-}&12&79&0.70&1055265.30&0.40&11.61&0.00\\%
\midrule%
Alg.~\mynameref{alg:arc_s1}&{-}&\ref{rule:all_adopt_scts}&1&15&0.55&1051031.55&0.00&0.00&0.00\\%
\midrule%
Alg.~\mynameref{alg:arc_s2}&{-}&\ref{rule:backbone_scts}, \ref{rule:all_adopt_scts}&1&16&0.57&1051031.55&0.00&0.00&0.00\\%
\bottomrule%
\end{tabular}%
}
\caption{Performance comparison between all algorithms on the Large-scale experiment.}
\label{table:atl_scooter_alg_perform}
\end{table}

\section{Conclusion}
\label{sect:conclusion}

This study aims at capturing the latent demand for the transit
agency such that the transit services can concurrently provide
high-quality services to existing riders and extend the services to
current underserved areas. The decision-making processes of both the transit agency and the latent riders form a game-like interaction. To model this dynamic interaction, this paper first introduces a bilevel optimization framework called Transit Network Design with Adoptions ({\sc TN-DA}). The framework essentially models a network design problem with an incorporated black-box-like mode choice model. Achieving success in solving the {\sc TN-DA} can identify the equilibrium point between the agency and the latent riders. However, even with an exact solution algorithm, due to the bilevel and combinatorial nature of the {\sc TN-DA}, practically solving it can lead to computational challenges, particularly for large-scale instances.

This paper then attempts to address this limitation and proposes five
heuristic algorithms to efficiently approximate the optimal solution
of the {\sc TN-DA} problem. The fundamental idea of these iterative
algorithms is to separate the previous bilevel framework into two
individual but simpler components, i.e., a regular network design problem with
no adoption awareness (the {\sc TN-DFD} problem in Figure
\ref{fig:tn_dfd_formulation}) and an evaluation procedure based on a mode
choice function.  These
heuristic algorithms can be further categorized into two general
classes based on their algorithmic designs: (1) trip-based algorithms
and (2) arc-based algorithms. In general, at each iteration, the
trip-based algorithms reconstruct the trip sets for the {\sc TN-DFD}
problem and the arc-based algorithms attempt to fix a certain amount
of bus arcs. The advantages and limitations of each algorithm are
discussed in this paper, and several theoretical properties are
derived from the optimal solutions in order to provide more guidelines
for the iterative algorithms under varied circumstances. In
particular, two key properties, correct rejection
(Property~\ref{property:reject}) and correct adoption
(Property~\ref{property:adopt}), are introduced.
Then, the {\sc false rejection
rate} and the {\sc false adoption rate} are proposed as metrics to evaluate the
ODMTS designs, providing additional insights on the performance of the
heuristic algorithms when an optimal solution becomes unavailable
for large-scale instances.

This paper further presents comprehensive case studies on two different transit systems to validate the heuristics. For the first transit system---On-Demand Multimodal Transit Systems (ODMTS), one case study is carried out with a
normal-sized transit dataset collected by the transit agency in the
Ann Arbor and Ypsilanti area in the state of Michigan. The results
demonstrate that the heuristic algorithms can discover the optimal
ODMTS design. Another case study is based on multi-sourced data and
includes a very large numbers of stops and trips in Atlanta in the state of Georgia. The main goal of
this case study is to show the advantages of the heuristic algorithms
on large-scale instances. Compared to the exact algorithm, the
heuristic algorithms achieve better objective values in considerably
less amount of time. The second transit system---Scooters-Connected Transit Systems (SCTS) are discussed in Section~\ref{sect:scts_da}. These case studies serve as evidence that the {\sc TN-DA} framework can be readily tailored to various transit systems and exhibit strong performance on large-scale scenarios while providing practical insights.

\ACKNOWLEDGMENT{
This work was partially supported by NSF Leap-HI Grant No. 1854684, NSF Grant No. 2112533, NSF CMMI Grant No.  2434302, and T-SCORE Grant 69A3552047141 from the Department of Transportation.
}



\bibliographystyle{informs2014trsc} 
\bibliography{ref} 


%
%
%
\newpage
\begin{APPENDICES}
\section{Additional Information For Heuristic Algorithms}
\label{sect:heuristic_algorithm_appendix}

This section completes Section~\ref{sect:algorithms} with algorithms. The corresponding propositions are represented in this section again for enhancing readability.

\subsection{Greedy Adoption Algorithm with Greedy Rejection Subproblem}

This section presents Algorithm~\mynameref{alg:rho_GAGR}.

\begin{algorithm}[!ht]
\SingleSpacedXI
\caption{$\rho$-GAGR}
\label{alg:rho_GAGR}
\begin{algorithmic}[1]
\STATE Set $k = 0$, $C = \emptyset$, $\overline{T}^0 = T_{core} \cup C$, and $o^* = \infty$.
\STATE \textbf{Require}: step-size parameters $\rho$ and $\eta$
\WHILE{True}
\STATE Set $\mathbf{z}^k = \mynameref{alg:eta_GRRE}$ (Algorithm \mynameref{alg:eta_GRRE} starts from $\overline{T}^k$ instead of $T_{core}$)
\IF{$\text{\sc Eval}(\mathbf{z}^k) < o^*$}
\STATE Set $o^* = \text{\sc Eval}(\mathbf{z}^k)$.
\STATE Set $\mathbf{z}_{min} = \mathbf{z}^k$.
\ENDIF
\STATE Carry out steps 5---21 of \mynameref{alg:rho_GRAD}
\ENDWHILE
\RETURN $\mathbf{z}_{min}$
\end{algorithmic}
\end{algorithm}

\subsection{Arc-based Extended Greedy Algorithm}
This section presents Algorithm~\mynameref{alg:arc_s2}.
\begin{algorithm}[!t]
\SingleSpacedXI
\caption{arc-S2}
\label{alg:arc_s2}
\begin{algorithmic}[1]
\STATE Set $k = 0$, $\mathbf{z}_{fixed} = \Vec{0}$, $\overline{T}^0 = T_{core}$, $B^0 = \infty$, $p= 1$
\STATE \textbf{Require}: {\sc Expand}$_i$ functions for increasing the set of users considered ($i =1,2$)
\WHILE{$ p\leq 2$}
\WHILE{True}
\STATE Step 4---22 from Algorithm~\ref{alg:arc_s1}.
\STATE Set $T_{next}^k$ = {\sc Expand}$_p$($T_{latent},\mathbf{z}^{k}$).
\STATE Step 24---25 from Algorithm~\ref{alg:arc_s1}.
\ENDWHILE
\STATE Set $p = p+1$.
\ENDWHILE
\RETURN $\mathbf{z}^k$
\end{algorithmic}
\end{algorithm}

\section{Additional Information on ODMTS-DA}
\label{sect:odmts_da_appendix}
This section completes Section~\ref{sect:odmts_da} by presenting the full formulations for {\sc ODMTS-DA} and {\sc ODMTS-DFD} and their modeling details with omitted proofs.

\subsection{ODMTS-DA Bilevel Formulation}
This section presents the bilevel formulation of {\sc ODMTS-DA}, where extensive discussions on this formulation can be found in \citet{basciftci2023capturing, guan2024path}. Under the context of {\sc ODMTS-DA}, the set $N$ represents stops for shuttles to pick-up and drop-off passengers. The set $H \subseteq N$ corresponds to the hubs that can be stops for high-frequency buses or rail services. In order to maintain conciseness and avoid complexity, this section focuses on the formalization of buses, as establishing new rail lines in practice can be challenging. However, existing rail lines or similar services (e.g., Bus rapid transit) can be preserved as \textit{backbone arcs} and integrated into {\sc ODMTS-DA} by utilizing a set $Z_{backbone}$. To have a weighted objective considering both convenience and monetary cost aspects of the ODMTS design, a parameter $\theta \in [0,1]$ is defined such that the trip duration associated with convenience is multiplied by $\theta$ and monetary spending is multiplied by $1
 - \theta$.  In particular, the transit agency invests a weighted cost $\beta_{hl} = (1 - \theta) b_{dist} \,
 n \, d_{hl}^{bus}$ if the leg between hubs $h, l \in H$ is open, where $b_{dist}$ is the cost of using a bus per kilometer and $n$ is the number of buses operating in each open leg within the planning
horizon. The weighted investment $\beta_{hl}$ can alternatively be modeled as $\beta_{hl} = (1 - \theta) b_{time} \, n \, t_{hl}^{bus}$ by the transit agency, where $b_{time}$ stands for the monetary cost of operating a bus per hour. The parameter $n$ is determined in advance prior to solving the problem, and a common approach is to compute it based on the headway of buses. For backbone arcs in $Z_{backbone}$, $\beta_{hl}$ can be simply set to 0 since they are preserved from the existing transit system. The bilevel optimization model uses binary variable $z_{hl}$ whose value is 1 if the leg between the hubs $h,l \in H$ is open.

For each trip $r \in T$, the model uses binary variables $x_{hl}^r$ and $y_{ij}^r$ to determine whether this trip utilizes the leg between the
hubs $h,l \in H$, and the shuttle leg between the stops $i, j \in N$, respectively. By combining the selected $x_{hl}^r$ and $y_{ij}^r$ variables, a path $\pi^r$ is provided to riders of trip~$r$ from their origin to destination. To this end, this section interchangeably uses $\pi^r$ and vectors $(\mathbf{x^r}, \mathbf{y^r})$ to denote the paths of the riders. The transit agency incurs a weighted inconvenience cost $\tau^r_{hl} = \theta (t_{hl}^{bus} + t_{hl}^{wait})$ for the parts of the trips traveled by buses, where $t_{hl}^{wait}$ is the average waiting time of a passenger for the bus between hubs $h$ and $l$. A similar equation can be applied if there exists backbone connection between hubs $h$ and $l$, e.g., $\tau^r_{hl} = \theta (t_{hl}^{rail} + t_{hl}^{wait})$. For the legs of the trips traveled by on-demand shuttles, the transit
agency incurs a weighted cost and inconvenience cost $\gamma_{ij}^r = (1 - \theta) \omega \, d_{ij}^{shuttle} + \theta t_{ij}^{shuttle}$ for each trip $r \in T$, where $\omega$ is the monetary cost of using a shuttle per kilometer.
Core riders are assumed to  utilize the ODMTS once the system is built, following the general assumptions in {\sc TN-DA}. This is justified by the case study in \citet{basciftci2023capturing} that demonstrated the improvement in
convenience of the suggested paths compared to the existing transit
system. On the other hand, the latent riders can decide whether to adopt
the ODMTS, based on the choice function ${\cal
  C}^r(\mathbf{x^r},\mathbf{y^r})$.  The
transit agency charges riders a fixed fare $\phi$ to use the ODMTS,
irrespective of their paths. Thus, the fixed value $\varphi = (1 -
\theta)\phi$ becomes an additional revenue to the transit agency for
riders adopting the ODMTS. Additionally, as compared to the general formulation of {\sc TN-DA}, the objective function of {\sc ODMTS-DA} can exclude the revenue from core riders since it remains as a constant (see the second term in the objective from both {\sc TN-DA} and {\sc ODMTS-DA}).

\begin{figure}[ht]
\begin{subequations} \label{eq:odmts_da_upperLevelProblem}
\begin{alignat}{1}
\min_{z_{hl}, \delta^r} \quad & \sum_{h,l \in H} \beta_{hl} z_{hl} + \sum_{r \in T_{core}} p^r (g^r - \varphi^r) + \sum_{r \in T_{latent}} p^r \delta^r (g^r - \varphi^r) \label{eq:odmts_da_upperLevelObj} \\
\text{s.t.} \quad & \sum_{l \in H} z_{hl} = \sum_{l \in H} z_{lh} \quad \forall h \in H \label{eq:odmts_da_upperLevelConstr1} \\
& z_{hl} = 1 \quad \forall (h,l) \in Z_{backbone} \label{eq:odmts_da_fixed_arcs} \\
& \delta^r = {\cal C}^r(\mathbf{x^r}, \mathbf{y^r}) \quad \forall r \in T_{latent} \label{eq:odmts_da_userChoiceModel} \\
& z_{hl} \in \{0,1\} \quad \forall h,l \in H  \label{eq:odmts_da_binaryConstraint} \\
& \delta^r \in \{0,1\} \quad \forall r \in T_{latent} \label{eq:odmts_da_continuousConstraint} 
\end{alignat}
\end{subequations}
where $(\mathbf{x^r}, \mathbf{y^r}, g^r)$ are a solution to the optimization problem
\begin{subequations}
\label{eq:odmts_da_lowerLevelProblem}
\begin{alignat}{1}
\lexmin_{x_{hl}^r, y_{ij}^r, g^r, t^r} \quad &  \langle g^r, t^r \rangle \label{eq:odmts_da_lowerLevelObj} \\
\text{s.t.} \quad
   & g^r = \sum_{h,l \in H} \tau_{hl}^r x_{hl}^r + \sum_{i,j \in N} \gamma_{ij}^r y_{ij}^r \label{eq:odmts_da_drDefinition} \\
   & t^r = \sum_{h,l \in H}  (t_{hl}^{bus/rail} + t_{hl}^{wait}) x_{hl}^r + \sum_{i,j \in N}  t_{ij}^{shuttle} y_{ij}^r \label{eq:odmts_da_frDefinition} \\
   & \sum_{\substack{h \in H \\ \text{if } i \in H}} (x_{ih}^r - x_{hi}^r) + \sum_{i,j \in N}  (y_{ij}^r - y_{ji}^r) = \begin{cases}
     1 & \text{if  } i = or^r \\
    -1 & \text{if  } i = de^r \\
    0 & \text{otherwise}
    \end{cases} \quad \forall i \in N \label{eq:odmts_da_minFlowConstraint} \\
  & x_{hl}^r \leq z_{hl} \quad \forall h,l \in H \label{eq:odmts_da_openFacilityOnlyAvailable} \\
  & x_{hl}^r \in \{0,1\} \quad \forall h,l \in H, y_{ij}^r \in \{0,1\} \quad \forall i,j \in N. \label{eq:odmts_da_integralityFlowConstr}
\end{alignat}
\end{subequations}
\caption{The Bilevel Optimization Model for ODMTS Design with Adoptions ({\sc ODMTS-DA}). }
\label{fig:odmts_da_formulation}
\end{figure}

Figure~\ref{fig:odmts_da_formulation} presents the bilevel formulation of {\sc ODMTS-DA}. The leader problem~\eqref{eq:odmts_da_upperLevelProblem} designs the
network between the hubs for the ODMTS, whereas the follower problem~\eqref{eq:odmts_da_lowerLevelProblem}
obtains paths for each trip $r \in T$ by considering the legs of the
network design and the on-demand shuttles to serve the first and last
miles of the trips. The objective of the leader problem~\eqref{eq:odmts_da_upperLevelObj} minimizes
the sum of (i) the investment cost of opening bus legs, (ii) cost of the trips of the existing riders, and (iii) the net-cost of the riders with choice that are adopting the ODMTS.
Constraint~\eqref{eq:odmts_da_upperLevelConstr1} ensures weak connectivity
between the hubs with the same number of incoming and outgoing open
legs. Constraint~\eqref{eq:odmts_da_fixed_arcs} assures all preserved arcs are open. Constraint~\eqref{eq:odmts_da_userChoiceModel} represents the mode
choice of the riders in $T_{latent}$ depending on the ODMTS paths. Since the key of this entire section is to demonstrate the
bilevel framework and heuristic algorithms on a particular transit network design problem, without loss of generality,
this section assumes that {\sc ODMTS-DA} only considers a unique bus frequency for potential bus arcs $z_{hl}$.
However, incorporating multiple frequencies can be achieved by introducing a frequency variable for each potential bus arc
and adjusting Constraint~\eqref{eq:odmts_da_upperLevelConstr1} to
ensure flow balance when considering multiple frequencies.

For each trip $r$, the objective of the follower problem
\eqref{eq:odmts_da_lowerLevelProblem} minimizes the lexicographic objective
function $\langle g^r, t^r \rangle$, where $g^r$ corresponds to the cost of trip $r$ and $t^r$ breaks potential ties
between the solutions with the same value of $g^r$ by returning the
most convenient path for the rider of trip $r$. Once again, when solving the shortest-path problem based on the cost metric $g^r$, the transit agency may not always provide the most convenient path to the riders, as their decisions are influenced by spending considerations as well. Since
sub-objective $g^r$ contains both monetary spending and duration of trip $r$,
$g^r$ includes sub-objective $t^r$ multiplied by $\theta$.  As
presented in \citet{basciftci2023capturing}, for any network design, a
lexicographic minimizer of the follower problem
\eqref{eq:odmts_da_lowerLevelProblem} exists, and the lexicographic minimum is
unique. Thus, once a network design $\mathbf{z}$ is given, the path of the riders in trip $r$ and their choices, i.e., $x_{hl}^r$, $y_{ij}^r$, and $\delta^r$, can all be deduced.
This aligns with the requirements in {\sc TN-DA}, as it ensures that the subproblem offers a unique solution. Constraint
\eqref{eq:odmts_da_minFlowConstraint} guarantees flow conservation used in trip $r$ for origin, destination and each of
the intermediate points utilized. Constraint
\eqref{eq:odmts_da_openFacilityOnlyAvailable} ensures that the path only takes
into account the open legs between the hubs.  As pointed out by \citet{basciftci2023capturing}, the follower
problem has a totally unimodular constraint matrix; it can be solved
as a linear program with an objective of the form $M \ g^r + t^r$ for
a suitably large value of $M$. This
can be achieved by replacing $\tau^r_{hl}$ and $\gamma^r_{ij}$ values
with $\hat{\tau}^r_{hl} := M \tau^r_{hl} + t_{hl}' + t_{hl}^{wait}$
and $\hat{\gamma}^r_{ij} := M \gamma^r_{ij} + t_{ij}$, respectively, where the value of big-M needs to be chosen according to the studied instances. To adjust the objective function value in~\eqref{eq:odmts_da_upperLevelObj}, $\beta_{hl}$ and $\varphi$ values can be further replaced with $\hat{\beta}_{hl} := M \beta_{hl}$ and $\hat{\varphi} := M \varphi$, respectively. Furthermore, to limit the maximum
number of transfers in each suggested path, a constraint can be added
to the follower problem. As the addition of this constraint violates
the totally unimodularity property, the underlying network can be
remodelled through a transfer expanded graph that captures the
transfer limit by design \citep{dalmeijer2020transfer}.
Note that the term \textit{direct shuttle trips} indicates ODMTS trips that are directly served by on-demand shuttles without the involvement of buses or rails.

\subsection{ODMTS-DFD}
Figure~\ref{fig:odmts_dfd_formulation} presents the formulation of {\sc ODMTS-DFD}.
\begin{figure}[!t]
\begin{subequations} 
\label{eq:odmts_dfd_upperLevelProblem}
\begin{alignat}{1}
\min_{z_{hl}} \quad & \sum_{h,l \in H} \beta_{hl} z_{hl} + \sum_{r \in \hat{T}} p^r (g^r - \varphi^r) \label{eq:odmts_dfd_upperLevelObj} \\
\text{s.t.} \quad & \eqref{eq:odmts_da_upperLevelConstr1}, \eqref{eq:odmts_da_fixed_arcs}, \eqref{eq:odmts_da_binaryConstraint} 
\end{alignat}
\end{subequations}
where $(\mathbf{x^r}, \mathbf{y^r}, g^r)$ are a solution to the optimization problem 
\eqref{eq:odmts_da_lowerLevelProblem}.
\caption{The Optimization Model for ODMTS Design with Fixed Demand $\hat{T}$: {\sc ODMTS-DFD}($\hat{T}$).}
\label{fig:odmts_dfd_formulation}
\end{figure}

\subsection{Proof for Proposition 5} 
The term $\text{UB}^r$ for each trip $r$ is defined in the following Proposition, which is proved by~\citet{basciftci2023capturing}.

\textbf{Proposition 7}: 
Consider an ODMTS design $\mathbf{z}^1$, and any ODMTS design
$\mathbf{z}^2 \geq \mathbf{z}^1$. Let the ODMTS travel time for trip
$r$ under $\mathbf{z}^1$ and $\mathbf{z}^2$ be denoted as $t^1$ and
$t^2$. Then $t^2$ is bounded by an upper bound $\text{UB}^r$ defined as follows:
\begin{equation*}
\text{UB}^r = \begin{cases}
    t^1 + \frac{(1 - \theta)}{\theta} \omega \left(d_{or^rm}^{shuttle}+ d_{nde^r}^{shuttle} - \min\limits_{h,l \in H} \left\{d_{or^rh}^{shuttle} + d_{lde^r}^{shuttle}\right\}\right), \\
    \text{if trip $r$ is assigned to a multimodal path under design $\mathbf{z}^1$ such that} \\
    \text{$m$ and $n$ are the first and last hub in the path, respectively.} \\
    \\
    \max \left\{t^1, t^1 + \frac{(1 - \theta)}{\theta} \omega \left(d_{or^r de^r}^{shuttle} - \min
    \limits_{h,l \in H} \left\{d_{or^rh}^{shuttle} + d_{lde^r}^{shuttle}\right\}\right)\right\}, \\
    \text{if trip $r$ is assigned to a direct shuttle path under design $\mathbf{z}^1$} \\
    \end{cases}
\end{equation*}

\textbf{Proposition 5}: The ODMTS design found by Algorithm~\mynameref{alg:arc_s1} with expansion rule~\ref{rule:upper_bound_odmts} satisfies the Correct Adoption property (Property~\ref{property:adopt}), i.e., $\%A_{false} = 0$.

\proof{Proof:} Consider the series of designs $\mathbf{z}^0,
\mathbf{z}^1..., \mathbf{z}^k$ found by Algorithm
\mynameref{alg:arc_s1} using expansion rule~\ref{rule:upper_bound_odmts}; then there are
$\mathbf{z}^0 \leq \mathbf{z}^1 \leq... \leq \mathbf{z}^k$ by
construction. Moreover, $\mathbf{z}^0$ satisfies the {\sc Correct Adoption}
property since there are no latent trips included in $\hat{T} =
\overline{T}^0$. For an arbitrary design $\mathbf{z}^i$ where $0 < i
\leq k$, it follows that all latent trips included in $\hat{T} =
\overline{T}^i = T_{core} \cup \bigg( \bigcup\limits_{j = 0,
  ..., i - 1} {T_{latent}}^{j}_{next} \bigg) $ adopt $\mathbf{z}^i$ by
definition of expansion rule~\ref{rule:upper_bound_odmts} and because trips in
${T_{latent}}^{j}_{next}$ for any $j = 0, 1, ..., i - 1$ adopt $\mathbf{z}^i$,
which is greater than any $\mathbf{z}^j$. Thus, all of these designs
satisfy the {\sc Correct Adoption} property.  \Halmos\endproof
\section{Additional Information on SCTS-DA}
\label{sect:scts_da_appendix}
This section completes Section~\ref{sect:scts_da} by providing the full formulations for {\sc SCTS-DA} and {\sc SCTS-DFD} and their modeling details.

\begin{figure}[ht]
\begin{subequations} \label{eq:scooter_da_upperLevelProblemUpdated2}
\begin{align}
\min_{z_{hl}, \delta^r} \quad & \sum_{h,l \in H} \beta_{hl} z_{hl} + \sum_{r \in T_{core}} p^r (g^r - \varphi^r) + \sum_{r \in T_{latent}} p^r \delta^r (g^r - \varphi^r) \label{eq:scooter_da_upperLevelObj} \\
\text{s.t.} \quad 
& \eqref{eq:odmts_da_upperLevelConstr1}-\eqref{eq:odmts_da_continuousConstraint} \notag
\end{align}
\end{subequations}
where $(\mathbf{x^r}, \mathbf{y^r}, g^r)$ are a solution to the optimization problem
\begin{subequations}
\label{eq:scooter_da_lowerLevelProblem}
\begin{align}
\lexmin_{x_{hl}^r, y_{ij}^r, g^r, t^r, s^r} \quad &  \langle t^r, s^r \rangle \label{eq:scooter_da_lowerLevelObj} \\
\text{s.t.} \quad
   & g^r = \sum_{h,l \in H}  \theta (t_{hl}^{bus/rail} + t_{hl}^{wait}) x_{hl}^r + \sum_{i,j \in N}  \theta t_{ij}^{scooter} y_{ij}^r  \label{eq:scooter_da_time_def} \\
   & t^r = \sum_{h,l \in H}   (t_{hl}^{bus/rail} + t_{hl}^{wait}) x_{hl}^r + \sum_{i,j \in N}  t_{ij}^{scooter} y_{ij}^r \\
   & s^r = \mathbbm{1} \bigg(\sum_{h, l \in H} x_{hl}^r = 0 \bigg) \\
   & s^r \in \{ 0 , 1 \} \\
   & \eqref{eq:odmts_da_minFlowConstraint}-\eqref{eq:odmts_da_integralityFlowConstr} \notag 
\end{align}
\end{subequations}
\caption{The Bilevel Formulation for Scooters-Connected Transit Systems Designs with Adoptions ({\sc SCTS-DA}). }
\label{fig:scts_da_formulation}
\end{figure}
The set $N$ represents the locations where scooters are available, while the set $H \subseteq N$ represents the hubs used for connecting rails and high-frequency buses. Similar to the {\sc ODMTS-DA} problem in Section~\ref{sect:odmts_da}, the focus in {\sc SCTS-DA} is on designing the high-frequency bus network, while the rail system is preserved in the set $Z_{backbone}$. The variables $x_{hl}^r$ represent whether a trip $r$ utilizes the hub-to-hub arc $(h,l)$ in its path, while $y_{ij}^r$ indicate the utilization of the scooter leg $(i,j)$ in the journey.

The bilevel formulation of the {\sc SCTS-DA} problem is presented in Figure~\ref{fig:scts_da_formulation}. Its leader problem is identical to {\sc ODMTS-DA} because the agency also aims to design a transit network. Unlike in {\sc ODMTS-DA}, where the transit agency provides paths to riders considering both operating cost and travel duration, in {\sc SCTS-DA}, riders are assumed to choose their own paths based on the fastest option. Therefore, for each trip $r$, the subproblem is formulated as a shortest path problem based on travel duration $t^r$. Binary variable $s^r$ indicates if riders exclusively use scooters and do not utilize the bus or rail network. To unify the objective in the objective function~\eqref{eq:scooter_da_upperLevelObj}, the cost term $g^r$ is defined as $t^r \theta$. The lexicographic objective~\eqref{eq:scooter_da_lowerLevelObj} prioritizes paths that use only scooters in case of ties in travel time.

\begin{figure}[!t]
\begin{subequations} 
\label{eq:scooter_dfd_upperLevelProblem}
\begin{alignat}{1}
\min_{z_{hl}} \quad & \sum_{h,l \in H} \beta_{hl} z_{hl} + \sum_{r \in \hat{T}} p^r (g^r - \varphi^r) \label{eq:scooter_dfd_upperLevelObj} \\
\text{s.t.} \quad & \eqref{eq:odmts_da_upperLevelConstr1}, \eqref{eq:odmts_da_fixed_arcs}, \eqref{eq:odmts_da_binaryConstraint} 
\end{alignat}
\end{subequations}
where $(\mathbf{x^r}, \mathbf{y^r}, g^r)$ are a solution to the optimization problem 
\eqref{eq:scooter_da_lowerLevelProblem}.
\caption{The Optimization Model for scooter Design with Fixed Demand $\hat{T}$: {\sc SCTS-DFD}($\hat{T}$).}
\label{fig:scooter_dfd_formulation}
\end{figure}

Following the construction of {\sc TN-DFD}, the {\sc SCTS-DFD} problem is presented in Figure~\ref{fig:scooter_dfd_formulation}, and the Evaluation function for design $\mathbf{z}$ is: 
\begin{equation}
    \text{\sc Eval-SCTS}(\mathbf{z}) =  \sum_{h,l \in H} \beta_{hl} z_{hl} + \sum_{r \in T_{core}} p^r g_\pi^r + \sum_{r \in T_{latent}} p^r \delta^r (g_\pi^r - \varphi^r) 
\end{equation}

\section{Additional Information for ODMTS-DA Instance in Ypsilanti, MI}
\label{sect:michigan_instance_appendix}
This section provides the supplementary information and additional computational studies over the medium-sized {\sc ODMTS-DA} instance (see Section~\ref{subsect:ypsi_case_study}) in Ypsilanti, Michigan, USA.

\subsection{Experimental Settings}
\begin{figure}[!ht]
    \centering
    \includegraphics[width=0.5\textwidth]{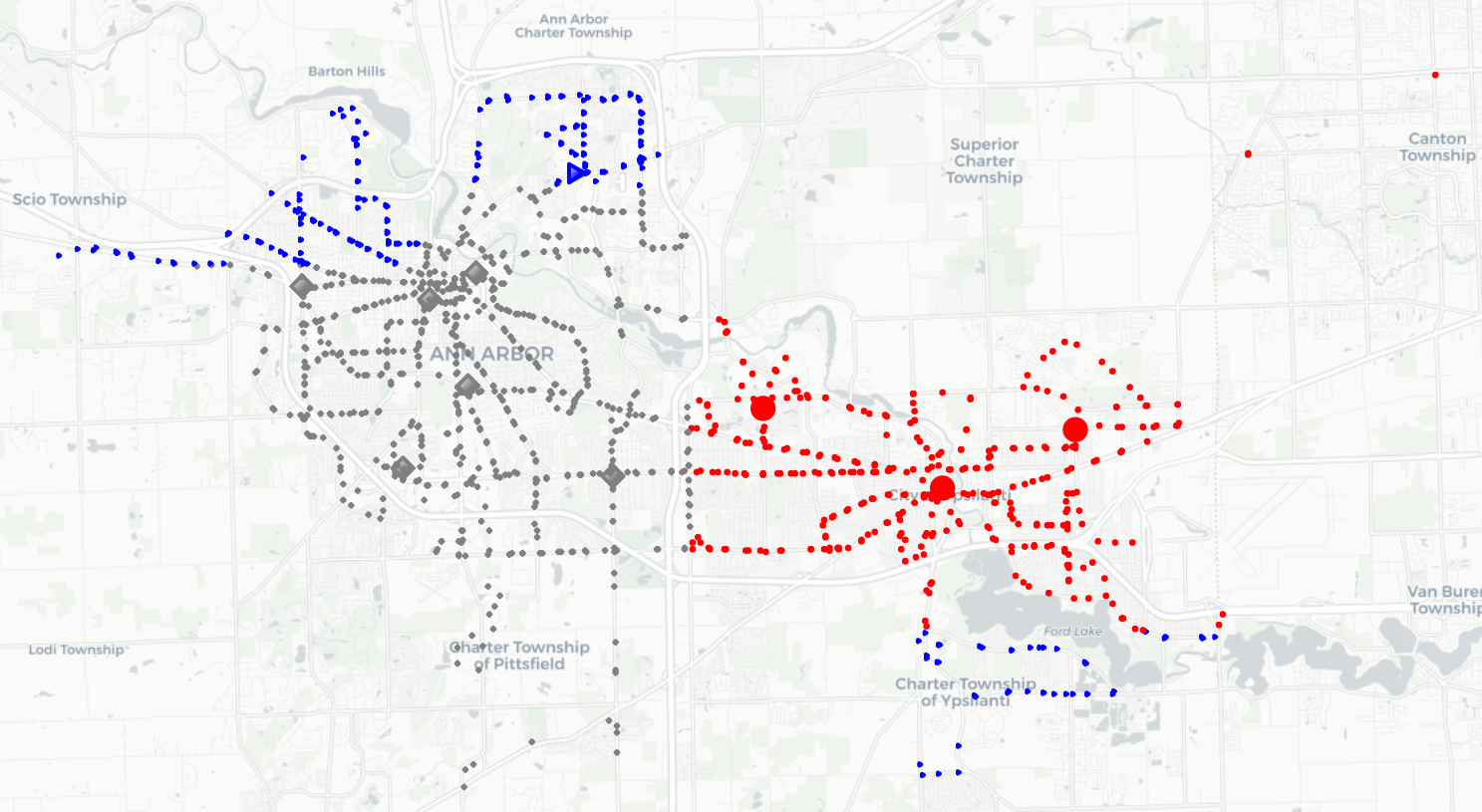}
    \caption{The 1267 bus stops used in this case study. The enlarged 10 stops are the designated hubs. The low-income, middle-income, and high-income stops are represented by red circles, grey squares, and blue triangles respectively. }
    \label{fig:yp_stop_loc}
\end{figure}
This case study is based on the AAATA transit system that operates over 1,267 stops. In order to design an ODMTS, 10 stops are designated as ODMTS hubs, and the other stops are only accessible to on-demand shuttles. There are 1,503 trips for a total of 5,792 riders, and each trip is represented by an origin stop and a destination stop. The time horizon for this case study is between 6 pm and 10 pm. Therefore, this dataset primarily consists of commuting trips from work locations to home. The mode preference of riders depends on their income levels, i.e, a rider with higher income level has lower tolerance on the amount of ODMTS travel time. Thus, the whole 1,503 trips were first divided into three different categories: low-income, middle-income, and high-income, based on the destination of the trips as it is associated with the home location of the riders. Figure~\ref{fig:yp_stop_loc} visualizes the locations of the stops and their corresponding income levels. If a trip is assigned to a certain income-level class, all riders belonging to this trip are then assumed to have the same income-level. Following this procedure, there are 476 low-income, 819 middle-income, and 208 high-income trips with 1,754, 3,316, and 722 riders respectively. Once the classification is carried out, an $\alpha$ value associated to the model {\sc Time-Based} (Equation~\eqref{eq:choice_model_time_based}) is assigned to each class. It must be noted that all trips in the low-income class are treated as members of the core trips set $T_{core}$; hence, no $\alpha$ value is required for them. For the middle-income and the high-income classes, 2.0 and 1.5 are employed as the $\alpha$ values, reflecting an assumption that the medium-income class has higher tolerance for travel time on transit.

In order to evaluate the performances of all five iterative algorithms
with different configurations, two values for each of the two
parameters---rider multiplier and core trips percentage are selected,
resulting in overall four experiments as described in Table~\ref{table:ypsi_exp_setup}. More specifically, the rider multiplier
parameter multiplies the number of riders for each trip. The core
trips percentage parameter is utilized to divide the dataset into core
trips and latent trips by varied partitions. Note that Experiment 4 is already presented in Section~\ref{subsect:ypsi_case_study}.

The on-demand shuttle price is fixed at \$1 per kilometer and each
shuttle is assumed to only serve one passenger. Shuttles between bus
hubs are not allowed in this case study. For buses, the operating fee
is \$3.38 per kilometer. Four buses per hour is the unique bus
frequency considered here; hence, each opened bus leg has an average
waiting time ($t_{hl}^{wait}$ in~\eqref{eq:odmts_da_frDefinition})
of 7.5 minutes. This assumes that passengers arrive at the bus stop randomly rather than attempting to synchronize with the bus schedule. For passengers, there is
no limit on the number of transfers. A constant ticket price that is
consistent with the existing AAATA system is applied. Thus, a full
ODMTS service charges \$2.5 for each individual regardless of the
travel length. Inconvenience is measured in seconds, and the
inconvenience and cost parameter $\theta$ is fixed at 0.001.

The iterative algorithms and the exact algorithm are programmed with
Python 3.7. The exact algorithm is set to terminate if it reaches a 0.1\% optimality gap; otherwise, it reports the upper bound as the ``optimal'' solution after a 6 hour time limit. The difference between the lower bound and upper bound identified by the exact algorithm when it terminates is denoted as `\% Opt. Gap'. The computational
experiments are carried out using Gurobi 9.5 without multi-threading on subproblems. {\sc P-Path} can obtain optimal solution in short amount of time, and it was programmed with Python 3.9 with the support of Gurobi 9.5. For the trip-based iterative algorithms, the $\rho$ value for Algorithm~\mynameref{alg:rho_GRAD} and Algorithm~\mynameref{alg:rho_GAGR}, and the $\eta$ value for Algorithms~\mynameref{alg:eta_GRRE} and~\mynameref{alg:rho_GAGR} are both set to
10. For arc-based iterative algorithms, two experimental runs are designed to test Algorithm~\mynameref{alg:arc_s1} with trip expansion rules~\ref{rule:all_adopt_odmts} and~\ref{rule:upper_bound_odmts}. Algorithm
\mynameref{alg:arc_s2} employs rules~\ref{rule:use_fixed_route_odmts} and~\ref{rule:all_adopt_odmts} and
rules~\ref{rule:upper_bound_odmts} and~\ref{rule:all_adopt_odmts} for the two stages. When solving
the ODMTS-DFD problem at each iteration, all iterative algorithms
employ a stand-alone solving tool developed by
\citet{dalmeijer2020transfer} that leverages Benders decomposition
algorithm. This tool is based on the solver CPLEX 12.9 and allows
multi-threading when solving subproblems in ODMTS-DFD. However, in
this case study, multi-threading techniques are disabled in order to
fairly compare the algorithmic performances.

\begin{table}[!t]
  \centering
\resizebox{0.8\textwidth}{!}{
\begin{tabular}{c c c c c c c c c c c}
\toprule
       & \multirow{3}{*}{ridership} & \multicolumn{3}{c}{low income core trips} & \multicolumn{3}{c}{medium income core trips} & \multicolumn{3}{c}{high income core trips}\\ \cmidrule{3-11}

       &  & \% trips & \# trips &  \# riders & \% trips & \# trips & \# riders & \% trips & \# trips &  \# riders\\ \midrule

Exp.1 & regular & 100\% & 476 & 1754 & 75\% & 614 & 2842 & 50\% & 104 & 434 \\ \midrule

Exp.2 & regular & 100\% & 476 & 1754 & 50\% & 409 & 2262 & 25\% & 52 & 258 \\ \midrule

Exp.3 & doubled & 100\%  & 476 & 3508 & 75\% & 614 & 5684 & 50\% & 104 & 868 \\ \midrule

Exp.4 & doubled & 100\% & 476 & 3508 & 50\% & 409 & 4524 & 25\% & 52 & 516 \\ \bottomrule
\end{tabular}
}
\caption{The  experimental  setups  for  the  four experiments. For doubled ridership, the number of riders for each O-D pair is multiplied by 2. The core trips percentages for each income level are [100\%, 75\%, 50\%] and [100\%, 50\%, 25\%] for low, medium, and high income trips, respectively. Experiment 4 is already presented in Section~\ref{subsect:ypsi_case_study}.}%
\label{table:ypsi_exp_setup}
\end{table}

\subsection{Additional Results for Experiments 1, 2, and 3}

\begin{table}[!ht]
    \centering%
    \resizebox{0.9\textwidth}{!}{

    \begin{tabular}{l r r r r r r r r r r}%
    \toprule
    Method&Step Size&Trip Rule&\# Outer Itr.&\# Total Itr.&Run Time (min)&Design Obj. & \% Opt Gap. &\% Gap &$\%R_{false}$&$\%A_{false}$\\%
    \midrule%
    Exact Alg.& - & - & - & - & 364.10 & 19012.91 & 0.65 & - & 0.00 & 0.00 \\%
    \midrule
    {\sc P-Path} & - & - & - & - & 2.19 & 19012.91 & - & 0.00 & 0.00 & 0.00\\%
    \midrule
    Alg.~\mynameref{alg:rho_GRAD}&10&{-}&1&17&1.66&19012.91&-& 0.00&0.00&0.00\\%
    \midrule%
    Alg.~\mynameref{alg:eta_GRRE}&10&{-}&1&17&1.67&19012.91&-&0.00&0.00&0.00\\%
    \midrule%
    Alg.~\mynameref{alg:rho_GAGR}&10&{-}&17&154&14.52&19012.91&-&0.00&0.00&0.00\\%
    \midrule%
    Alg.~\mynameref{alg:arc_s1}&{-}&\ref{rule:all_adopt_odmts}&1&5&0.37&19012.91&-&0.00&0.00&29.13\\%
    \midrule%
    Alg.~\mynameref{alg:arc_s1}&{-}&\ref{rule:upper_bound_odmts}&1&5&0.38&19012.91&-&0.00&19.42&0.00\\%
    \midrule%
    Alg.~\mynameref{alg:arc_s2}&{-}&\ref{rule:use_fixed_route_odmts},~\ref{rule:all_adopt_odmts}&1&6&0.43&19012.91&-&0.00&0.00&2.27\\%
    \midrule%
    Alg.~\mynameref{alg:arc_s2}&{-}&\ref{rule:upper_bound_odmts},~\ref{rule:all_adopt_odmts}&1&6&0.43&19012.91&-&0.00&0.00&0.00\\%
    \bottomrule%
    \end{tabular}%
    }
\caption{Performance comparison on Experiment 1. The exact algorithm reports optimality gap after 6 hours. The heuristics report the gap between themselves and the optimal solution ({\sc P-Path}).}
\label{tab:yp_exp_1_alg_perform}
\end{table}
\begin{table}[!ht]
    \centering%
    \resizebox{0.9\textwidth}{!}{
    
    \begin{tabular}{l r r r r r r r r r r}%
    \toprule
    Method&Step Size&Trip Rule&\# Outer Itr.&\# Total Itr.&Run Time (min)&Design Obj. & \% Opt Gap. &\% Gap &$\%R_{false}$&$\%A_{false}$\\%
    \midrule%
    Exact Alg.& - & -& - & - & 367.52 & 16635.73 & 2.48 &- & 0.00 & 0.00\\%
    \midrule
    {\sc P-Path} & - & -& - & - & 7.80 & 16635.73 & - & 0.00 & 0.00 & 0.00\\%
    \midrule
    Alg.~\mynameref{alg:rho_GRAD}&10&{-}&1&33&4.28&16635.73&-&0.00&0.00&5.12\\%
    \midrule%
    Alg.~\mynameref{alg:eta_GRRE}&10&{-}&1&30&4.18&16635.73&-&0.00&12.90&5.12\\%
    \midrule%
    Alg.~\mynameref{alg:rho_GAGR}&10&{-}&30&480&54.78&16635.73&-&0.00&0.00&5.30\\%
    \midrule%
    Alg.~\mynameref{alg:arc_s1}&{-}&\ref{rule:all_adopt_odmts}&1&5&0.59&16635.73&-&0.00&0.00&26.50\\%
    \midrule%
    Alg.~\mynameref{alg:arc_s1}&{-}&\ref{rule:upper_bound_odmts}&1&4&0.41&16983.87&-&2.09&33.57&0.00\\%
    \midrule%
    Alg.~\mynameref{alg:arc_s2}&{-}&\ref{rule:use_fixed_route_odmts},~\ref{rule:all_adopt_odmts}&1&6&0.76&16635.73&-&0.00&0.00&14.66\\%
    \midrule%
    Alg.~\mynameref{alg:arc_s2}&{-}&\ref{rule:upper_bound_odmts},~\ref{rule:all_adopt_odmts}&1&6&0.61&16635.73&-&0.00&0.00&13.60\\%
    \bottomrule%
    \end{tabular}%
}
\caption{Performance comparison on Experiment 2. The exact algorithm reports optimality gap after 6 hours. The heuristics report the gap between themselves and the optimal solution ({\sc P-Path}).}
    \label{tab:yp_exp_2_alg_perform}
\end{table}
Table~\ref{tab:yp_exp_1_alg_perform} presents an overview of the algorithmic performances of all algorithms on Experiment 1. All heuristic algorithms achieve the optimal solution in a short duration. Table~\ref{tab:yp_exp_2_alg_perform} compares the computational results obtained from all algorithms on Experiment 2. Similar to Experiment 4 (see Section~\ref{subsect:ypsi_case_study}), the ratio of core trips to the set of all trips is significantly lower in this experiment; hence, more iterations are conducted in most of the algorithms. Furthermore, $\%R_{false}$ and $\%A_{false}$ are not necessarily equal to 0 although the optimal design $\mathbf{z}^*$ is found by the iterative algorithms. This result can be explained by the statement summarized in Remark~\ref{rmk:opt_but_not_0p}. Under this setting, most of the iterative algorithms successfully obtain the optimal design found by the exact algorithm and {\sc P-Path}, except Algorithm~\mynameref{alg:arc_s1} with rule~\ref{rule:upper_bound_odmts}. This exception could be attributed to the design of rule~\ref{rule:upper_bound_odmts} which strictly selects  trips for the following iterations, resulting in an earlier termination. The algorithmic performances of Experiment 3 are summarized in Table~\ref{tab:yp_exp_3_alg_perform}. Under this setting, none of the proposed algorithms successfully discover the optimal solution. The trip set $\hat{T}^*$ constructed by the optimal solution $\mathbf{z}^*$ includes a few trips that reject all ODMTS designs explored by the heuristic algorithms. For trip-based algorithms, this can be a result of the expansion rules as all of them only select candidates from adopting trips at certain iterations, and it is difficult to justify the addition of currently rejected trips. For arc-based algorithms, these results are likely to be related to their greedy searching rules. The arc-based algorithms only fix the cycle that provides the greatest design objective decrements, leading to local, but not global, optimality. Although the proposed algorithms may not discover the optimal solution, these observations may support the hypothesis that the proposed iterative algorithms can find a high-quality solution in short duration. Indeed, the reported designs only have a 0.03\% optimality gap.  Furthermore, some iterative algorithms achieve zero {\sc False Rejection Rate} and zero {\sc False Adoption Rate} while not obtaining the optimal solution, as discussed in Remark~\ref{rmk:0p_but_not_opt}.

\begin{table}[!ht]
    \centering%
    \resizebox{0.9\textwidth}{!}{
    \begin{tabular}{l r r r r r r r r r r}%
    \toprule
    Method&Step Size&Trip Rule&\# Outer Itr.&\# Total Itr.&Run Time (min)&Design Obj. & \% Opt Gap. &\% Gap &$\%R_{false}$&$\%A_{false}$\\%
    \midrule%
    Exact Alg.& - & - & - & - & 363.22 & 34732.09 & 0.99 & - & 0.00 & 0.00 \\%
    \midrule
    {\sc P-Path} & - & - & - & - & 2.62 & 34732.09 & - & 0.00 & 0.00 & 0.00 \\%
    \midrule
    Alg.~\mynameref{alg:rho_GRAD}&10&{-}&1&15&2.12&34743.29&-&0.03&0.00&0.00\\%
    \midrule%
    Alg.~\mynameref{alg:eta_GRRE}&10&{-}&1&15&2.15&34743.29&-&0.03&0.00&0.00\\%
    \midrule%
    Alg.~\mynameref{alg:rho_GAGR}&10&{-}&15&121&17.10&34743.29&-&0.03&0.00&0.00\\%
    \midrule%
    Alg.~\mynameref{alg:arc_s1}&{-}&\ref{rule:all_adopt_odmts}&1&8&0.91&34743.29&-&0.03&0.00&42.39\\%
    \midrule%
    Alg.~\mynameref{alg:arc_s1}&{-}&\ref{rule:upper_bound_odmts}&1&8&0.80&34743.29&-&0.03&12.94&0.00\\%
    \midrule%
    Alg.~\mynameref{alg:arc_s2}&{-}&\ref{rule:use_fixed_route_odmts},~\ref{rule:all_adopt_odmts}&1&9&0.86&34743.29&-&0.03&0.00&5.50\\%
    \midrule%
    Alg.~\mynameref{alg:arc_s2}&{-}&\ref{rule:upper_bound_odmts},~\ref{rule:all_adopt_odmts}&1&9&0.86&34743.29&-&0.03&0.00&0.00\\%
    \bottomrule%
    \end{tabular}%
    }
\caption{Performance comparison on Experiment 3. The exact algorithm reports optimality gap after 6 hours. The heuristics report the gap between themselves and the optimal solution ({\sc P-Path}).}
\label{tab:yp_exp_3_alg_perform}
\end{table}

\section{Data and Experimental Settings for Atlanta Case Studies}
\label{sect:atlanta_data_appendix}

This section introduces the procedures to construct the dataset used in the large-scale case studies, both {\sc ODMTS-DA} and {\sc SCTS-DA}. This dataset focuses on the 6 am--10 am morning peak of a regular weekday. ODMTS or SCTS operate across 2,426 locations. For hubs, 58 and 66 are selected in different case studies: ODMTS-DA large instance and SCTS-DA have 58 hubs, while ODMTS-DA extra-large instance has 66 hubs. Among the hubs, 38 of them are adopted from rail stations. There are 15,478 core trips, 36,283 latent trips with 55,871 existing riders, and 54,902 riders with choices. As morning rush hours are considered, congestion is considered in computing the travel times of on-demand shuttles in case studies.

\subsection{Core Demand}
The core demand of this case study is based on two realistic datasets provided by the transit agency \textit{Metropolitan Atlanta Rapid Transit Authority} (MARTA). One dataset includes the transit pass transaction records collected by the Automated Fare Collection (AFC) system, and the other dataset is constructed using Automated Passenger Counter (APC) system in the MARTA buses. These two datasets are combined to generate a set of trips that represent a regular weekday in March 2018. For practical reasons, all MARTA stops within a \SI{1500}{feet} (\SI{457.2}{\meter}) radius are first clustered into 1,585 \textit{cluster stops}. The clustering process initially guarantees the reservation of all rail stations as cluster stops. Subsequently, it progressively incorporates individual bus stops into clusters or designates them as new clusters. In cases where a bus stop is eligible for inclusion in multiple clusters, it is allocated to the closest cluster regarding to the haversine distance. The core trip set $T_{core}$ is then constructed after narrowing the time window to the morning peak and applying trip chaining techniques, which are identical to the data employed by \citet{auad2021resiliency} for the baseline scenarios.

\subsection{Latent Demand}
\label{subsect:atlanta_latent_demand}
This section introduces the source of the latent demand for this case study. The ABM-ARC dataset is simulated by the Atlanta Regional Commission (ARC) based on the ARC Regional Household Travel Survey conducted in 2011 following a
rigorous activity-based model (ABM) (\url{http://abmfiles.atlantaregional.com/} Last Visited Date: September 30, 2025). In this dataset, there are more than 4.5 million simulated agents who are characterized by demographic data corresponding to themselves and their households. It also contains more than 13 million individual tours that are completed in the Atlanta metropolitan area in a typical weekday. All tours are described by detailed travel information including origin-destination (O-D) locations, travel period, travel distance, travel mode, and purposes. In particular, the data with horizon year 2020 is used for this case study to align with the 2018 core demand.

The latent demand for this case study consists of all individual tours that satisfy the following five conditions: (i) drive alone tours,  (ii) commuting tours, (iii) the tours start between 6 am and 10 am, (iv) tour agents whose households are located inside the \textit{East Atlanta} region---a large and generally underserved residential area, and (v) tour agents whose work locations located inside or around the I-285 interstate highway loop or inside the city of Sandy Springs---both are MARTA's primary service regions. The first two conditions regulate tour mode and tour purposes such that commuters who drive alone form the potential market of an ODMTS service. Particularly, the second condition regulates tour purpose because a commuter's household and work location are not regularly changed; hence, these tours are more representative of a regular weekday. The third condition further refines the dataset by narrowing down the time window to the morning peak. The last two conditions control the origin and destination locations. It is worth pointing out that this study keeps the tour destination rather than the first trip-leg destination in order to generalize the home-to-work demand for people residing in East Atlanta. Once the valid tours are isolated from the ABM-ARC dataset, their tour origins and tour destinations are then extracted to form the latent demand for this case study.

All the geographical information in the ABM-ARC dataset are at Traffic Analysis Zone (TAZ) level; however, instead of directly utilizing the TAZ centroids to represent the exact locations, origins and destinations are further sampled to lower geographical levels to construct a more representative dataset. In particular, a TAZ can be decomposed into multiple Census Blocks that belong to or intersect with it. A Census Block consists of one or several adjacent neighbourhoods and records the number of residents. Since all selected trips originate from the commuters' home locations, each origin is sampled from its corresponding TAZ to a Census Block by the population distribution of the census blocks. Similarly, the destinations are sampled from TAZ level to Points of Interest (POI) level. Points of interest are specific locations such as restaurants, schools, and  factories. POIs are associated with a TAZ and their area size is also recorded; thus, in this study, each destination is sampled from its corresponding TAZ to a POI by the area distribution of the POIs.

\subsection{Available Locations, Transit Hubs,  O-D Pairs, and Congestion Modeling}

\begin{figure}[!ht]
    \centering
    \begin{subfigure}[b]{0.19\textwidth}
        \includegraphics[width=\textwidth]{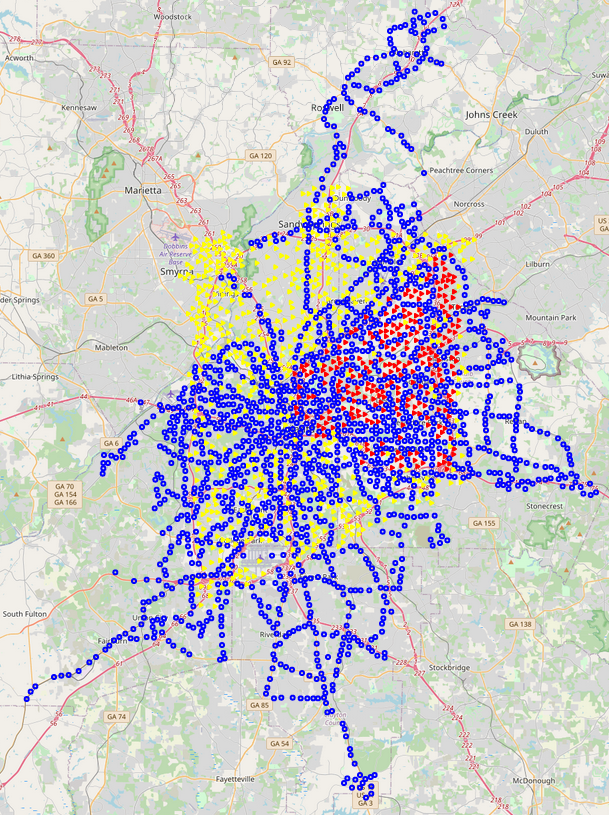}
        \caption{All Locations (set $N$)}
        \label{subfig:atl_stops}
    \end{subfigure}
    \begin{subfigure}[b]{0.19\textwidth}
        \includegraphics[width=\textwidth]{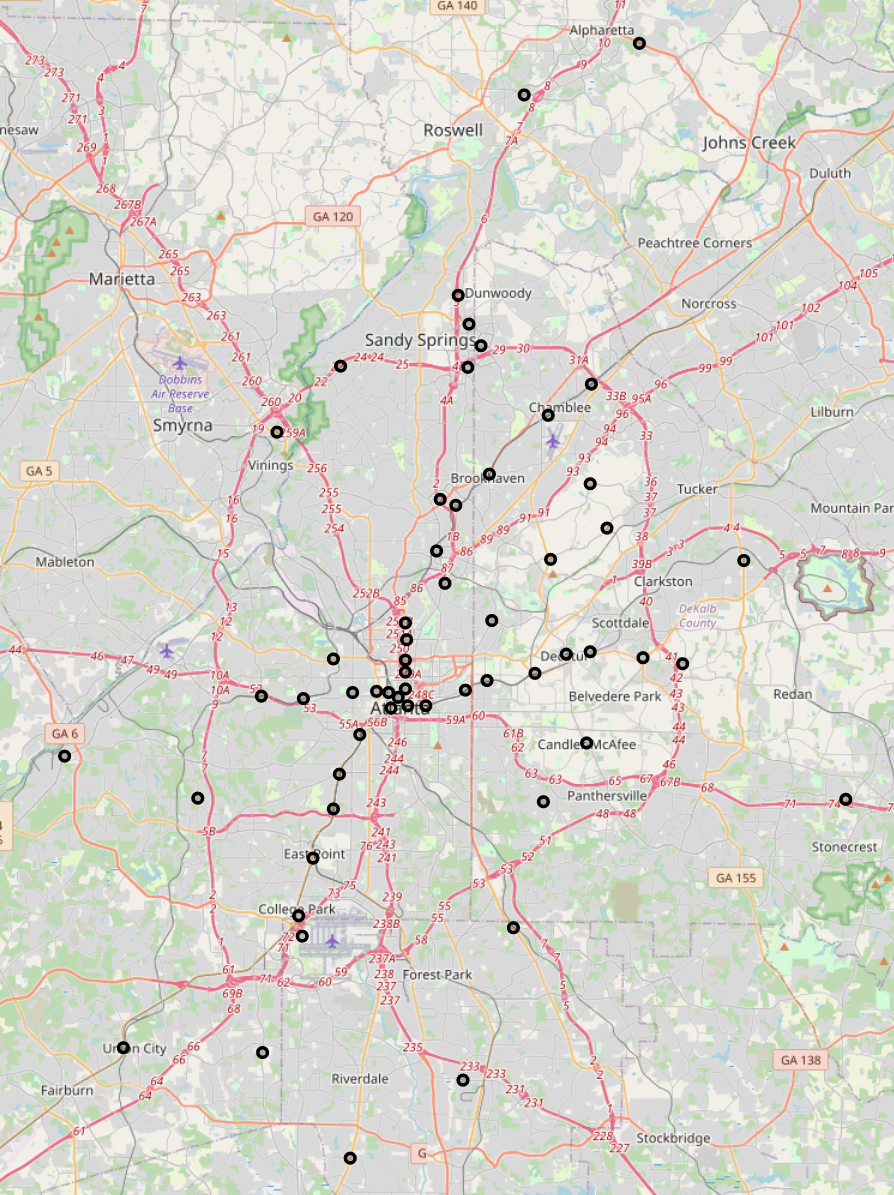}
        \caption{58 Transit Hubs (set $H$)}
        \label{subfig:atl_hubs_58}
    \end{subfigure}
    \begin{subfigure}[b]{0.19\textwidth}
        \includegraphics[width=\textwidth]{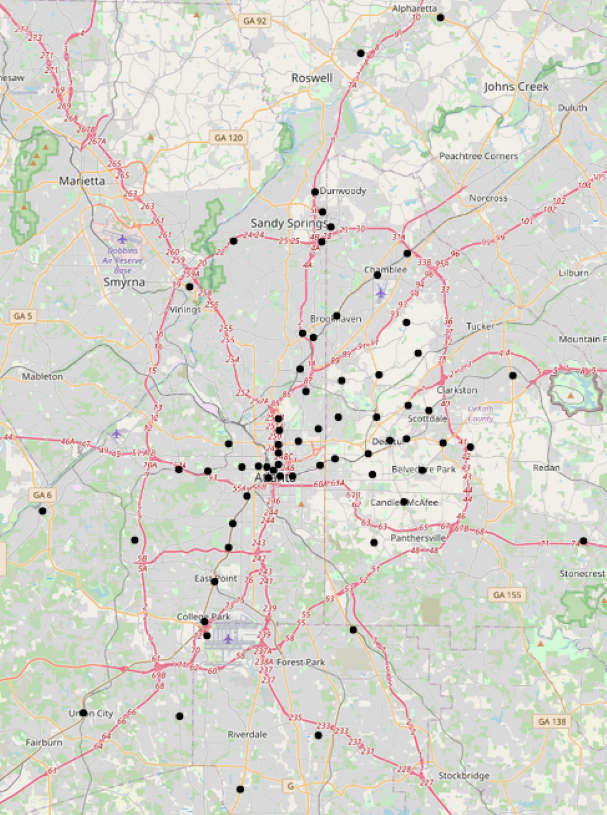}
        \caption{66 Transit Hubs (set $H$)}
        \label{subfig:atl_hubs_66}
    \end{subfigure}
    \begin{subfigure}[b]{0.19\textwidth}
        \includegraphics[width=\textwidth]{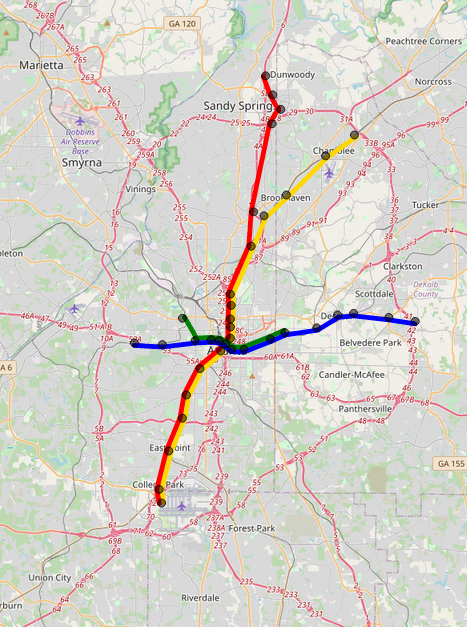}
        \caption{MARTA Rail Lines}
        \label{subfig:atl_rails}
    \end{subfigure}
    \begin{subfigure}[b]{0.19\textwidth}
        \includegraphics[width=\textwidth]{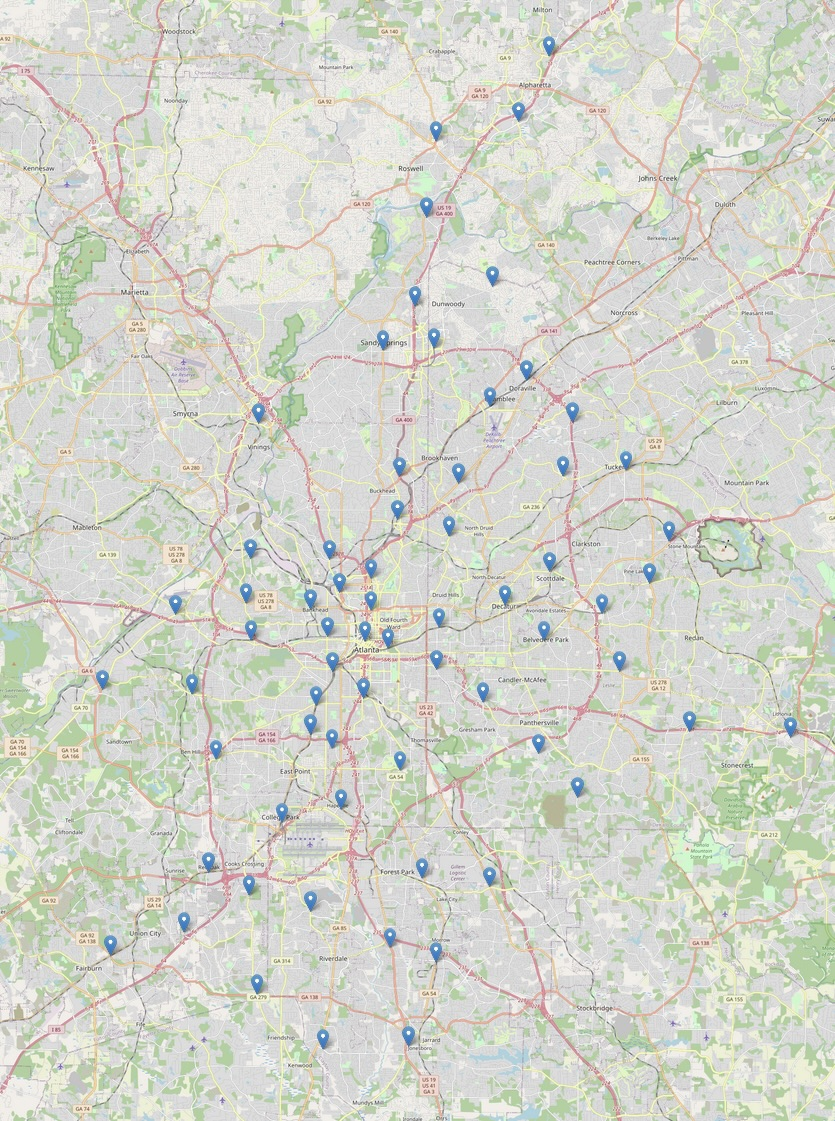}
        \caption{Congestion Query Point}
        \label{subfig:atl_congestion_point}
    \end{subfigure}

\caption{The 2,426 locations are shown in (a). The blue circles, red triangles, and yellow triangles denote the locations of the MARTA cluster stops, stops adopted from Census Block centers, and stops adopted from Points of Interest centers, respectively. The ODMTS hubs and the MARTA rail system are visualized in (b), (c), and (d) respectively. All 38 MARTA rail stations are reserved as ODMTS hubs. The 68 congestion query points are shown in (e). }
\end{figure}
\begin{figure}[htbp]
    \centering
    \begin{subfigure}[b]{0.4\textwidth}
        \includegraphics[width=\textwidth]{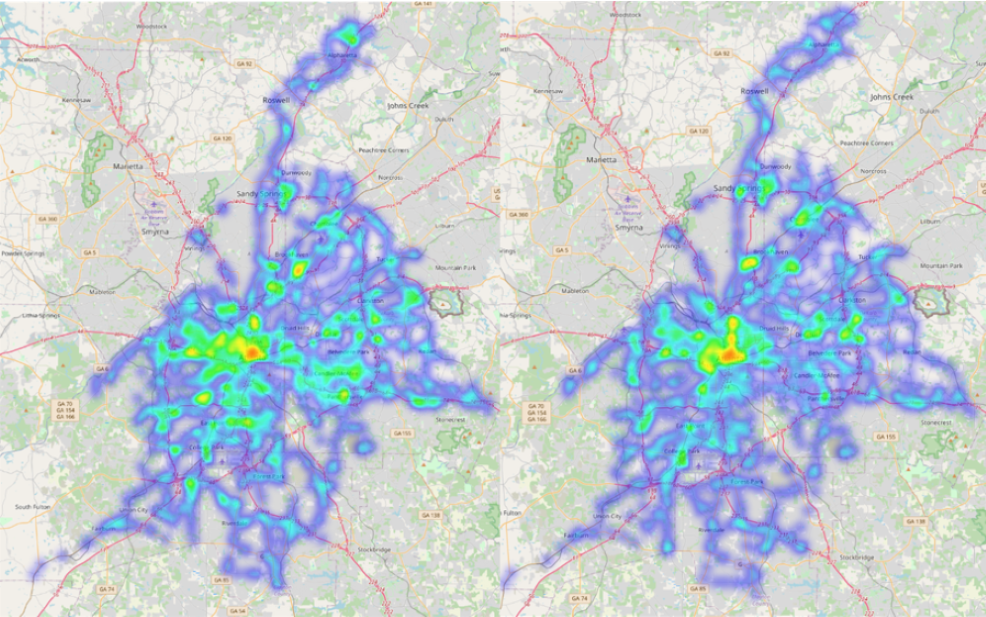}
        \caption{O-Ds in Core Trip Set $T_{core}$}
        \label{subfig:atl_core_trips_heatmap}
    \end{subfigure}
    \begin{subfigure}[b]{0.4\textwidth}
        \includegraphics[width=\textwidth]{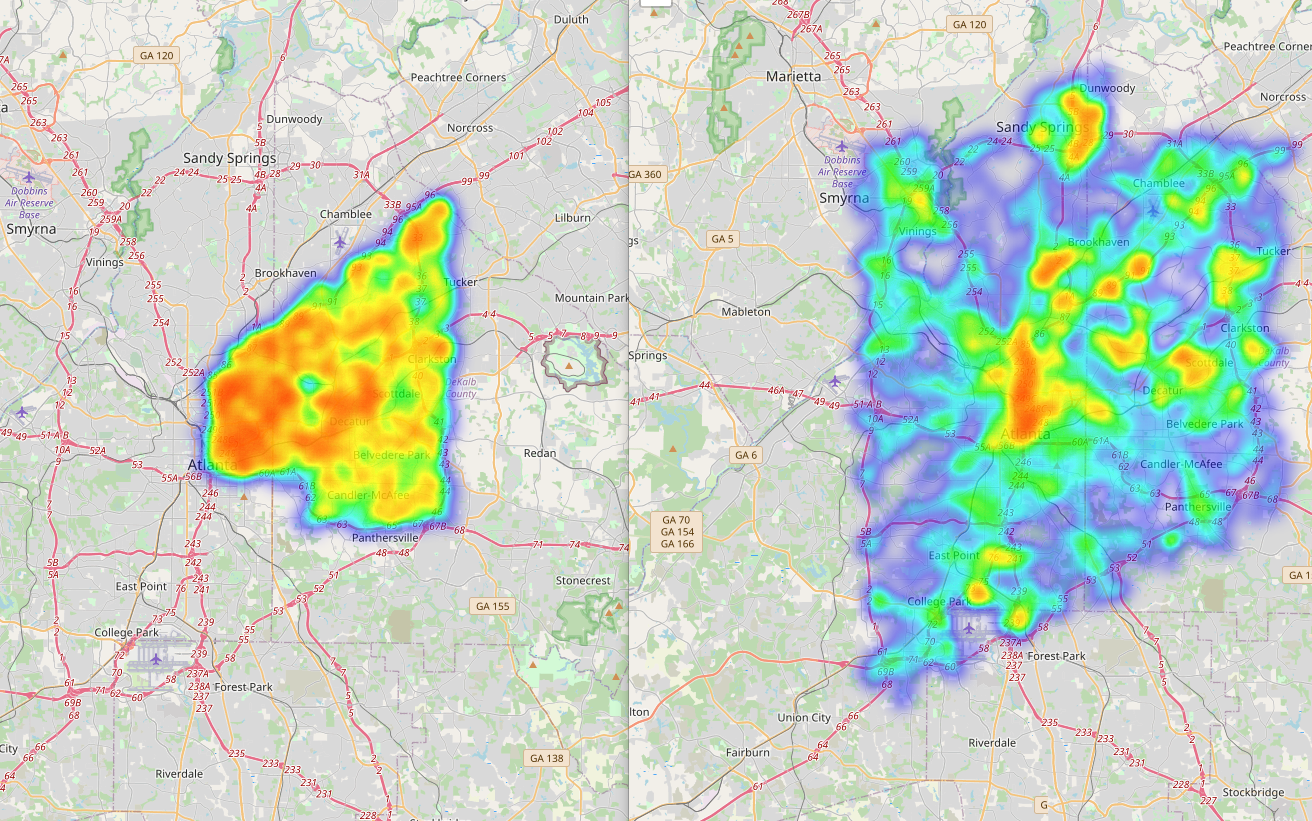}
        \caption{O-Ds in Latent Trip Set $T_{latent}$}
        \label{subfig:atl_additional_trips_heatmap}
    \end{subfigure}
\caption{The heat-maps in (a) visualize the origin and destination stops of the core trip set $T_{core}$. All origins and destinations are MARTA cluster stops. The heat-maps in (b) visualize the origin and destination stops of the Latent Trip set $T_{latent}$. All origins are from the East Atlanta area, and all destinations are either inside or around the I-285 interstate or inside the city of Sandy Springs. The origins and destinations are first sampled to lower geographical level, then mapped to their nearest ODMTS stops.}
\end{figure}
The 2,426 stops (set $N$) are allocated on the basis of three sources: (i) the 1,585 MARTA cluster stops, (ii) the centers of Census Blocks in East Atlanta, (iii) the centers of POIs inside or around the I-285 interstate highway loop or inside the city of Sandy Springs. Specifically, all 1,585 MARTA cluster stops are first reserved as the ODMTS stops. To improve the granularity of stops such that the maximum haversine distance between a stop and its nearest stop is \SI{457.2}{\meter} (\SI{1500}{feet}),  247 Census Block centers and 594 POI centers are then appointed as the ODMTS stops. All accessible locations are illustrated in Figure~\ref{subfig:atl_stops}.

Following the allocation of all stops, for each trip in the latent demand, its origin and destination are mapped to their nearest location with regards to haversine distance. At this stage, the Latent Trip set $T_{latent}$ is constructed. It should be recalled that the trips in the core trip set $T_{core}$ originate from and arrive at MARTA cluster stops, which are already reserved in set $N$. The origins and destinations of $T_{core}$ and $T_{latent}$ are illustrated as heat-maps in Figure \ref{subfig:atl_core_trips_heatmap} and Figure \ref{subfig:atl_additional_trips_heatmap}, respectively.

As shown in Figure~\ref{subfig:atl_hubs_58}~and~\ref{subfig:atl_hubs_66}, 58 and 66 stops are designated as transit hubs, used for different case studies. Among these hubs, 38 of them are rail stations and the rail system is illustrated in \ref{subfig:atl_rails}. The other 28 stops are appointed as hubs following these four rules: (1) the stops are on major roads, (2) there are a considerable amount of riders who depart from or arrive at this stop, (3) the haversine distance between two hubs is at least \SI{1}{km}, and (4)
a higher density of hubs should be allocated in the East Atlanta region.

To estimate congested shuttle travel time during the morning peak, the most convenient method would be to query the Google Maps Directions API (\url{https://developers.google.com/maps/documentation/directions/overview} Last Visited Date: September 30, 2025). However, due to the high cost associated with this approach, this study follows the methodology proposed by \citet{lu2024impact} using query points as illustrated in Figure~\ref{subfig:atl_congestion_point}. The 68 query points include all rail stations and have a significant overlap with the transit hubs. To estimate the congested travel time between two stops, $i$ and $j$, Equation \eqref{eq:scale_congestion} is utilized. In this equation, $q_i$ and $q_j$ correspond to the nearest query points to stops $i$ and $j$ in terms of haversine distance, respectively. If stops $i$ and $j$ are mapped to the same query point, then $q_i$ and its nearest query point are selected. The congested travel time between any two query points is obtained by querying the Google Map Directions API, while the free flow travel time and travel distance between any two stops are queried from GraphHopper (\url{https://www.graphhopper.com/} Last Visited Date: September 30, 2025), which is constructed based on OpenStreetMap data. Here, the same travel time and distance are applied to both buses and shuttles. Since buses run only on the main corridors (typically major roads or multi-lane highways), the difference in travel time between the two modes is relatively small, unlike on smaller local roads.

\begin{equation}
\label{eq:scale_congestion}
t_{i,j}^{shuttle} = t_{i,j}^{shuttle, free} \times \cfrac{t_{q_i,q_j}^{shuttle, congested}}{t_{q_i,q_j}^{shuttle, free}}
\end{equation}

\subsection{Experimental Setting for ODMTS-DA in Atlanta}
\label{subsubsect:atl_exp_setting}
This section presents the experimental settings for the two {\sc ODMTS-DA} case studies in Atlanta. The MARTA rail system, which is the rail system of the Atlanta Metropolitan area, is completely preserved in this experiment and the
passenger travel time between any two rail stations is derived from
GTFS data. Therefore, these decision variables form the set $Z_{backbone}$ in {\sc ODMTS-DA}
(Figure~\ref{fig:odmts_da_formulation}) and  {\sc ODMTS-DFD} (Figure
\ref{fig:odmts_dfd_formulation}). The unique bus frequency considered in this
experiment is six buses per hour; hence, the average waiting time
($t_{hl}^{wait}$ in~\eqref{eq:odmts_da_drDefinition})
is five
minutes. Unlike the Ypsilanti case study, shuttles are allowed to move between hubs in this case study.
Additionally, bus arcs are restricted from overlapping rail arcs, and for each hub, a bus arc can only connect the hub to its three closest rail stations in terms of travel time to simplify the problem setting. For the extra-large instance, this number is increased to 10 nearby rail stations and 10 nearby bus-only hubs (i.e., a hub but not a rail station), to add complexity to the problem. For the adoption models {\sc Time-Based} and {\sc Time-Transfer-Based}, $\alpha = 1.5$ is applied to all latent trips. The transfer tolerance term $l_{ub}^r$ in {\sc Time-Transfer-Based} is set to 3 for all latent trips. A \$2.5 ticket is charged for each ODMTS rider and this constant price
is identical to the current service provided by MARTA. The operating
fees for on-demand shuttle and buses are fixed at \$0.62 per km (\$1
per mile) and \$72.15 per hour, respectively. The inconvenience and
cost parameter $\theta$ is at $7.25 / (60 + 7.25)$ and 0.08 for the large and extra-large instances, respectively. This case
study employs minute as the unit of time in computation. These
parameter values are directly adopted from an Atlanta-based ODMTS
study presented in \citet{auad2021resiliency}. Congested travel time is only applied to the extra-large case study.

For {\sc P-Path}, it ran until it reaches the optimal solution. For the exact algorithm, the running time is limited at 24 hours. Due
to the complexity of the master problem in the {\sc ODMTS-DA} problem
(Figure~\ref{fig:odmts_da_formulation}), for each iteration of the exact algorithm,
the master problem terminates when reaching to 1\% MIP gap or an 8
hours time limit. For the trip-based iterative algorithms, the
iterative step sizes $\rho$ and $\eta$ are tested with three
values---1000, 2000, and 3000. The step size that leads to the best
design objective is reported together with the best discovered ODMTS
design. In particular, Algorithm~\mynameref{alg:rho_GRAD} and Algorithm
\mynameref{alg:eta_GRRE} only demand the $\rho$ value and the $\eta$
value respectively, and $\eta$ is set to equal to $\rho$ for Algorithm
\mynameref{alg:rho_GAGR}. The combined trip-based iterative algorithm,
i.e., Algorithm~\mynameref{alg:rho_GAGR}, is also terminated after 24
hours of running time. For the arc-based iterative algorithms,
expansion rule~\ref{rule:all_adopt_odmts} and~\ref{rule:upper_bound_odmts} are separately utilized
by Algorithm~\mynameref{alg:arc_s1}. The two-stage Algorithm
\mynameref{alg:arc_s2} employs two groups of rules---(i) rules
\ref{rule:use_fixed_route_odmts} and~\ref{rule:all_adopt_odmts} and (ii) rules~\ref{rule:upper_bound_odmts} and
\ref{rule:all_adopt_odmts}. Furthermore, at the starting point, $\mathbf{z}_{fixed}$
represents the MARTA rail system. Lastly, the same set of
computational tools and solvers are used as in the Ypsilanti case study (see Appendix~\ref{sect:michigan_instance_appendix}).

\subsection{Experimental Setting for SCTS-DA in Atlanta}
\label{subsect:scts-da-settings}

Since the transit agency will not receive any fare if riders only use the scooters, and the scooters are relatively slow, the theta value is chosen to be 0.6, to encourage riders to use high-frequency services such as the rail. A \$2.5 constant fare is incurred to the transit agency for each rider. For each hub, similar to the large-scale {\sc ODMTS-DA} instance, bus connections are restricted to the three nearest rail stations based on travel time. Due to variations in time matrices used across case studies, the network size ($|Z|$) for this case study is 878.

The iterative step sizes $\rho$ and $\eta$ are tested as 500. The step size that leads to the best design objective is reported together with the best discovered SCTS design. In particular, Algorithm~\mynameref{alg:rho_GRAD} and Algorithm~\mynameref{alg:eta_GRRE} only demand the $\rho$ value and the $\eta$ value respectively, and $\eta$ is set to equal to $\rho$ for Algorithm~\mynameref{alg:rho_GAGR}.

For the arc-based heuristic algorithms, to expand the network design over iterations, the type of sub-network considered by the algorithms is the cycles formed by bus arcs. Expansion rule~\ref{rule:all_adopt_scts}~and~\ref{rule:backbone_scts} are separately utilized by Algorithm~\mynameref{alg:arc_s1}. The two-stage Algorithm \mynameref{alg:arc_s2} employs a group of rules---rules~\ref{rule:backbone_scts}~and~\ref{rule:all_adopt_scts}. Furthermore, at the starting point, $\mathbf{z}_{fixed}$ represents the MARTA rail system, which is already included in the set $Z_{backbone}$.
\section{The Benefits of ODMTS-DA}
\label{sect:atlanta_real_study_appendix}
In this section, the benefits of the ODMTS-DA are presented under realistic assumptions, providing a complementary perspective to the content covered in Section~\ref{subsect:atlanta_extra_large_instance}.

\subsection{Latent Demand Adoption}

\begin{table}[!ht]
\centering%
\resizebox{0.8\textwidth}{!}{
\begin{tabular}{l r r l r l r l r l }%
\toprule
Algorithm&Trip Rule&\multicolumn{2}{c}{Adoption}&\multicolumn{2}{c}{Only Use Shuttle}&\multicolumn{2}{c}{Use Bus / Rail}&\multicolumn{2}{c}{Profitable Adoption}\\%
\midrule%
Exact Alg.&{-}&31421&57.23\%&15920&50.67\%&15501&49.33\%&19919&63.39\%\\%
\midrule%
Alg.~\mynameref{alg:rho_GRAD}&{-}&30657&55.84\%&15228&49.67\%&15429&50.33\%&20617&67.25\%\\%
\midrule%
Alg.~\mynameref{alg:eta_GRRE}&{-}&33080&60.25\%&16481&49.82\%&16599&50.18\%&21338&64.50\%\\%
\midrule%
Alg.~\mynameref{alg:rho_GAGR}&{-}&30674&55.87\%&15456&50.39\%&15218&49.61\%&20507&66.85\%\\%
\midrule%
Alg.~\mynameref{alg:arc_s1}&\ref{rule:all_adopt_odmts}&30487&55.53\%&15286&50.14\%&15201&49.86\%&20643&67.71\%\\%
\midrule%
Alg.~\mynameref{alg:arc_s1}&\ref{rule:upper_bound_odmts}&39460&71.87\%&20784&52.67\%&18676&47.33\%&20408&51.72\%\\%
\midrule%
Alg.~\mynameref{alg:arc_s2}&\ref{rule:use_fixed_route_odmts}, \ref{rule:all_adopt_odmts}&31043&56.54\%&15235&49.08\%&15808&50.92\%&20990&67.62\%\\%
\midrule%
Alg.~\mynameref{alg:arc_s2}&\ref{rule:upper_bound_odmts}, \ref{rule:all_adopt_odmts}&30187&54.98\%&14758&48.89\%&15429&51.11\%&20779&68.83\%\\%
\bottomrule%
\end{tabular}%
}
\caption{Adoption results of the reported ODMTS designs. To compute \% only use shuttle, \% Use Bus or Rail, and \% Profitable adoptions,  the denominator is \# adopted riders.}%
\label{tab:atl_adoption}
\end{table}

Table~\ref{tab:atl_adoption} summarizes the adoption results of the
ODMTS designs. The first observation is that all ODMTS designs provide
at least a 54.98\% adoption rate, which is more than half of the riders with choice. Among these adopted riders, significant percentage
of them only use the on-demand shuttle services. That is, the
on-demand shuttles pick up the rider at the origin stop and drops her
at the exact destination stop.  There are several possible
explanations for this result. First, there are many local trips inside
the study area; thus, a direct on-demand shuttle is the best solution
when considering both of cost and convenience. A rider always adopts
the direct on-demand shuttle service since its travel time is assumed
to be identical to the direct driving time. Secondly, the $\theta$
value might exaggerate the importance of cost. A rider is then
assigned to an ODMTS path with relatively high inconvenience,
resulting in rejection. Lastly, the choice function used in this study
might be too strict for some riders. For example, a rider rejects a
12 minutes long ODMTS service when the direct driving time is
only 5 minutes. Note that this challenge can be overcome by
considering more flexible choice functions. Although the proportion of direct shuttle users is high, there is
still a massive adoption value in terms of bus or rail
users. Recalling that there are 55,871 existing riders in the system;
thus, these newly adopted bus riders and rail riders can solely
increase the number of public transit users by
27.2\%--33.4\%. These outcomes can be ascribed to the realistic assumptions on congested travel time. Given that major corridors experience congestion during peak hours, ODMTS emerges as a more favorable commuting alternative. It guides individuals toward high-frequency fixed routes, particularly the rail systems, which remain slightly affected or totally unaffected by congestion. Furthermore, the shuttle aspect of ODMTS trips is usually marginally impacted by congestion, primarily because they involve local trip legs.

Inspections on both Table~\ref{table:atl_alg_perform_extra_large} (see Section~\ref{subsect:atlanta_extra_large_instance})
and Table~\ref{tab:atl_adoption} show that the design with the
smallest objective generally has a lower additional demand adoption
rate. This might be explained by the nature of {\sc ODMTS-DFD} (Figure
\ref{fig:odmts_dfd_formulation}) since it tries to provide an optimized network
with the input trip set $\hat{T}$ regardless of the existence of the
latent trips, resulting in a relatively inconvenient network but
with lower objective. Once the {\sc ODMTS-DFD} problem is solved, the
evaluation procedure is carried out using the entire trip set $T$, and
some of the latent trips reject the discovered ODMTS although they
were involved in $\hat{T}$.  Another interesting observation from
Table~\ref{tab:atl_adoption} is at least 51.72\% percent of the newly
adopted riders can bring profit to the transit agency. It is worth
pointing out that only the shuttle operating cost and the \$2.5 revenue
from the transit ticket are involved in the computation.

\subsection{Travel Distance, Travel Time, and Operating Cost}
\label{subsubsect:atl_other_results}

\begin{table}[!t]
\centering%
\resizebox{0.9\textwidth}{!}{
\begin{tabular}{l r c c c c }%
\toprule
Algorithm&Trip Rule& Car Dist. after ODMTS (km)& Car Dist. Reduced (km)  &\% Car Dist. Reduced&Bus Dist. (km)\\%
\midrule%
Exact Alg.&{-}&434.43k&147.58k&25.36&8.50k\\%
\midrule%
Alg.~\mynameref{alg:rho_GRAD}&{-}&428.30k&153.71k&26.41&8.02k\\%
\midrule%
Alg.~\mynameref{alg:eta_GRRE}&{-}&423.62k&158.39k&27.21&7.35k\\%
\midrule%
Alg.~\mynameref{alg:rho_GAGR}&{-}&429.58k&152.44k&26.19&7.85k\\%
\midrule%
Alg.~\mynameref{alg:arc_s1}&\ref{rule:all_adopt_odmts}&429.28k&152.73k&26.24&8.30k\\%
\midrule%
Alg.~\mynameref{alg:arc_s1}&\ref{rule:upper_bound_odmts}&414.93k&167.08k&28.71&5.46k\\%
\midrule%
Alg.~\mynameref{alg:arc_s2}&\ref{rule:use_fixed_route_odmts}, \ref{rule:all_adopt_odmts}&421.89k&160.12k&27.51&7.62k\\%
\midrule%
Alg.~\mynameref{alg:arc_s2}&\ref{rule:upper_bound_odmts}, \ref{rule:all_adopt_odmts}&427.89k&154.13k&26.48&8.58k\\%
\bottomrule%
\end{tabular}%
}
\caption{Travel Distance of additional riders. The \textit{car distance before ODMTS} aggregates the travel distance covered by the additional riders when they drive to commute. On the other hand, the \textit{car distance after ODMTS} sums the travel distance on shuttles of the adopting riders and the
  travel distance of the non-adopting riders.}
\label{tab:atl_result_distance}
\end{table}

This section further analyzes the ODMTS design using other performance
measures. Table \ref{tab:atl_result_distance} summarizes the impact of
the ODMTS on travelled distance by car. Note that the core riders are
existing transit users; hence, they are not included in this
table. What stands out of the table is the substantial reduction in
travel distance after a group of self-driving commuters adopts the
ODMTS, which is replaced by the bus or rail services.  Table
\ref{tab:atl_result_distance} also reports the travel distance of the
ODMTS buses. As shown in Figure~\ref{subfig:atl_stops}, the existing bus system in Atlanta exhibits a dense cluster of stops, resulting in significant emissions and more congested traffic due to the regular pull in and pull out movements within traffic. This situation is further exacerbated by the utilization of diesel fuels. It's important to note that this illustration showcases clustered bus stops; the original stop network is even more extensive. In contrast, ODMTS only operate buses on a few key corridors that connect more widely spaced stops, offering a solution to alleviate congestion and emissions. The observations above together suggest that ODMTS can provide a positive
impact on morning peak traffic by reducing congestion and pollution.

\begin{table}[!t]
\centering%
\resizebox{0.7\textwidth}{!}{
\begin{tabular}{l r c c c c c c}%
\toprule
&&\multicolumn{2}{c}{Riders Adopting ODMTS}&\multicolumn{2}{c}{Core Riders}&\multicolumn{2}{c}{Riders not adopting ODMTS}\\%
\midrule%
Algorithm&Trip Rule&ODMTS&Direct&ODMTS&Direct&ODMTS&Direct\\%
\midrule%
Exact Alg.&{-}&18.86&15.64&24.00&20.94&38.96&20.17\\%
\midrule%
Alg.~\mynameref{alg:rho_GRAD}&{-}&18.91&15.55&24.03&20.94&38.46&20.15\\%
\midrule%
Alg.~\mynameref{alg:eta_GRRE}&{-}&18.62&15.44&23.99&20.94&39.50&20.83\\%
\midrule%
Alg.~\mynameref{alg:rho_GAGR}&{-}&18.84&15.52&24.04&20.94&38.51&20.19\\%
\midrule%
Alg.~\mynameref{alg:arc_s1}&\ref{rule:all_adopt_odmts}&18.79&15.47&24.01&20.94&38.33&20.21\\%
\midrule%
Alg.~\mynameref{alg:arc_s1}&\ref{rule:upper_bound_odmts}&19.06&16.14&23.83&20.94&39.68&21.25\\%
\midrule%
Alg.~\mynameref{alg:arc_s2}&\ref{rule:use_fixed_route_odmts}, \ref{rule:all_adopt_odmts}&19.16&15.73&23.87&20.94&37.99&19.98\\%
\midrule%
Alg.~\mynameref{alg:arc_s2}&\ref{rule:upper_bound_odmts}, \ref{rule:all_adopt_odmts}&18.97&15.55&24.02&20.94&37.81&20.05\\%
\bottomrule%
\end{tabular}%
}
\caption{Average Travel Time of Different Groups.}%
\label{table:atl_result_travel_time}
\end{table}

Table \ref{table:atl_result_travel_time} presents the ODMTS travel
time of three different groups: (i) riders with choice who adopt the
system, (ii) core riders, and (iii) riders with choice who continue to
drive alone. For each group, the travel times of two modes, ODMTS and
direct driving, are compared. For the riders who adopt ODMTS, the
average ODMTS travel time is only slightly greater than the average
direct driving time. This can be explained by recalling the direct
shuttle percentage reported in Table~\ref{tab:atl_adoption}. On the
other hand, the non-adopting commuters
experience a larger time gap, where the average duration of the routes
suggested by ODMTS is still less than two times of the average direct
trip duration. For core riders, ODMTS can offer a level of service nearly comparable to driving a car, despite their lack of personal vehicles. It is also important to compare the travel time of the
existing system for the core riders. As reported in
\citet{auad2021resiliency}, when evaluating on the same data (core
riders in this study) with the existing system, on average, the core
riders spend 40 minutes on the transit system if they are assumed to
only use the bus services, and this value drops to 26 minutes when
rail system is available. Thus, from the core riders' perspective, all
of the ODMTS systems found by the iterative algorithms outperform the
existing transit system. The advantages of ODMTS stem from its effective utilization of the rail system, where many shuttle legs serve as feeders to hubs, thereby circumventing congestion issues. In contrast, driving alone often entails contending with traffic congestion, especially during the morning peak of a workday.

\begin{table}%
\centering%
\resizebox{0.8\textwidth}{!}{
\begin{tabular}{l r r r r r r }%
\toprule
Algorithm&Trip Rule& Revenue (\$) &Bus Investment (\$) &Shuttle Opt. Cost (\$) &Rail Opt. Cost (\$) & NP/rider (\$)\\%
\midrule%
Exact Alg.&{-}&218.23k&13.64k&136.51k&24.87k&0.49\\%
\midrule%
Alg.~\mynameref{alg:rho_GRAD}&{-}&216.32k&12.75k&124.51k&24.87k&0.63\\%
\midrule%
Alg.~\mynameref{alg:eta_GRRE}&{-}&222.38k&10.97k&133.58k&24.87k&0.60\\%
\midrule%
Alg.~\mynameref{alg:rho_GAGR}&{-}&216.36k&12.26k&125.01k&24.87k&0.63\\%
\midrule%
Alg.~\mynameref{alg:arc_s1}&\ref{rule:all_adopt_odmts}&215.90k&12.98k&123.58k&24.87k&0.63\\%
\midrule%
Alg.~\mynameref{alg:arc_s1}&\ref{rule:upper_bound_odmts}&238.33k&7.59k&177.91k&24.87k&0.29\\%
\midrule%
Alg.~\mynameref{alg:arc_s2}&\ref{rule:use_fixed_route_odmts}, \ref{rule:all_adopt_odmts}&217.28k&12.36k&127.21k&24.87k&0.61\\%
\midrule%
Alg.~\mynameref{alg:arc_s2}&\ref{rule:upper_bound_odmts}, \ref{rule:all_adopt_odmts}&215.15k&13.90k&121.48k&24.87k&0.64\\%
\bottomrule%
\end{tabular}%
}
\caption{Revenue and Cost from the Agency Perspective. Opt. and NP stand for Operating and Net Profit.}%
\label{tab:atl_result_cost}
\end{table}

Finally, the cost analysis from the perspective of the transit agency
is shown in Table \ref{tab:atl_result_cost}. Ticket is the only source
of revenue considered in this table, and \textit{ticket revenue}
stands for the revenue collected from the core riders and the adopting
riders. The \textit{bus investment} and \textit{shuttle operating
  cost} are directly computed from the equations explained in Section
\ref{sect:original_problem}. This study assumes that the ODMTS is
operated by MARTA; thus, a fixed \textit{rail operating cost} between
6 am---10 am is also considered in this table. The rail operating cost
is estimated based on the MARTA rail schedule and the approximated
\$240 operating cost per revenue hour is reported in fiscal year 2018
(shown in Section 4.7.1 in \citet{marta2019report}). After including
the 24.87k rail operating cost, the \textit{net profit per rider} is
computed. Table \ref{tab:atl_result_cost} shows that the ODMTS system
is currently earning between \$0.29-\$0.64 when serving a customer,
depending on the selected design. This is a result of the on-demand shuttles and the rapid bus services efficiently replacing the traditional buses, especially the buses that operate in regions with low ridership. In addition, the current potential riders only reside in one specific area, i.e., East Atlanta. Given that the rail operating cost is fixed, a case study with an augmented additional demand region might increase the number of adopting riders, resulting in a citywide profitable system.
\end{APPENDICES}

\end{document}